\documentclass{article}
\usepackage{amsmath,amssymb,amsthm}
\usepackage{csquotes}
\usepackage{tikz}
\usepackage{tikz-cd}
\usepackage{geometry}
\newtheorem{theorem}{Theorem}[section]
\newtheorem{lemma}[theorem]{Lemma}
\newtheorem{corollary}[theorem]{Corollary}

\newtheorem{proposition}[theorem]{Proposition}
\theoremstyle{definition}
\newtheorem{definition}[theorem]{Definition}
\newtheorem{example}[theorem]{Example}
\theoremstyle{remark}
\newtheorem{remark}[theorem]{Remark}
\newtheorem{notation}[theorem]{Notation}
\usepackage{graphicx}

\makeatletter
\providecommand*{\Dashv}{%
	\mathrel{%
		\mathpalette\@Dashv\vDash
	}%
}
\newcommand*{\@Dashv}[2]{%
	\reflectbox{$\m@th#1#2$}%
}
\makeatother
\usepackage{color}
\usepackage{bussproofs}
\usepackage{comment}
\usepackage{ stmaryrd }
\usepackage{graphicx}
\usepackage{amsmath,amssymb,amsthm}
\usepackage{color}
\usepackage{csquotes}
\usepackage[scale=2]{ccicons}
\usepackage{tikz-cd}
\usepackage{tikz}
\usepackage{pgfplots}
\usepackage{bussproofs}
\usepackage{comment}
\usepackage{ stmaryrd }
\usepackage{ wasysym }
\usepackage{graphicx}
\usepackage{stmaryrd}
\usetikzlibrary{positioning, shapes}
\usetikzlibrary{decorations.text,calc,arrows.meta}
\usepgfplotslibrary{dateplot}
\usepackage{adjustbox}
\usepackage{txfonts}
\usepackage{pxfonts}
\usepackage{mathrsfs}  

\newcommand*{\twoheadrightarrowtail}{\mathrel{\rightarrowtail\kern-1.9ex\twoheadrightarrow}}
\newcommand{\counterfactual}{\ensuremath{%
		\Box\kern-1.5pt
		\raise1pt\hbox{$\mathord{\rightarrow}$}}}
\newcommand{\counterfactuals}{\ensuremath{%
		\lozenge\kern-1.5pt
		\raise1pt\hbox{$\mathord{\rightarrow}$}}}
\newcommand{\alg}[1]{{\textbf{\upshape #1}}}  %
\newcommand{\vv}[1]{\mathsf {#1}}

\newcommand{\ba}{\boxright}
\newcommand{\BA}{\boxRight}
\newcommand{\BAS}{\boxRight_S}
\newcommand{\form}{Fm_{\mathcal{L}}}
\newcommand{\Cl}{\mathbf{LV}}	
\newcommand{\Cg}{\mathbf{GV}}
\newcommand{\Cul}{\mathbf{LVC}}	
\newcommand{\Cug}{\mathbf{GVC}}
\newcommand{\Cdl}{\mathbf{LVCS}}	
\newcommand{\Cdg}{\mathbf{GVCS}}
\newcommand{\LCA}{\mathsf{VA}}
\newcommand{\LC}{\mathsf{VCA}}
\newcommand{\CA}{\mathsf{VCSA}}

\renewcommand{\d}{\delta}
\newcommand{\f}{\varphi}

\newcommand{\e}{\varepsilon}

\newcommand{\sse}{\subseteq}

\newcommand{\app}{\approx}

\newcommand{\cc}[1]{\mathcal{#1}}
\newcommand{\stone}{\mathfrak{s}}
\newcommand{\STONE}{\mathscr{S}}
\newcommand{\COUNT}{\mathscr{C}}
\newcommand{\ALPHA}{\mathscr{A}}
	
\newcommand{\CATA}{\mathsf{AM}}
\newcommand{\Sphere}{\mathcal{L}}	
\newcommand{\Alpha}{\mathcal{F}}
\newcommand{\TOPSPH}{\mathsf{SPH}}

\usetikzlibrary{positioning, shapes}
\usetikzlibrary{decorations.text,calc,arrows.meta}

\usepackage[backend=biber]{biblatex}

\title{The Algebras of Lewis's Counterfactuals}
\author{Giuliano Rosella and Sara Ugolini}

\begin{document}
	
	\maketitle
	\begin{abstract}
		
		The logico-algebraic study of Lewis's hierarchy of variably strict conditional logics has been essentially unexplored, hindering our understanding of their mathematical foundations, and the connections with other logical systems. This work aims to fill this gap by providing a comprehensive logico-algebraic analysis of Lewis's logics. We begin by introducing novel finite axiomatizations for varying strengths of Lewis's logics, distinguishing between global and local consequence relations on Lewisian sphere models. We then demonstrate that the global consequence relation is strongly algebraizable in terms of a specific class of Boolean algebras with a binary operator representing the counterfactual implication; in contrast, we show that the local consequence relation is generally not algebraizable, although it can be characterized as the degree-preserving logic over the same algebraic models. Further, we delve into the algebraic semantics of Lewis's logics, developing two dual equivalences with respect to particular topological spaces. In more details, we show a duality with respect to the topological version of Lewis's sphere models, and also with respect to Stone spaces with a selection function; using the latter, we demonstrate the strong completeness of Lewis's logics with respect to sphere models. Finally, we draw some considerations concerning the limit assumption over sphere models. 
	\end{abstract}
	
	\section{Introduction}
	
	A counterfactual conditional (or simply a counterfactual) is a conditional statement of the form ``If {\em antecedent} were the case, then {\em consequent} would be the case'', where the antecedent is usually assumed to be false. Counterfactuals have been studied in different fields, such as linguistics, artificial intelligence, and philosophy. The logical analysis of counterfactuals is rooted in the work of Lewis \cite{Lewis1973, Lewis71} and Stalnaker \cite{Stalnaker1968, STALNAKER2008} who have introduced what has become the standard semantics for counterfactual conditionals based on particular Kripke models equipped with a similarity relation among the possible worlds. In Lewis's language, a counterfactual is formalized as a formula of the kind $\varphi \ba \psi$ which is intended to mean that if $\varphi$ were the case, then $\psi$ would be the case. 
	Lewis develops a hierarchy of logics meant to deal with different kinds of counterfactual conditionals that have had a notable impact, and are quite well-known and studied; surprisingly, the literature on these logics (quite vast: Lewis's book counts thousands of citations at the present date) essentially lacks a logico-algebraic treatment. 
	
	It is important to stress that the role of algebra has been pivotal in the formalization and understanding of reasoning; indeed, modern logic really flourishes with the rise of the formal methods of mathematical logic, which moves its first steps with George Boole's intuition of using the symbolic language of algebra as a mean to formalize how sentences connect together via logical connectives \cite{Boole1854}. More recently, the advancements of the discipline of (abstract) algebraic logic have been one of the main drivers behind the surge of systems of {\em nonclassical logics} in the 20th century, in particular via the notion of {\em equivalent algebraic semantics} of a logic introduced by Blok and Pigozzi \cite{Block1989}. In this framework, the deductions of a logic and its properties are fully and faithfully interpreted by the semantical consequence of the related algebras.

	While a few works \cite{weiss, Nute1980-NUTTIC-5, Segerberg1989, Rosella2023-ROSCAM-10} present a semantics in terms of algebraic structures for Lewis's conditional logics, the results therein are either partial or fall outside the framework of the abstract algebraic analysis. On a different note, the proof-theoretic perspective on Lewis's conditional logics is instead more developed, in particular it is carried out in a series of recent works \cite{Girlando2016, Girlando2019, Negri2016-NEGPAF-2}. 
	However, although the research on Lewis's conditional logics has been and still is very prolific, a foundational work that carries Lewis's hierarchy within the realm of the well-developed discipline of (abstract) algebraic logic is notably missing in the literature; the present manuscript aims at filling this void. 
	
	To this end, we start by considering Lewis's logics as consequence relations, instead of just sets of theorems; this brings us to consider two different kinds of derivation, depending on whether the deductive rules are applied only to theorems (giving a relatively weaker calculus) or to all derivations (i.e. yielding a stronger calculus). We stress that this distinction, although relevant, is often blurred in the literature.
	As it is the case for modal logic (see \cite{Blackburn2001, Wen2021}), these two choices turn out to correspond to considering two different consequence relations on the intended sphere models: a \emph{global} and a \emph{local} one. The finite strong completeness of the local consequence with respect to the weaker calculi follows by Lewis's work; at the end of this paper we will extend the latter to strong completeness (i.e., accounting for deductions from infinite theories), also providing the analogous result for the global consequence relation with respect to the stronger calculi. 
	
	The parallel with modal logic will guide the groundwork for this investigation; following some results connecting modal logic and Lewis counterfactuals (see \cite{Lewis1973, Williamson2010-WILMLW}) we show that one can term-define in the language a modal operator $\Box$ that can be used to show a deduction theorem for the strong calculus (whereas the weak calculus is known to have the usual deduction theorem).
	The reader shall notice that the binary operator $\boxright$ does not straightforwardly inherit the plethora of results on modal operators; for the reader more expert on the algebraic perspective, the models are not Boolean algebras with an operator in the usual sense (see \cite{JonssonTarski}), since $\ba$ is not additive on both arguments (more precisely, it only distributes over infima on the right) and it cannot be recovered from a unary modal operator. 
	
	Our next main results show that the stronger calculi, associated to the global consequence relation, are strongly algebraizable in the sense of Blok-Pigozzi, i.e. there is a class of algebras (a variety of Boolean algebras with an extra binary operator $\ba$) that serve as the equivalent algebraic semantics; the weaker calculi, associated to the local consequence relation, are shown to not be algebraizable (there is no class of algebras whose consequence relation ``corresponds'' to the deduction of the logic), but they are the logics preserving the degree of truth of the same algebraic models. Thus the same class of algebras can be meaningfully used to study both versions of Lewis's logics.
	We also initiate the study of the structure theory of the algebraic models, showing in particular that the congruences of the algebras can be characterized by means of the congruences of their modal reducts.
	
	In the last sections, we circle back to the original intended models; it will become clear that to properly consider infinite models one should not simply consider sets, but {\em topological spaces}; i.e., the subsets of the universe that are meant to represent the formulas (the clopen sets of the topology) play a special role. In more details, we show two different dual categorical equivalences of our algebraic structures with respect to topological spaces based respectively on Lewis's spheres and (Stalnaker's inspired) selection function. The dualities we show are enrichments of Stone duality between Boolean algebras with homomorphisms and Stone spaces with continuous maps, where the operator $\ba$ is interpreted first by means of a selection function, and then by a map associating a set of nested spheres to each element of the space. We note that the constructions used to move between the different structures are based on those introduced by Lewis in \cite{Lewis71}; we here then provide also detailed proofs of some facts that have been used without proof in the literature. 
	The formal work developed for the dualities will allow us to demonstrate the strong completeness of sphere models with respect to Lewis's logics. Finally, we will also clarify the role of the {\em limit assumption}, a condition on sphere models that has been extensively discussed in the literature. In particular, we will see that both the strong and weak calculi are strongly complete with respect to models that do satisfy the limit assumption; in this sense, models without the limit assumption are not really ``seen'' by Lewis's logics.
	
	\section{Logico-algebraic preliminaries}
	This work aims to provide a logico-algebraic treatment of Lewis's logics. In order to make the results more significant to the reader, in this section we recall the basics of the Blok-Pigozzi machinery \cite{Block1989}, that connects {\em algebraizable logics} and particular classes of algebras ({\em quasivarieties of logic}). A quasivariety is a class of algebras that can be axiomatized by means of quasiequations (formulas where a conjunction of equations imply one single equation, intuitively playing the role of the rules of the corresponding logic). The underlying idea is indeed that the deducibility of the logic is fully and faithfully interpreted in the equational consequence relation of the quasivariety. In particular, we show that the corresponding algebraic models to Lewis's logics are {\em varieties}, i.e. equations suffice for the axiomatization of the class of algebras. 
	For the omitted details we refer the reader to \cite{Block1989, Fontbook}. 

Let us set some notation; given an algebraic language $\rho,$ and a denumerable set of variables $Var$, we write $Fm_{\rho}(Var)$ for the {\em algebra of formulas} written over the language $\rho$ and variables in $Var$. An {\em equation} is a pair $p,q$ of $\rho$-terms (i.e. elements of $Fm_{\rho}(Var)$) that we write suggestively as $p \app q$.
An {\em assignment} of $Var$ into an algebra $\alg A$ of type $\rho$ is a function $h$ mapping each variable to an element of $\alg A$, that extends (uniquely) to a function that respects all the operations (i.e., an homomorphism) from the algebra of formulas $Fm_{\rho}(Var)$ to $\alg A$.
An algebra $\alg A$ satisfies an equation $p \app q$ with an assignment $h$ (and we write $\alg A, h \vDash p \approx q$) if $h(p) = h(q)$ in $\alg A$.  An equation $p \app q$ is {\em valid} in $\alg A$ (and we write $\alg A \vDash p \approx q$) if for all assignments $h$ in $\alg A$, $\alg A, h \vDash p \approx q$; if $\Sigma$ is a set of equations then
 $\alg A \vDash \Sigma$ if $\alg A \vDash \sigma$ for all $\sigma \in \Sigma$. An equation  $p \app q$ is valid in a class of algebras $\vv K$ if $\alg A \vDash p \approx q$ for all $\alg A \in \vv K$.
 
Let us now be more precise about what we consider a logic.	
A {\em consequence relation} on the set of terms $Fm_{\rho}(Var)$ of some algebraic language $\rho$ is a relation $\vdash\, \sse \cc P(Fm_{\rho}(Var)) \times Fm_{\rho}(Var)$ (and we write $\Sigma \vdash \gamma$ for $(\Sigma, \gamma) \in \,\vdash$)  such that:
\begin{enumerate}
	\item if $\alpha \in \Gamma$ then $\Gamma \vdash \alpha$;
	\item if $\Gamma \vdash \delta$ for all $\delta \in \Delta$ and $\Delta \vdash \beta$, then $\Gamma \vdash \beta$.
\end{enumerate}
We call {\em substitution} any endomorphism of $Fm_{\rho}(Var)$ (i.e., any function on itself that respects all operations); $\vdash$ is {\em substitution invariant} (also called {\em structural}) if $\Gamma \vdash \alpha$ implies $\{\sigma(\gamma): \gamma \in \Gamma\} \vdash \sigma(\alpha)$ for each substitution $\sigma$. Finally, $\vdash$ is {\em finitary} if $\Gamma \vdash \alpha$ implies that there is a finite $\Gamma' \sse \Gamma$ such that $\Gamma' \vdash \alpha$. 
By a \emph{logic} $\cc L$ in what follows we mean a substitution-invariant finitary consequence relation $\vdash_{\cc L}$ on  $Fm_{\rho}(Var)$ for some algebraic language $\rho$, $\vdash_{\cc L}\, \sse \cc P(Fm_{\rho}(Var)) \times Fm_{\rho}(Var)$.

	In loose terms, to establish the algebraizability of a logic $\cc L$ with respect to a quasivariety of algebras $\vv Q_{\cc L}$ over the same language $\rho$, one needs a finite set of one-variable equations $$\tau(x) =\{\d_i(x) \app \e_i(x): i = 1,\dots,n\}$$ over terms of type $\rho$ and a finite set of formulas of $\cc L$
	in two variables $$\Delta(x,y)=\{\f_1(x,y),\dots,\f_m(x,y)\}$$ that allow to transform  equations in $\vv Q_{\cc L}$ into formulas of $\cc L$ and vice versa; moreover this transformation must respect both the consequence relation of the logic and the semantical consequence of the quasivariety. More precisely, for all sets of formulas $\Gamma$ of $\cc L$ and formulas $\f \in Fm_{\rho}(Var)$
	\begin{equation}\label{eq:alg1}
		\Gamma \vdash_{\cc L} \f\quad\text{iff}\quad \tau(\Gamma) \vDash_{\vv Q_{\cc L}} \tau(\f)
	\end{equation}
	where $\tau(\Gamma)$ is a shorthand for $\{\tau(\gamma): \gamma \in \Gamma\}$, and also
	\begin{equation}\label{eq:alg2}
		(x \app y) \Dashv \vDash_{\vv Q_{\cc L}}\tau(\Delta(x,y)).
	\end{equation}
	
	A quasivariety $\vv Q$ is a {\em quasivariety of logic} if it is the equivalent algebraic semantics for some logic $\cc L_\vv Q$.
	\begin{example}
		Classical logic is algebraizable with respect to the variety of Boolean algebras, as testified by $\tau(x) = x \approx 1$ and $\Delta(x,y) = x \leftrightarrow y$.
	\end{example}
	While the two conditions (\ref{eq:alg1}) and (\ref{eq:alg2}) above are necessary and sufficient to show algebraizability of a logic with respect to a quasivariety, in some cases there are easier ways to check that a logic is algebraizable. In fact, many of the well-known algebraizable logics belong to the class of {\em implicative logics}, that is, they have a well-behaved binary connective $\to$ which allows to show that (\ref{eq:alg1}) and (\ref{eq:alg2}) hold.
	\begin{definition}\label{def:implicativelogic}
		An {\em implicative logic} is a logic $\vdash$ with a binary term $\to$ such that:
		\begin{enumerate}
			\item $\vdash x \to x$
			\item $x \to y, y \to z \vdash x \to z$
			\item\label{item:implic} $x_1 \leftrightarrow y_1, \ldots, x_n \leftrightarrow y_n \vdash \lambda(x_1, \ldots, x_n) \leftrightarrow \lambda(y_1, \ldots, y_n)$ for each term $\lambda \in \mathcal{L}$ of arity $n > 0$
			\item $x, x \to y \vdash y$
			\item $x \vdash y \to x$.
		\end{enumerate}
	\end{definition}
	In an implicative logic that does not have a constant $1$ that is a theorem, one  can always define $1:= x \to x$ for a fixed variable $x$, and $1$ is a theorem by the above definition.
	\begin{theorem}[\cite{Fontbook}]
		All implicative logics are algebraizable, with defining equation
		$\tau(x) := \{x \approx x \to x\}$ and equivalence formulas $\Delta(x,y) := \{x\to y, y\to x\}$. The quasivariety can be presented by the equations and quasi-equations that result by applying the transformation $\tau$ to the axioms and rules of any Hilbert-style presentation of the logic.
	\end{theorem}
	Classical logic is an example of an implicative logic. The class of algebras that we will be interested in are extensions of Boolean algebras with an extra operator $\ba$ meant to interpret the counterfactual implication. 
	
	From the more strictly algebraic perspective, we mention for the interested reader that in a quasivariety that is the equivalent algebraic semantics of a logic $\cc L$ over a language $\rho$, congruences are in one-one correspondence with the {\em deductive filters} induced by the logic. A deductive filter $F$ of an algebra $\alg A$ is a subset of the domain of $\alg A$ that is closed under the interpretation of the axioms and rules of the logic; that is, consider a Hilbert style presentation of the logic, then for every rule $\Gamma \vdash \varphi$ and every homomorphism $f$ from $Fm_{\rho}(Var)$ to $\alg A$, if $f[\Gamma] \sse F$, then $f(\varphi) \in F$. 
	\begin{theorem}[{\cite[Theorem 3.51]{Fontbook}}]\label{thm:isomorphismfilters}
		Let $\vdash$ be a logic with equivalent algebraic semantics $\vv K$, and let $\alg A$ be an algebra of the same type of $\vv K$. Then the $\vdash$- deductive filters of $\alg A$ are in bijection with the $\vv K$-relative congruences of $\alg A$.
	\end{theorem}
	We mention that the deductive filters of the algebras of formulas are the {\em theories} of the logic. In particular, if $\cc L$ is an implicative logic, for every $\alg A$ in $\vv Q_{\cc L}$ the correspondence between congruences and deductive filters is given by the following maps:
	$$\theta \mapsto F_\theta = \{a \in A: (a, a \to a) \in \theta\}, \qquad F \mapsto \theta_F =\{(a, b): a \to b, b \to a \in F\}.$$
	
	where $\theta$ is any congruence of $\alg A$ and $F$ is any deductive filter of $\alg A$. If there is a constant $1$ in the language of $\cc L$ that is a theorem, as it is the case for classical logic, then congruences are totally determined by their $1$-blocks, i.e. the first map above becomes:
	$$\theta \mapsto F_\theta = \{a \in A: (a, 1) \in \theta\}.$$
	Precisely, $\vv Q_{\cc L}$ is said to be {\em 1-regular}, or {\em ideal-determined} with respect to $1$ (see \cite{AglianoUrsini,GummUrsini}).

	\section{Lewis's Logics: axioms and sphere models}
	This section lays out the groundwork for a logico-algebraic study of the hierarchy of logics introduced by Lewis in his seminal book \cite{Lewis1973} to reason with counterfactuals conditionals, and their intended models, i.e. sphere models. All the logics in the hierarchy are axiomatic extensions of the system $\mathbf{V}$, which according to Lewis is the ``weakest system that has any claim to be called the logic of conditionals'' \cite[p. 80]{Lewis71}; therefore our investigation starts with $\mathbf{V}$, and all our results will carry through its axiomatic extensions\footnote{The system $\mathbf{V}$ from \cite{Lewis1973} is equivalent to the system $\mathbf{C0}$ in \cite{Lewis71}, essentially differing in minor differences in the language.}.
	
	In particular, we will give a new and simpler axiomatization of $\mathbf{V}$ with respect to the original ones (\cite{Lewis1973}); we will take the counterfactual connective as a primitive symbol in the language, and we will distinguish between two different consequence relations: a weak one, where the rules of the calculus only apply to theorems (which is the one usually considered in the literature), and a strong one, where the rules instead apply to all deductions. We will see that these two choices correspond to considering two different consequence relations over sphere models: a local and a global one, in complete analogy with the case of modal logic. This analogy will continue and guide the investigation throughout the rest of the section. In particular, we will use the tool of generated submodels, borrowed from modal logic (see \cite{Blackburn2001}), and apply it to sphere models to first characterize the global consequence relation via the local one, secondly to prove a deduction theorem, and finally to prove the  finite strong completeness of the global consequence with respect to the strong version of the presented Hilbert-style calculus. 
	
	\subsection{Axiomatizations}	
	We fix the language $\mathcal{L}$ obtained from a denumerable set of variables $Var$, and expanding the language of propositional classical logic $\{\land, \lor, \to, 0, 1\}$ with a binary connective $\ba$, where $\varphi \ba \psi$ should be read as \begin{center}``if it were the case that $\varphi$, then it would be the case that $\psi$''.  \end{center} 
	As usual, one can define further connectives by:  $$\neg x : = x \to 0, \;\; x \leftrightarrow y : = (x \to y) \land (y \to x).$$ The following derivative connectives will be also considered: \begin{center} $(x \Diamondright y):=\neg(x \boxright \neg y)$ \\ $ \square x := \neg x \boxright x$ \\ $\Diamond x:=\neg\square\neg x$ \\  $x \preccurlyeq y := ((x \vee y)\Diamondright(x \vee y)) \to ((x\vee y)\Diamondright x)$ \\ $x \prec y := \neg (y \preccurlyeq x)$ \\ $x \approx y := (x \preccurlyeq y) \wedge (y \preccurlyeq x)$ \end{center}

	Let us denote with $\form$ the set of all $\mathcal{L}$-formulas over variables in $Var$.

	While often in the literature $\mathbf{V}$ is presented as a set of theorems, we are interested in studying logics as consequence relations.
	We will hence distinguish two different logics, $\mathbf{GV}$ and $\mathbf{LV}$, which arise depending on whether we apply the rules of Lewis calculus only to theorems (for the weaker logic $\mathbf{LV}$) or to all derivations (for the stronger logic $\mathbf{GV}$). We remark that this distinction, although significant, is often blurred in the literature.
	
	The two systems $\Cg$ and $\Cl$ share the same axioms, that is, given $\varphi, \psi, \gamma \in \form$ we have:
	\begin{itemize}
		\item[(L0)] the reader's favorite Hilbert-style axioms of classical logic \footnote{The reader can find some in \cite{ukasiewicz1963-UKAEOM}.};
		\item[(L1)] $\vdash\varphi \boxright \varphi$;
		\item[(L2)] $\vdash((\varphi \boxright \psi) \wedge (\psi \boxright \varphi)) \to  ((\varphi \boxright \gamma)\leftrightarrow(\psi \boxright \gamma))$;
		\item[(L3)] $\vdash((\varphi \vee \psi)\boxright \varphi) \vee ((\varphi \vee \psi)\boxright \psi) \vee (((\varphi \vee \psi)\boxright \gamma);\leftrightarrow((\varphi \boxright \gamma)\wedge(\psi \boxright \gamma))$;
		\item[(L4)] $\vdash(\varphi \ba (\psi \land \gamma)) \leftrightarrow ((\varphi \ba \psi) \land (\varphi \ba \gamma))$;
	\end{itemize}
	Moreover, both $\Cg$ and $\Cl$ satisfy modus ponens: 
	\begin{itemize}
		\item[(MP)] $\varphi, \varphi \to \psi\vdash \psi$.
	\end{itemize}
	While $\Cg$ satisfies the following rule involving the counterfactual implication:
	\begin{itemize}
		\item[(C)] $\varphi \to \psi \vdash (\gamma \ba \varphi) \to (\gamma \ba \psi)$,
	\end{itemize}
	$\Cl$ satisfies the following weaker version of the rule:
	\begin{itemize}
		\item[(wC)] $ \vdash \varphi \to \psi$ implies  $\vdash (\gamma \ba \varphi) \to (\gamma \ba \psi)$.
	\end{itemize}
	\begin{definition}
		The logic $\Cg$, denoted by $\vdash_{\Cg}$,  is the smallest finitary consequence relation satisfying all axioms [(L1)]--[(L4)], (MP), and (C). 
		
		The logic $\Cl$, denoted by $\vdash_{\Cl}$, is the smallest finitary consequence relation satisfying all axioms [(L1)]--[(L4)], (MP), and (wC).
	\end{definition}
	
    Lewis conditional logics consists in axiomatic extensions of the system $\mathbf{V}$ with the following axioms:
	
	\begin{itemize}

		\item[($\mathbf{W}$)] $\vdash(\varphi\boxright \psi)\to(\varphi\to\psi);$ (weak centering)
		\item[($\mathbf{C}$)] $\vdash ((\varphi\boxright \psi)\to(\varphi\to\psi))\wedge(\varphi \wedge \psi)\to(\varphi\boxright\psi)$ (centering)
		\item[($\mathbf{N}$)] $\vdash \square \varphi \to \Diamond \varphi$ (normality)
		\item[($\mathbf{T}$)] $\vdash \square \varphi \to \varphi$ (total reflexivity)
		\item[($\mathbf{S}$)] $\vdash(\varphi \ba \psi) \lor (\varphi \ba \neg \psi)$ (Stalnakerian)
		\item[($\mathbf{U}$)] $\vdash (\Diamond\varphi \square \to \square\Diamond\varphi) \wedge (\square \varphi \to \square\square \varphi)$ (uniformity)
		\item[($\mathbf{A}$)] $\vdash(\varphi\preccurlyeq \psi)\to\square(\varphi\preccurlyeq\psi) \wedge ((\varphi\prec\psi)\to\square(\varphi\prec\psi))$ (absoluteness)
		\end{itemize}
	
	We indicate a certain system in the family of Lewisian conditional logics by just juxtaposing to $\mathbf{GV}$ or $\mathbf{LV}$ the corresponding letter for axioms. For instance, $\mathbf{LVCA}$ indicates the axiomatic extension of the logic $\mathbf{LV}$ with the axioms $\mathbf{C}$ and $\mathbf{A}$. 

	Among these axiomatic extensions, it is worth mentioning the system $\mathbf{LVC}$ which is the considered by Lewis the ``correct logic of counterfactual conditionals'' \cite{Lewis71}, while $\mathbf{LVCS}$ essentially corresponds to Stalnaker's logic of conditionals \cite{Stalnaker1968, STALNAKER2008}.
	
	It is clear from the definition that $\Cg$ is a stronger deductive system than $\Cl$, i.e.:
	\begin{lemma}\label{lemma:Cgstronger}
		For any set of $\cc L$-formulas $\Gamma$ and $\cc L$-formula $\varphi$, $\Gamma \vdash_{\Cl} \varphi$ implies $\Gamma \vdash_{\Cg} \varphi$.
	\end{lemma}
	While $\Cg$ is strictly stronger than $\Cl$, e.g. the latter does not satisfy $(C)$, they do have the same theorems; the following proof is standard.
	\begin{theorem}\label{thm:samethm}
		$\Cg$ and $\Cl$ have the same theorems.
	\end{theorem}
	\begin{proof}
		It follows from Lemma \ref{lemma:Cgstronger} that if a formula $\varphi$ is a theorem of $\Cl$, it is also a theorem of $\Cg$; we show the converse.
		
		Let $\varphi$ be a theorem of $\Cg$, we show by induction on the length of the proof that $\varphi$ is a theorem of $\Cl$. The base case is for $\varphi$ being an axiom, thus the thesis holds given that $\Cg$ and $\Cl$ share the same axioms. Assume now that $\varphi$ is obtained by an application of a rule of $\Cg$, i.e., either by modus ponens or (DWC). But such rule is applied to theorems or axioms of $\Cg$, that by inductive hypothesis are theorems of $\Cl$; therefore, one can obtain the same conclusion by applying (MP) or (wDWC).
	\end{proof}
	We will now see that the axiomatization we have given is equivalent to the one given by Lewis in \cite{Lewis71}; with respect to the latter,
	we have added axioms [(L4)] and removed the denumerable set of rules describing ``deductions within conditionals''. Let us present the latter in the two versions, the strong ones:
	\begin{align}
		\psi \vdash \varphi \ba \psi,\tag{DWC$_0$}\\
		(\varphi_1 \wedge\dots\wedge\varphi_n)\to\psi\vdash ((\gamma \boxright\varphi_1)\wedge\dots\wedge(\gamma\boxright\varphi_n))\to(\gamma\boxright\psi)\; \tag{DWC$_n$}
	\end{align}
	for each $n \in \mathbb{N}, n \geq 1$, and the weaker versions: 
	\begin{align}
		\vdash \psi \mbox{ implies } \vdash \varphi \ba \psi,\tag{wDWC$_0$}\\
		\vdash (\varphi_1 \wedge\dots\wedge\varphi_n)\to\psi \mbox{ implies } \vdash((\gamma \boxright\varphi_1)\wedge\dots\wedge(\gamma\boxright\varphi_n))\to(\gamma\boxright\psi)\;\tag{wDWC$_n$}
	\end{align}
	for each $n \in \mathbb{N}, n \geq 1$.
	We start by noting that the monotonicity of $\ba$ on the consequent can be shown to be a consequence of the axioms.
	\begin{lemma}\label{lemma:vdashs}
		The following holds for all $ \cc L$-formulas $\varphi, \psi, \gamma$:
		\begin{enumerate}
			\item\label{lemma:vdashs-4}$\vdash_\mathbf{GV} (\varphi \ba \psi) \to (\varphi \ba (\psi \lor \gamma))$.
		\end{enumerate}
	\end{lemma}
	\begin{proof}
		Observe first that by the axioms and rules of classical logic, it holds that $ \vdash_\mathbf{GV} \psi \leftrightarrow (\psi \land (\psi \lor \gamma))$.
		Therefore, by using (C) and (L4) we get $$\vdash_\mathbf{GV} (\varphi \ba \psi) \leftrightarrow (\varphi \ba (\psi \land (\psi \lor \gamma))) \vdash_\mathbf{GV} (\varphi \ba \psi) \leftrightarrow ((\varphi \ba \psi) \land (\varphi \ba (\psi \lor \gamma))) $$
		from which we can derive $\vdash_\mathbf{GV} (\varphi \ba \psi) \to (\varphi \ba (\psi \lor \gamma))$, which concludes the proof.
	\end{proof}
	Moreover:
	\begin{lemma}\label{lemma:vdashs2}
		Consider a logical system $\vdash'$ in the language $\mathcal{L}$ satisfying the axioms of classical logic, (MP), and (DWC$_2$). Then:
		\begin{enumerate}
			\item $\varphi \to \psi \vdash' (\gamma \ba \varphi) \to (\gamma \ba \psi)$;
			\item $\varphi \leftrightarrow \psi \vdash' (\gamma \ba \varphi) \leftrightarrow (\gamma \ba \psi)$;
			\item $\vdash' (\varphi \ba \psi) \to (\varphi \ba (\psi \lor \gamma))$.
		\end{enumerate}
	\end{lemma}
	\begin{proof}
		For (1), from (DWC$_2$) it is easily shown that classically equivalent formulas can be substituted in the consequent of $\ba$; indeed we get:
		$$(\varphi \land \varphi) \to \psi \vdash' ((\gamma \ba \varphi) \land (\gamma \ba \varphi)) \to (\gamma \ba \psi).$$
		(2) is a direct consequence of (1).	For (3) we have the following:
		$$\vdash'(\psi \land \psi) \to (\psi \lor \gamma) \vdash' [(\varphi \ba \psi) \land (\varphi \ba \psi)] \to (\varphi \ba (\psi \lor \gamma)),$$
		which, given that $\vdash' \psi \leftrightarrow (\psi \land \psi)$, proves the claim via modus ponens.
	\end{proof}
	\begin{theorem}\label{thm:firstax}
		Consider a logical system $\vdash'$ in the language $\mathcal{L}$ satisfying the axioms of classical logic and (MP). The following are equivalent.
		\begin{enumerate}
			\item $\vdash'$ satisfies \emph{(L1)}--\emph{(L4)} and \emph{(C)};
			\item $\vdash'$ satisfies \emph{(L1)}--\emph{(L3)} and the rule \emph{(DWC$_2$)}.
		\end{enumerate}
		The same holds replacing \emph{(C)}, and \emph{(DWC$_2$)} with their weaker versions \emph{(wC)}, and \emph{(wDWC$_2$)}.
	\end{theorem}
	\begin{proof}
		Let us first show that (1) implies (2), i.e., we derive the rule (DWC$_2$) using (L4) and (C). 
		Using (C) we obtain that 
		$$(\varphi_1 \land \varphi_2) \to \psi \vdash' [\gamma \ba (\varphi_1 \land \varphi_2)] \to (\gamma \ba \psi);$$
		using now (L4) we get
		$$(\varphi_1 \land \varphi_2) \to \psi \vdash' [(\gamma \ba \varphi_1) \land (\gamma \ba \varphi_2)] \to (\gamma \ba \psi)$$
		which is exactly (DWC$_2$).
		
		Conversely, let us prove that (DWC$_2$) implies (C) and (L4); (C) can be derived by (DWC$_2$) by setting $\varphi_1 = \varphi_2 := \varphi$; for (L4), 
		note that by the axioms and rules of classical logic and (DWC$_2$) we get:
		\begin{align*}
			&\vdash' (\psi \land \gamma) \to (\psi \land \gamma) \;\vdash'\; ((\varphi \ba \psi)\land (\varphi \ba \gamma)) \to (\varphi \ba (\psi \land \gamma));
		\end{align*}
		moreover, one can use monotonicity of $\ba$ (Lemma \ref{lemma:vdashs}) and the fact that classically equivalent formulas can be substituted in the consequent of $\ba$ (Lemma \ref{lemma:vdashs2}) and obtain
		\begin{align*}
			&\vdash' (\varphi \ba (\psi \land \gamma)) \to  (\varphi \ba \psi),\\
			&\vdash' (\varphi \ba (\psi \land \gamma)) \to  (\varphi \ba \gamma),
		\end{align*}
		and therefore $\vdash' (\varphi \ba (\psi \land \gamma)) \to  ((\varphi \ba \psi) \land (\varphi \ba \gamma))$, which shows (L4). 
		
		The proof for the statement involving the weaker rules goes similarly, with the proviso that every derivation starts from theorems and axioms.
	\end{proof}
	The next theorem shows that both the alternative axiomatizations we have introduced are equivalent to the one presented by Lewis's in \cite{Lewis71}.
	\begin{theorem}\label{thm:secondax}
		Consider a logical system $\vdash'$ in the language $\mathcal{L}$ satisfying the axioms of classical logic and (MP). The following are equivalent.
		\begin{enumerate}
			\item $\vdash'$ satisfies \emph{(L1)}--\emph{(L3)} and \emph{(DWC$_n$)} for all $n \in \mathbb{N}$;
			\item $\vdash'$ satisfies 	\emph{(L1)}--\emph{(L3)} and \emph{(DWC$_2$)};
			\item $\vdash'$ coincides with $\vdash_\mathbf{GV}$.
		\end{enumerate}
		The same holds replacing \emph{(DWC$_2$)}, \emph{(DWC$_n$)}, and $\vdash_\mathbf{GV}$ with their weakened versions  \emph{(wDWC$_2$)}, \emph{(wDWC$_n$)}, and $\vdash_\mathbf{LV}$.
	\end{theorem}
	\begin{proof}
		(2) and (3) are equivalent by Theorem \ref{thm:firstax}; moreover, it is obvious that (1) implies (2); let us show the converse. Consider $n \in \mathbb{N}, n \geq 1$, then with (DWC$_2$) we obtain immediately:
		$$ ((\varphi_1 \wedge\dots\wedge \varphi_{n-1}) \wedge \varphi_n)\to\psi\vdash [((\gamma \boxright (\varphi_1\ldots\varphi_{n-1})) \land (\gamma \ba \varphi_n )]\to(\gamma\boxright\psi),$$
		which using the fact that (DWC$_2$) implies (L4) by Theorem \ref{thm:firstax}, yields that (DWC$_n$) holds for all $n \geq 1$. 
		In order to show that (DWC$_0$) holds, we first observe that $\varphi \ba 1$ is a theorem, indeed by (L4) we get
		$$\vdash'(\varphi \ba (\varphi \land 1)) \leftrightarrow ((\varphi \ba \varphi) \land (\varphi \ba 1)) \vdash' (\varphi \ba \varphi) \leftrightarrow ((\varphi \ba \varphi) \land (\varphi \ba 1)) \vdash' 1 \leftrightarrow (\varphi \ba 1)$$
		where in the derivations we used (L1) and substitution of classically equivalent formulas in the consequent of $\ba$ (Lemma \ref{lemma:vdashs2}, which can be used since (2) and (3) are equivalent).
		Finally, (DWC$_0$) is then a consequence of applying (DWC$_2$) with $\varphi_1 = \varphi_2 := 1$, $\psi := \psi$, $\gamma := \varphi$.
		
		The proof for the weaker calculus is completely analogous.
	\end{proof}	
	Lastly, let us observe that both $\Cg$ and $\Cl$ satisfy the substitution of logical equivalents, in the following sense.
	\begin{lemma}\label{lemma:equiv}
		The following hold for any $\varphi, \psi, \gamma \in Fm_{\cc L}$:
		\begin{enumerate}
			\item $(\varphi \leftrightarrow \psi) \vdash_\mathbf{GV} (\gamma \ba \varphi) \leftrightarrow (\gamma \ba \psi)$,
			\item $(\varphi \leftrightarrow \psi) \vdash_\mathbf{GV} (\varphi \ba \gamma) \leftrightarrow (\psi \ba \gamma)$,
			\item $\vdash_\mathbf{LV} (\varphi \leftrightarrow \psi) \mbox{ implies } \vdash_\mathbf{LV} (\gamma \ba \varphi) \leftrightarrow (\gamma \ba \psi)$ and $\vdash_\mathbf{LV} (\varphi \ba \gamma) \leftrightarrow (\psi \ba \gamma)$.
		\end{enumerate}
	\end{lemma}
	\begin{proof}
		First, notice that (1) follows by (C). Moreover, since (DWC$_1$) holds by Theorem \ref{thm:secondax}, one gets 
		 $$\varphi \to \psi \vdash_\mathbf{GV} (\varphi \ba \varphi) \to (\varphi \ba \psi) \vdash_\mathbf{GV} \varphi \ba \psi$$
		and its symmetric copy;
		thus it follows $\varphi \leftrightarrow \psi \vdash_\mathbf{GV} (\varphi \ba \psi) \land (\psi \ba \varphi),$
		which via (L2) gives (2). Finally, (3) follows from the previous points, given that $\Cg$ and $\Cl$ have the same theorems (Theorem \ref{thm:samethm}).
	\end{proof}

	\subsection{Sphere models, local and global semantics}
	Lewis bases his interpretation of the counterfactual connective $\ba$ on a neighbourhood-style semantics.
	The intuitive idea to evaluate the connective $\ba$ is that $\varphi \ba \psi$ is true at some world $w$ if in the closest worlds to $w$ in which $\varphi$ is true, also $\psi$ is true. This results in the definition of what Lewis calls a \enquote{variably strict conditional}, where the word \enquote{variably} stresses the fact that to evaluate counterfactuals with different antecedents at some world $w$, one might need to evaluate the corresponding classical implication in different worlds; from another point of view, this also means that in general a counterfactual $\varphi \ba \psi$ does not arise as some $\Box(\varphi \to \psi)$, for some modality $\Box$.
	Lewis formalizes this intuition by assigning to each possible world $w$ a nested set of {\em spheres}, which are subsets of possible worlds, meant to describe a similarity relationship with $w$; the smaller is the sphere to which a world $w'$ belongs, the closer, and therefore more similar, it is to $w$. Let us now be more precise.
	
	\begin{definition}\label{def:lewissphere} A sphere model $\mathcal{M}$ is a tuple $\mathcal{M}=\langle W, \mathcal{S}, v \rangle$ where:
		\begin{enumerate}
			\item $W$ is a set of possible worlds;
			\item $\mathcal{S}: W\to\mathcal{P}(\mathcal{P}(W))$ is a function assigning to each $w \in W$ a set $S(w)$ of nonempty nested subsets of $W$, i.e., for all $X, Y \in S(w)$, either $\emptyset \neq X\subseteq Y$ or $\emptyset \neq Y\subseteq X$.
			\item $v: Var \to \mathcal{P}(W)$ is an assignment of the variables to subsets of $W$, extending to all $\cc L$-formulas as follows:
				\[\begin{array}{lcl}
					v(0) & = & \emptyset \\
					v(1) & = & W \\
					v(\varphi \wedge \psi) & = & v(\varphi)\cap v(\psi)\\
					v(\varphi \vee \psi) & = & v(\varphi)\cup v(\psi)\\
					v(\varphi \to \psi) & = & (W\setminus v(\varphi)) \cup v(\psi)\\
					v(\varphi\boxright\psi) & = &\{w \in W : (\bigcup S(w)\cap v(\varphi)) =\emptyset\text{ or } \exists S \in S(w) \mbox{ such that } \emptyset \neq (S\cap v(\varphi)) \subseteq v(\psi))\}
				\end{array}\] 
				%
			\end{enumerate}		
		\end{definition}	
		Given a sphere model $\mathcal{M}=\langle W, \mathcal{S}, v \rangle$, and a set of $\cc L$-formulas $\Gamma$, we set:
		\begin{align}
			w \Vdash \Gamma \mbox{ iff } w \in \bigcap\{v(\Gamma):\gamma \in \Gamma\};\\
			\mathcal{M}\Vdash \Gamma \text{ iff for all $w \in W$, } w\Vdash \Gamma.
		\end{align}
		\begin{notation}If $\Gamma = \{\gamma\}$, we drop the parentheses and write $w \Vdash \gamma$ (or $\cc M \Vdash \gamma$) instead of $w \Vdash \{\gamma\}$ (or $\cc M \Vdash \{\gamma\}$). Moreover, since in what follows we will be dealing with some {\em submodels}, it is sometimes convenient to stress to which universe a world belongs;  given a sphere model $\cc M = (W, \cc S, v)$, $w \in W$, we then write
			\begin{equation}
				\mathcal{M}, w \Vdash \varphi \mbox{ iff } w \Vdash \varphi.
			\end{equation}
		\end{notation}

		The theorems of $\Cg$ and $\Cl$ (since they are the same by Theorem \ref{thm:samethm}) are exactly the set of formulas true in all sphere models, i.e. the set of $\cc L$-formulas $\varphi$ such that $\cc M \Vdash \varphi$ for all sphere models $\cc M$.
		We now introduce two natural notions of semantical consequence, in close analogy with the local and global consequence relations of modal logic, and we will see by the end of this section that they are sound and (finitely strongly) complete with respect to, respectively, $\Cl$ and $\Cg$. In the last section of the paper we will actually show the strong completeness for both consequence relations.
		
		\begin{definition}
			Let $\vv K$ be a class of sphere models. 
			\begin{enumerate}
				\item  We define the {\em global} $\vv K$-consequence relation on sphere models as: $\Gamma \vDash_{\vv K,g} \varphi$ if and only if  for all sphere models $\mathcal{M} \in \vv K$, $\mathcal{M}\Vdash \Gamma$ implies $\mathcal{M}\Vdash\varphi$.
				\item  We define the \emph{local} $\vv K$-consequence relation on sphere models as: $\Gamma \vDash_{\vv K,\ell} \varphi$ if and only if for all sphere models $\mathcal{M}=\langle W, \mathcal{S}, v \rangle \in \vv K$ and all $w \in W$, $w\Vdash\Gamma$ implies $w\Vdash\varphi$.
			\end{enumerate}
		\end{definition}

		\begin{notation}
			When $\vv K$ is the class of all sphere models we write $\vDash_{g}$ for $\vDash_{\vv K,g}$ and $\vDash_{\ell}$  for $\vDash_{\vv K,\ell}$.	
		\end{notation}
		The following is a direct consequence of the definitions.
		\begin{theorem}\label{thm:samevalid}
			Given any $\cc L$-formula $\varphi$, $\vDash_{g} \varphi$ if and only if $\vDash_{l} \varphi$.
		\end{theorem}
		
		Lewis \cite{Lewis1973} considers the classes of sphere models corresponding to the main axiomatic extensions of the system $\mathbf{V}$; those classes are listed in the following definition:
		\begin{definition}\label{def:classes}
			Let $\mathcal{M}=\langle W, \mathcal{S}, v \rangle$ be a sphere model. 
			\begin{enumerate}
				\item $\mathcal{M}$ is \emph{normal} if each $\mathcal{S}(w)$ is such that $\bigcup\mathcal{S}(w)\neq \emptyset$.
				\item $\mathcal{M}$ is \emph{totally reflexive} if for each $\mathcal{S}(w)$, $w \in  \bigcup\mathcal{S}(w)$.
				\item $\mathcal{M}$ is \emph{weakly centered} if for each $\mathcal{S}(w)$ there is a $S \in \mathcal{S}(w)$. such that $S\neq\emptyset$, and for each $S \in \mathcal{S}(w)$ if $S \neq\emptyset$ then $w \in S$.
				\item $\cc M$ is \emph{centered} if each $\mathcal{S}(w)$ is centered on $w$, i.e., $\{w\} \in \mathcal{S}(w)$.
				\item $\cc M$ is \emph{Stalnakerian} if for any $w \in W$, and any $\cc L$-formula $\varphi$ such that $v(\varphi) \cap \bigcup S(w) \neq \emptyset$, there is some $X \in S(w)$ and $y \in W$ such that $v(\varphi) \cap X = \{y\}$.
				\item $\mathcal{M}$ is \emph{uniform} if for each $\mathcal{S}(w)$ and $\mathcal{S}(v)$, $\bigcup\mathcal{S}(w)=\bigcup\mathcal{S}(v)$. 
				\item  $\mathcal{M}$ is \emph{absolute} if for each $\mathcal{S}(w)$ and $\mathcal{S}(v)$, $\mathcal{S}(w)=\mathcal{S}(v)$. 
				\item  $\mathcal{M}$ is \emph{weakly trivial} if for each $w \in W$, $W$ is the only non-empty member of $\mathcal{S}(w)$.
				\item $\mathcal{M}$ is \emph{trivial} if $W=\{w\}$ and $\mathcal{S}(w)=\{\{w\}, \emptyset\}$.
			\end{enumerate}
		\end{definition}
		Notice that if $\cc M=\langle W, \mathcal{S}, v \rangle$ is \emph{centered}, $\{w\}$ is the smallest sphere of $S(w)$ for all $w \in W$. Stalnakerian sphere models are complete with respect to Stalnaker logic of conditionals \cite{Stalnaker1968, Lewis1973}.
\begin{notation}
		In some of the following result it will be convenient to clarify the following notation. Given  an axiomatic extension of $\mathbf{LV}$ by $\Sigma$ be a subset of the axioms $\{\alg W, \alg C, \alg N, \alg T, \alg S, \alg U, \alg A\}$, we denote by $\vv K_{\Sigma}$ the corresponding class of sphere models, defined by the corresponding conditions in Definition \ref{def:classes}.
		\end{notation}
In particular then, the class of normal spheres by $\vv K_N$, totally reflexive by $\vv K_T$, weakly centered by $\vv K_W$, centered by $\vv K_C$, Stalnakerian by $\vv K_S$, uniform by $\vv U$, absolute by  $\vv K_A$. Moreover, Lewis shows that weak triviality corresponds to the combination of axioms $\mathbf{WA}$ thus we denote the class by $\vv K_{WA}$, and triviality is linked to $\mathbf{CA}$ so the class of models is $\vv K_{CA}$.
\begin{theorem}[\cite{Lewis1973}]\label{thm:lewiscompleteness}
	Let $\Sigma$ be a subset of the axioms $\{\alg W, \alg C, \alg N, \alg T, \alg S, \alg U, \alg A\}$. Then for any $\cc L$-formula $\varphi$, $\vdash_{\mathbf{LV\Sigma}} \varphi$ iff $\models_{\vv K_{\Sigma}} \varphi$.
\end{theorem}

		We now show that the global consequence relation can be characterized by means of the local one. In order to do that, we will introduce a useful tool. In close analogy to the case of Kripke models, we will see how to manipulate a sphere model in order to obtain a new model, preserving the satisfaction of formulas. Namely, we will prove some invariance results for Lewis's sphere semantics of counterfactuals.
		
		\begin{definition}\label{submodel}
			Let $\Sigma=\langle W, \mathcal{S}, v \rangle$ and $\Sigma'=\langle W', \mathcal{S}', v' \rangle$ two sphere models. We say that $\Sigma'$ is a {\em submodel} of $\Sigma$ if and only if:
			\begin{enumerate}
				\item $W' \subseteq W$
				\item $\mathcal{S}'$ is the restriction of $\mathcal{S}$ to $W'$, i.e. for all $w \in W'$, $\mathcal{S}'(w)=\mathcal{S}(w) \cap \mathcal{P}(\mathcal{P}(W'))$.
				\item $v'$ is the restriction of $v$ to $W'$, i.e. for any $\cc L$-formula $\varphi$, $v'(\varphi) = v(\varphi) \cap W'$. 
			\end{enumerate}
		\end{definition}
		
		We now consider a special class of submodels, namely {\em generated submodels}.
		\begin{definition}\label{gensubmodel}
			Let $\mathcal{M}=\langle W, \mathcal{S}, v \rangle$ and $\mathcal{M}'=\langle W', \mathcal{S}', v' \rangle$ be two sphere models. We say that $\mathcal{M}'$ is a {\em generated submodel} of $\mathcal{M}$ if $\mathcal{M}'$ is a submodel of $\mathcal{M}$ such that for all $w' \in W'$, $\mathcal{S}'(w')=\mathcal{S}(w')$.
		\end{definition}
		In other words, in order to obtain a generated submodel of some sphere model $\mathcal{M}=\langle W, \mathcal{S}, v \rangle$ one needs to select a subset $W' \sse W$ in such a way that, for each $w' \in W'$, all the worlds belonging to the spheres of $w'$ also belong to $W'$. This particular type of generated submodel will play a key role in our analysis. Let us show how one can effectively construct such a submodel.
		
		Consider a sphere model $\mathcal{M}=\langle W, \mathcal{S}, v \rangle$ and a subset $X \subseteq W$. Let us define a binary relation $R_{\cc S} \subseteq W\times W$ as follows: for all $w, u \in W$,
		\begin{equation}\label{eq:RS}
			w \,R_{\cc S}\, u\Leftrightarrow u\in \bigcup\mathcal{S}(w)
		\end{equation}
		Thus, $wR_{\cc S}u$ if and only if $u$ appears in the system of spheres associated to $w$.
		Now, for all $n \in \mathbb{N}$, we inductively define a relation $R_{\cc S}^n \subseteq W \times W$ in the following way:
		\begin{itemize}
			\item $wR_{\cc S}^0u \Leftrightarrow w=u$
			\item $wR_{\cc S}^{n+1}u \Leftrightarrow$ there is $z \in W$ such that $wR_{\cc S}^nz$ and $zR_{\cc S}u$.
		\end{itemize}
		We refer to $R_{\cc S}$ as the {\em accessibility relation of $\cc M$}. Intuitively, the relation $wR_{\cc S}^nv$ indicates that there are $n$ steps needed to reach the world $u$ starting from $w$, where every steps is given by checking a set of spheres. 
		Now, consider the subset $M_X \subseteq W$ defined as follows: 
		\begin{equation}\label{eq:MX}
			M_X= \{w \in W : uR_{\cc S}^nw \text{ for some } n \in \mathbb{N} \text{ and } u \in X\}
		\end{equation}
		Namely, $M_X$ is the set of all worlds in $W$ that are reachable from a member of $X$ by a finite number of steps via the accessibility relation $R_{\cc S}$. We shall now define the sphere model 
		\begin{equation}\label{eq:MMX}
			\mathcal{M}_X=\langle M_X, \mathcal{S}_X, v_X \rangle
		\end{equation} where:
		\begin{itemize}
			\item $\mathcal{S}_X$ is the restriction of $\mathcal{S}$ to $M_X$, i.e. for all $w \in W_X$, $\mathcal{S}_X(w)=\mathcal{S}(w) \cap \mathcal{P}(\mathcal{P}(M_X))$.
			\item $v_X$ is the restriction of $v$ to $M_X$, i.e. $v_X(\varphi) = v(\varphi) \cap M_X$ for all $\cc L$-formulas $\varphi$.
		\end{itemize} 	
		We say that a submodel $\cc M'$ of $\cc M$ is {\em smaller} than another submodel $\cc M^*$ if the domain of $\cc M'$ is contained in the domain of $\cc M^*$.
		Then:
		\begin{lemma}\label{lem:smallestsm}
			Let $\mathcal{M}=\langle W, \mathcal{S}, v \rangle$ be a sphere model, $X \sse W$, and let $\cc M_X$ be defined as above. Then 
			$\mathcal{M}_X$ is the smallest generated submodel of $\mathcal{M}$ whose domain contains $X$.
		\end{lemma}
		
		\begin{proof}
			First observe that $\mathcal{M}_X$ is a submodel of $\mathcal{M}$ by definition; it is also easily seen to be a generated submodel of $\mathcal{M}$. Indeed by the definition of $M_X$ the following holds: \[\text{if } w \in M_X \text{ and } wR_{\cc S}u, \text{ then } u \in M_X\]
			Equivalently,
			\[\text{if } w \in M_X \text{ and } u \in \bigcup\mathcal{S}(w), \text{ then } u \in M_X,\]
			i.e. $\mathcal{S}_X(w)=\mathcal{S}(w) \cap \mathcal{P}(\mathcal{P}(M_X)) = \cc S_X(w)$ and then $\mathcal{M}_X$ is a generated submodel of $\mathcal{M}$ by definition.
			By the very same equivalence and the definition of $M_X$, it follows that if  $\mathcal{M}^*=\langle W^*, \mathcal{S}^*, \vDash^* \rangle$ is any other generated submodel of $\cc M$ such that $X \sse W^*$, necessarily $M_X$ is contained in $W^*$.
			Therefore, $\mathcal{M}_X$ is the smallest submodel of $\mathcal{M}$ generated by $X$.
		\end{proof}
		
		\begin{definition}\label{pgensubmodel}
			Consider a sphere model $\mathcal{M}=\langle W, \mathcal{S}, v \rangle$, and let $X \subseteq W$; we call {\em submodel generated by $X$} the smallest submodel of $\mathcal{M}$ whose domain contains $X$. Moreover, we call {\em centered} or {\em point-generated} a submodel of $\mathcal{M}$ generated by a singleton.
		\end{definition}
		Notice that by Lemma \ref{lem:smallestsm} above, the submodel of $\cc M$ generated by $X$ is exactly $\cc M_X$.
		Importantly, all generated submodels preserve the validity of formulas, as the following lemma shows.
		
		\begin{lemma}\label{invariance}
			Let $\mathcal{M}=\langle W, \mathcal{S}, v \rangle$ be a sphere model, and let $\mathcal{M}'=\langle W', \mathcal{S}', v' \rangle$ be a generated submodel of $\mathcal{M}$. The following holds for all $w \in W'$, and all $\cc L$-formulas $\varphi$:
			\[\mathcal{M}, w \vDash \varphi \Leftrightarrow \mathcal{M}', w\vDash \varphi \]
		\end{lemma}
		\begin{proof}
			The statement can be easily proved by induction on the construction of the formula $\varphi$. In particular, the base case where $\varphi$ is a variable and the inductive cases given by the classical connectives (i.e. $\varphi = \psi * \gamma$ for $* \in \{\land, \lor, \to\}$) directly follow from the fact that $v'$ is the restriction of $v$ to $W'$. The inductive case $\varphi = \psi \ba \gamma$ follows from the fact that that $v'$ is the restriction of $v$ to $W'$ and that $\mathcal{S}'(w) = \cc S(w)$ for all $w \in W'$. 
		\end{proof}
		Moreover, the following is a direct consequence of the definitions.
		\begin{lemma}\label{lemma:invariance2}
			All the classes of spheres in Definition \ref{def:classes} are closed under generated submodels.
		\end{lemma}
		Before showing the characterization of the global consequence relation by means of the local consequence, we need another technical result.
		Observe that, given a sphere model $\mathcal{M}=\langle W, \mathcal{S}, v \rangle$ and $w \in W$, $w \vDash \neg \varphi \boxright \varphi$ if and only if $\bigcup\mathcal{S}(w) \cap v(\neg \varphi)=\emptyset$, or equivalently $\bigcup\mathcal{S}(w)\subseteq v(\varphi)$. Recall that $ \square \varphi := \neg \varphi \boxright \varphi$.
		It is then straightforward that, given a sphere model $\mathcal{M}=\langle W, \mathcal{S}, \vDash \rangle$, $\Box$ can be characterized by means of the relation $R_{\cc S}$ defined in (\ref{eq:RS}) by: $w R_{\cc S} v$ if and only if $v \in \bigcup\mathcal{S}(w)$. 
		\begin{lemma}\label{lem:semantic clauses}
			Let $\mathcal{M}=\langle W, \mathcal{S}, v \rangle$ be a sphere model; given any $w \in W$ and $\cc L$-formula $\varphi$, the following are equivalent:
			\begin{enumerate}
				\item $w \vDash \Box \varphi$;
				\item $\bigcup\mathcal{S}(w)\subseteq v(\varphi)$;
				\item $wR_{\cc S}u$ implies $u\vDash \varphi$.
			\end{enumerate} 
		\end{lemma}
		One can then easily show that $\Box$ is a modal operator, in the following sense.
		\begin{proposition}\label{prop:modalbox}
			The following hold for all $\varphi,\psi \in Fm_{\cc L}$:
			\begin{enumerate}
				\item $\models_g \Box(\varphi \to \psi) \to  (\Box \varphi \to \Box \psi)$;
				\item $\varphi \models_g \Box \varphi$;
				\item $\models_l \varphi$ implies $\models_l \Box \varphi$;
				\item $\models_g \Box(\varphi \land \psi) \leftrightarrow (\Box \varphi \land \Box \psi)$.
			\end{enumerate}	
		\end{proposition}
		\begin{proof}
			Let us start with (1); consider any sphere model $\mathcal{M}=\langle W, \mathcal{S}, v \rangle$, and let $w \in W$. Suppose $w 	\vDash \Box(\varphi \to \psi)$, i.e. by Lemma \ref{lem:semantic clauses} $\bigcup\mathcal{S}(w)\subseteq v(\varphi \to \psi)$. Therefore, if $w \vDash \Box \varphi$, or equivalently $\bigcup\mathcal{S}(w)\subseteq v(\varphi)$, it follows that $\bigcup\mathcal{S}(w)\subseteq v(\psi)$; applying Lemma \ref{lem:semantic clauses} again, we get that $w \vDash \Box \psi$, which proves the claim.
			
			Let us now prove (2); one needs to show that for all sphere models $\cc M = (W, S, v))$, $\cc M \Vdash \varphi$ implies $\cc M \Vdash \Box \varphi$, which is an easy consequence of Lemma \ref{lem:semantic clauses}, since if  $\cc M \Vdash \varphi$ every $w \in v(\varphi)$. (3) can be proved analogously, while (4) follows from the previous points.
		\end{proof}
		Let us  define inductively an operator $\Box^n$, that iterates $\Box$, for $n\in \mathbb{N}$:
		\begin{align}
			\Box^0\varphi&:=\;\varphi\\
			\Box^{n+1}\varphi&:=\; \Box\Box^{n}\varphi
		\end{align}	
		We are now ready to characterize the connection between the local and global consequence relations.
		\begin{theorem}\label{th:chgloballocal}
		Let $\vv K$ be a class of sphere models closed under generated submodels.
			For all sets of $\cc L$-formulas $\Gamma$ and $\cc L$-formula $\varphi$ the following are equivalent:
			\begin{enumerate}
				\item $\Gamma \models_{\vv K,g} \varphi$;
				\item $\{\Box^n\gamma: n \in \mathbb{N}, \gamma \in \Gamma\}\models_{\vv K,l} \varphi$.
				\item There exists $n_0 \in \mathbb{N}$ such that $\{\Box^n \gamma: \gamma \in \Gamma, n \leq n_0\} \models_{\vv K,l} \varphi$
			\end{enumerate}
		\end{theorem}
		\begin{proof}
		Observe that (3) trivially implies (2); we verify (2) $\Rightarrow$ (1)  by contrapositive. Assume $\Gamma \not\models_{\vv K,g} \varphi$; i.e., there is a sphere model $\mathcal{M}=\langle W, \mathcal{S}, v \rangle \in \vv K$ such that $w \Vdash \gamma$ for all $w \in W$,  $\gamma \in \Gamma$, and for some $u \in W$, $u \nVdash \varphi$. By the definition of $\Box^n$ and Lemma \ref{lem:semantic clauses}, it follows that $w \Vdash \Box^n \gamma$ for all $\gamma \in \Gamma, n \in \mathbb{N}$. Thus in particular $u \Vdash \Box^n\gamma$ for all $\gamma \in \Gamma, n \in \mathbb{N}$, but $u \nVdash \varphi$. Therefore
			$\{\Box^n\gamma\mid n \in \mathbb{N}\text{ and } \gamma \in \Gamma\}\not\models_l \varphi$.
			
			Lastly, we prove (1) implies (3), again by contrapositive; assume that for all $m \in \mathbb{N}$, $\{\Box^n\gamma\mid n \leq m\text{ and } \gamma \in \Gamma\}\not\models_l \varphi$. Thus, fixing some $m \in \mathbb{N}$, there is a sphere model $\mathcal{M}=\langle W, \mathcal{S}, v \rangle \in \vv K$ and $x \in W$ such that $x \Vdash \Box^n\gamma$ for all $n \leq m$ and  $\gamma \in \Gamma$ but $x \nVdash \varphi$. 
			Consider the submodel generated by $\{x\}$, $\mathcal{M}_x=\langle M_x, \mathcal{S}_x, v_x \rangle$, where $$M_x = \{w \in W: x R^n_{\cc S}w \mbox{ for some } n \in \mathbb{N}\}.$$ 
			
			By Lemma \ref{lem:smallestsm}, $\mathcal{M}_x$ is a sphere model and it is in $\vv K$ since $\vv K$ is closed under generated submodels by assumption;
			 moreover, by Lemma \ref{invariance} we have that for all $w \in W_x$ $$\mathcal{M}, w \Vdash \varphi \mbox{ if and only if } \mathcal{M}_x, w \Vdash \varphi.$$ 
			Hence, in particular, $\mathcal{M}_x, x \Vdash \Box^n\gamma$ for all $n \leq m,\gamma \in \Gamma$ but $\mathcal{M}_x, x \nVdash \varphi$. 
			We now prove that for all $w \in W_x$, $w \Vdash \Gamma$, which will conclude the proof by showing that $\cc M_x \Vdash \Gamma$ but $\cc M_x \not\Vdash \varphi$ (since $\mathcal{M}_x, x \nVdash \varphi$). 
			By definition, all elements $w \in W_x$ are such that $x R^m_{\cc S} w$ for some $m \in \mathbb{N}$; we show by induction on $k$ that $x R^k_{\cc S} w$ implies $w \vDash \Box^n \gamma$ for all $n \leq m, \gamma \in \Gamma$. 
			\begin{itemize}
				\item 	If $k = 0$, we get $w = x$ and thus by assumption $x \Vdash \Box^n\gamma$ for all $\gamma \in \Gamma,n\leq m$.
				\item Assume that the inductive hypothesis holds for $k < m$, we show it for $k+1$. Suppose $x R^{k+1}_{\cc S}w$, i.e. by definition of $R^{k+1}_{\cc S}$, there is some $z \in W_x$ such that  $x \,R^{k}_{\cc S}\,z \,R_{\cc S}\,w$. By inductive hypothesis $z \Vdash \Box^n \gamma$ for all $n\leq m, \gamma \in \Gamma$. Thus Lemma \ref{lem:semantic clauses} implies that also $w \Vdash \Box^n \gamma$ for all $n\leq m, \gamma \in \Gamma$. 
			\end{itemize}
			Therefore, we have shown that, in particular, all elements $w \in W_x$ are such that $w \vDash \Box^0 \gamma = \gamma$ for all $\gamma \in \Gamma$, which concludes the proof.
		\end{proof}
		In the next section we will use the last result to prove a deduction theorem and a finite strong completeness result for the global consequence relation.

		\subsection{Completeness and deduction theorem}
		Lewis proves the completeness of what we called the local consequence relation with respect to the logic $\vdash_\mathbf{LV}$ (Theorem \ref{thm:lewiscompleteness}); he also shows a deduction theorem with respect to the classical implication:
		\begin{theorem}[Local Deduction Theorem]
			For all $\Gamma \cup \{\varphi, \psi\} \subseteq For_{\cc L}$, the following are equivalent:
			\begin{enumerate}
				\item $\Gamma, \psi \models_l \varphi $,
				\item $\Gamma \models_l \psi \to \varphi$.
			\end{enumerate}	
		\end{theorem}
		It follows that there is a finite strong completeness result for the local consequence, i.e.:
		\begin{theorem}[\cite{Lewis71}]\label{thm:lewiscomplete}
		Let $\Sigma$ be a subset of the axioms $\{\alg W, \alg C, \alg N, \alg T, \alg S, \alg U, \alg A\}$. Then for any finite set of $\cc L$-formulas $\Gamma$ and $\cc L$-formula $\varphi$, $\Gamma \vdash_{\mathbf{LV\Sigma}} \varphi$ iff $\Gamma \models_{\vv K_{\Sigma},\,l} \varphi$.
		\end{theorem}
		
		In this subsection we prove the analogous result for the stronger deductive systems and the corresponding global consequence relation; but first, some technical results.
		\begin{lemma}\label{lem:derivedrule1}
			Consider any $\cc L$-formula $\varphi$, then $\varphi \vdash_\mathbf{GV} \Box^n\varphi$ for all $n \in \mathbb{N}$.
		\end{lemma}	
		\begin{proof}
			The claim is easily shown by induction on $n$; indeed the case $n = 0$ is obvious, and the inductive case is given by one application of (DWC$_0$): $ \varphi \vdash_\mathbf{GV} \neg \varphi \to \varphi$, which holds for $\Cg$ by Theorems \ref{thm:firstax} and \ref{thm:secondax}.	
		\end{proof}
		\begin{proposition}\label{lem:derivedrule}
		Let $\mathbf{L}$ be any axiomatic extension of $\Cg$.
			For all $\Gamma \cup \{\varphi\}\subseteq For_{\cc L}$ the following are equivalent: 
			\begin{enumerate}
				\item $\Gamma \vdash_{\alg L} \varphi$;
				\item $\{\Box^n\gamma\mid \gamma \in \Gamma \text{ and }n \in \mathbb{N}\} \vdash_{\alg L} \varphi$.
				\item There exist a finite subset $\Gamma_0 \sse \Gamma$ and $n_0 \in \mathbb{N}$ such that $\{\Box^n\gamma\mid \gamma \in \Gamma_0 \text{ and }n\leq n_0\} \vdash_{\alg L} \varphi$.
			\end{enumerate}
		\end{proposition}
		\begin{proof}
			The fact that (1) implies (2) is obvious, since $\Gamma \subseteq \{\Box^n\gamma\mid \gamma \in \Gamma \text{ and }n \in \mathbb{N}\}$. 
			For the converse, let us assume that $\{\Box^n\gamma\mid \gamma \in \Gamma \text{ and }n \in \mathbb{N}\} \vdash_\mathbf{GV} \varphi$;
			By Lemma \ref{lem:derivedrule1}, we have that $\Gamma \vdash_\mathbf{GV} \Box^n \gamma$ for all $\gamma \in \Gamma$ and $n \in \mathbb{N}$, and thus also $\Gamma \vdash_\mathbf{GV} \varphi$. 
			Lastly, (2) and (3) are equivalent since $\vdash_{\alg L}$ is a finitary consequence relation.
		\end{proof}
		We now have all the ingredients to prove our completeness result.	
		\begin{theorem}[Soundess and finite strong completeness]\label{th: globalompleteness}
		Let $\Sigma$ be a subset of the axioms $\{\alg W, \alg C, \alg N, \alg T, \alg S, \alg U, \alg A\}$.
			For all finite subsets $\Gamma \cup \{\varphi\}\subseteq For_{\cc L}$, the following are equivalent: 
			\begin{enumerate}
				\item $\Gamma \vdash_{\mathbf{GV\Sigma}} \varphi$;
				\item $\Gamma \models_{\vv K_{\Sigma}, g} \varphi$.
			\end{enumerate}
		\end{theorem}
		
		\begin{proof}
			
			The soundness follows from the facts that: (MP) and (C) are easily seen to be sound with respect to sphere models, and the axioms of $\mathbf{GV\Sigma}$ are the same axioms of $\mathbf{LV\Sigma}$, which is sound with respect to the same class of models with respect to the local consequence relation (Theorem \ref{thm:lewiscomplete}), and the latter has the same valid formulas as the global one (Theorem \ref{thm:samevalid}).
			
			We prove completeness by contrapositive; assume $\Gamma \nvdash_{\mathbf{GV\Sigma}} \varphi$. By Lemma \ref{lem:derivedrule}, we have that for all $m \in \mathbb{N}$, $\{\Box^n \gamma \mid \gamma \in \Gamma \text{ and }n \leq m\} \nvdash_s \varphi$. By the fact that all deductions of $\vdash_{\mathbf{LV\Sigma}} $ are deductions of $\vdash_{\mathbf{GV\Sigma}}$ (Lemma \ref{lemma:Cgstronger}), and by the finite strong completeness of $\vdash_{\mathbf{LV\Sigma}} $ in Theorem \ref{thm:lewiscomplete}, we have that $\{\Box^n \gamma\mid \gamma \in \Gamma \text{ and }n \leq m\} \not\models_{\vv K_{\Sigma},\, l} \varphi$ for all $m \in \mathbb{N}$. Theorem \ref{th:chgloballocal} then yields that $\Gamma \not\models_{\vv K_{\Sigma}, g} \varphi$ and the proof is complete.
		\end{proof}
	
		We will now show another useful fact, i.e., that the global consequence relation has a deduction theorem. However, it generally does not have the classical deduction theorem, as the following example show.
		\begin{example}
		By Lemma \ref{lem:derivedrule1} and the finite strong completeness in Theorem \ref{th: globalompleteness}, $\varphi\models_g \Box \varphi$ for any $\cc L$-formula $\varphi$. However, it is easily seen that in general  $\not\models_g \varphi \to \Box \varphi$. 
			Consider the following sphere model $\cc M = (W, \cc S, v)$ such that 
			\begin{itemize}
				\item $W = \{w_1, w_2\}$;
				\item $\cc S(w_1) = \{\{w_2\}\}$; $\cc S(w_2) = \{w_2, \{w_1, w_2\}$;
				\item $v$ is such that it maps a propositional variable $p \in Var$ to $v(p) = \{w_1\}$.
			\end{itemize}
			Note that then $v(\Box p) = \emptyset$, and therefore $w \Vdash p$ but $w \nVdash \Box p$; hence $\cc M \nVdash p \to \Box \varphi$. One can easily adapt this same example for any class of sphere models that allow at least two different worlds in $W$ (thus, except for the trivial class of models). 
		\end{example}
		Nonetheless:	
		\begin{theorem}[Global Deduction Theorem]\label{thm:globaldeduction}
		Let $\vv K$ be a class of sphere models closed under generated submodels.
			For all $\Gamma \cup \{\varphi, \psi\} \subseteq For_{\cc L}$, the following are equivalent:
			\begin{enumerate}
				\item $\Gamma, \psi \models_{\vv K,g} \varphi $,
				\item there is $n \in \mathbb{N}$ such that $\Gamma \models_{\vv K,g} \left(\bigwedge\limits_{m \leq n} \Box^m\psi\right)\to \varphi$,
			\end{enumerate}	
		\end{theorem}
		\begin{proof}
			Let us start by proving that (1) implies (2); assume  $\Gamma, \psi \models_g \varphi$. By Theorem \ref{th:chgloballocal}, this is equivalent to the fact that there exists $n_0 \in \mathbb{N}$ such that $$\{\Box^n \gamma: \gamma \in \Gamma, n \leq n_0\} \cup\{\Box^n\psi: n \in \mathbb{N}\} \models_{\vv K,l} \varphi$$  
			By the deduction theorem for the local consequence, it follows that  $$\{\Box^n \gamma: \gamma \in \Gamma, n \leq n_0\} \models_{\vv K, l} \left (\bigwedge\limits_{k \leq n_0} \Box^k\psi\right ) \to  \varphi.$$ Using Theorem \ref{th:chgloballocal} again, we have that $\Gamma \models_{\vv K,g} \left (\bigwedge\limits_{k \leq n_0} \Box^k\psi\right ) \to  \varphi$. 
			
			We now show that (2) implies (1); assume that there is a $n \in \mathbb{N}$ such that $\Gamma \models_{\vv K,g} \left(\bigwedge\limits_{m \leq n} \Box^m\psi\right)\to \varphi$. Equivalently, for all sphere models $\mathcal{M}=\langle W, \mathcal{S}, v \rangle \in \vv K$, we have that if $\mathcal{M}\Vdash \gamma$ for all $\gamma \in \Gamma$, then $\mathcal{M}\Vdash \left(\bigwedge\limits_{m \leq n} \Box^m\psi\right)\to \varphi$.		
			Consider then a sphere model $\mathcal{M} = (W, \cc S, v) \in \vv K$ such that $\mathcal{M}\Vdash \gamma$ for all $\gamma \in \Gamma$ and $\mathcal{M} \Vdash \psi$. By assumption, we then have that for any world $w \in W$, $w \Vdash \psi$ and $w \Vdash \left(\bigwedge\limits_{m \leq n} \Box^m\psi\right)\to \varphi$. By Lemma \ref{lem:semantic clauses}, this implies that for all $w \in W$, $w \Vdash \bigwedge\limits_{m \leq n} \Box^m\psi$. Therefore, by modus ponens, we have that for all $w \in W$, $w \Vdash \varphi$ as well; i.e. $\cc M \Vdash \varphi$. Hence $\Gamma, \psi \models_{\vv K,g} \varphi$. 
		\end{proof}
		
		\noindent
		
		
		We are now ready to proceed our investigation towards an algebraic study of Lewis's hierarchy of logics for counterfactual conditionals.
		
		\section{Algebraic semantics}
		In this section we show that the global consequence relation is algebraizable in the sense of Blok-Pigozzi, and we study the equivalent algebraic semantics, that are varieties of Boolean algebras with an extra operator $\ba$. Moreover, we show that such varieties of algebras give a semantics for the logics with a local consequence relation as well, in the sense that the latter are the logics preserving the degrees of truth of these algebras.

		\subsection{Global equivalent algebraic semantics}
		The algebraizability of $\Cg$ follows from the fact that it is an implicative logic.
		\begin{theorem}
			Let $\cc L$ be any axiomatic extension of $\Cg$. Then $\cc L$ is an implicative logic.
		\end{theorem}
		\begin{proof}
			We need to show that the conditions of Definition \ref{def:implicativelogic} hold; (\ref{item:implic}) follows from Lemma \ref{lemma:equiv} and the others follow from the fact that $\to$ is a Boolean implication. 
		\end{proof} 
		Moreover, since $1$ is a theorem of $\Cg$ (recall the discussion in the preliminaries), we get the following.
		\begin{theorem}\label{thm:Cgalg}
			$\Cg$ is algebraizable with defining equation $\tau(x) = \{x \approx 1\}$ and equivalence formula $\Delta(x,y) = \{ x \leftrightarrow y\}$. 
		\end{theorem}
		An important consequence of algebraizability is that axiomatic extensions of algebraizable logics are also algebraizable with the same equivalence formulas and defining equations. Moreover, the lattice of axiomatic extensions of the logic is dually isomorphic to the subvariety lattice of the algebraic semantics whenever the latter is a variety, as it is the case here.
		\begin{corollary}\label{cor:alg}
			Let $\cc L$ be any axiomatic extension of $\Cg$, axiomatized relatively to $\Cg$ by the set of axioms $\Phi$. Then the equivalent algebraic semantic of $\cc L$ is the subvariety of $\LCA$ axiomatized by $\tau(\Phi) = \{\varphi \approx 1 : \varphi \in \Phi\}$. In particular:
			$$\Gamma \vdash_{\cc L} \f\quad\text{iff}\quad \tau(\Gamma) \vDash_{\LCA + \tau(\Phi)} \tau(\f).$$
		\end{corollary}
		We can then get an axiomatization of the equivalent algebraic semantics of $\Cg$ (and its extensions); it turns out that equations suffice, in particular because the rules (MP) and (C) are translated to quasi-equations that can be shown to hold as a consequence of the other axioms. In other words, the equivalent algebraic semantics is a {\em variety} of algebras.
			\begin{definition}
			A \emph{Lewis variably strict  conditional algebra}, or {\em $\vv V$-algebra} for short, is an algebra $\mathbf{C}= (C, \wedge, \vee, \to, \boxright, 0, 1)$ where $(C, \wedge, \vee, \to, 0, 1)$ is a Boolean algebra and $\boxright$ is a binary operation such that, for all $x, y, z \in C$:
			\begin{enumerate}
				\item\label{C1} $x \boxright x=1$
				\item\label{C2} $((x \boxright y) \wedge (y \boxright x)) \leq ((x \boxright z)\leftrightarrow(y \boxright z))$
				\item\label{C3} $((x \vee y)\boxright x) \vee ((x \vee y)\boxright y) \vee (((x \vee y)\boxright z)\leftrightarrow((x \boxright z)\wedge(y \boxright z))=1$
				\item\label{C4} $x \boxright (y \wedge z) = (x \boxright y)\wedge (x \boxright z)$
			\end{enumerate}	
			We denote the variety of $\vv V$-algebras with $\LCA$.
		\end{definition}
		\begin{theorem}
			$\LCA$ is the equivalent algebraic semantics of $\Cg$.
		\end{theorem}
		\begin{proof}
			By Corollary \ref{cor:alg}, the equivalent algebraic semantics of $\Cg$ is axiomatized by:
			\begin{itemize}
				\item $x \boxright x=1$
				\item $((x \boxright y) \wedge (y \boxright x)) \to  ((x \boxright z)\leftrightarrow(y \boxright z)) = 1$
				\item $((x \vee y)\boxright x) \vee ((x \vee y)\boxright y) \vee (((x \vee y)\boxright z)\leftrightarrow((x \boxright z)\wedge(y \boxright z))=1$
				\item $x \boxright (y \wedge z) = (x \boxright y)\wedge (x \boxright z)$
				\item $x = 1, x \to y = 1$ implies $y =1$
				\item $x \to y = 1$ implies $(z \ba x) \to (z \ba y) = 1$.
			\end{itemize}
			The result then follows from the fact that given $\alg C \in \LCA$, $x \to y = 1$ iff $x \leq 1$, and moreover the last two quasi-equations are easily seen to hold; in particular, note that the fact that $\ba$ is order-preserving on the right here follows from the distributivity over the meet operation on the right: if $x \leq y$ then $x = x \land y$, and so $z \ba x = z \ba (x \land y) = (z \ba  x) \land (z \ba y)$ and so $z \ba x \leq z \ba y$.
		\end{proof}
		Moreover, set again $x \Diamondright y := \neg (x \boxright \neg y)$ and $\square x := \neg x \boxright x$.
		Let $\LC$ be the subvariety of $\LCA$ further satisfying:
		\begin{equation}
		x \land y \leq x \ba y \leq x \to y	
		\end{equation}
		and $\CA$ the subvariety of $\LC$ of \emph{Lewis conditional algebras}, satisfying:
		\begin{equation}
			(x \ba y) \lor (x \ba \neg y) = 1.
		\end{equation}

		\begin{corollary}\label{cor:algebraizability}
			$\LC$ and $\CA$ are the equivalent algebraic semantics of, respectively, $\Cug$ and $\Cdg$.
	\end{corollary}
		
		The following properties can be shown via easy calculations.
		\begin{lemma}
			In every $\vv V$-algebra $\mathbf{C}=(C, \wedge, \vee, \to, \boxright, 0, 1)$, the following hold for all $x, y, z \in C$:
			\begin{enumerate}
				\item  $(x \boxright z) \wedge (y \boxright z) \leq (x \vee y)\boxright z$;
				\item  $x \boxright y \leq x \boxright (y \vee z)$;
				\item $x \to y =1$ iff $x \boxright y =1$;
				\item $0 \boxright x = 1$;
				\item $x \boxright 1 = 1$;
				\item $x \boxright 0 \leq \neg x$.
			\end{enumerate}
		\end{lemma}	
		
		Let us consider again the unary connective $\Box$ in the language as $\Box\varphi := \neg \varphi \boxright \varphi$, and its iteration
		$\Box^n\varphi$. We can show that, analogously to the case of modal algebras, the operator $\Box$ can be used to characterize congruence filters.
		\begin{definition}
			Let $\alg A \in \LCA$; a nonempty lattice filter $F \sse A$ is said to be {\em open} if $x \in F$ implies $\Box x \in F$.
		\end{definition}
		\begin{proposition}\label{prop:filters}
			Let $\alg A \in \LCA$; a nonempty lattice filter $F \sse A$ is a congruence filter if and only if it is open.
		\end{proposition}
		\begin{proof}
			The proof is based on the fact that, as a consequence of the fact that $\LCA$ is the equivalent algebraic semantics of $\Cg$, congruence filters coincide with the deductive filters induced by the logic. In other words, for every rule $\Gamma \vdash \varphi$ and every homomorphism $f$ from $Fm_{\cc L}$ to $\alg A$, if $f[\Gamma] \sse F$, then $f(\varphi) \in F$. 
			
			It is clear that every deductive filter is an open lattice filter. For the converse, consider an open lattice filter $F$; $F$ respects the axioms because it is nonempty (i.e. it contains $1$, where all the instances of the axioms are mapped as a direct consequence of the algebraizability result, Corollary \ref{cor:alg}), and it respects modus ponens since it is a lattice filter. We only need to prove that it respects the rule $\varphi \to \psi \vdash (\gamma \ba \varphi) \to (\gamma \ba \psi)$.
			
			Suppose there is an assignment $h$ to $\alg A \in \LCA$ such that $h(\varphi) = a, h(\psi) = b, h(\gamma) = c$, with $a,b,c \in A$, and assume that $a \to b \in F$. From the fact that $\varphi \to \psi \vdash (\gamma \ba \varphi) \to (\gamma \ba \psi)$ holds, by the global deduction theorem (Theorem \ref{thm:globaldeduction}) and the finite strong completeness (Theorem \ref{th: globalompleteness}) we obtain that there is some $n \in \mathbb{N}$ such that $$\vdash_{\Cg} \left(\bigwedge\limits_{m \leq n} \Box^m(\varphi \to \psi)\right)\to ((\gamma \ba \varphi) \to (\gamma \ba \psi)).$$
			This implies that the element $\left(\bigwedge\limits_{m \leq n} \Box^m(a \to b) \right)\to ((b \ba a) \to (c \ba b)) = 1 \in F$. Now, since $a \to b \in F$ and $F$ is open, $\Box^k(a \to b) \in F$ for all $k \in \mathbb{N}$. Thus, since filters are closed under finitary meets, the element $\left(\bigwedge\limits_{m \leq n} \Box^m(a \to b) \right) \in F$. Since lattice filters respect modus ponens, also $((b \ba a) \to (c \ba b)) \in F$, which shows that open lattice filters respect the rule (C). 
			
			We have then shown that open filters coincide with deductive filters, and therefore with congruence filters.
		\end{proof}
		The following is an interesting observation.
		\begin{corollary}
			Let $\alg A \in \LCA$; then the congruence filters of $\alg A$ are exactly the congruence filters of its modal reduct $(A, \land, \lor, \to, \Box, 0, 1)$.
		\end{corollary}
		
		\subsection{Local algebraic semantics}
		In this section we focus on the weaker logic $\Cl$; in particular, we will show that the latter is not algebraizable, however it still can be studied by means of $\vv V$-algebras. We will indeed see that $\vdash_\mathbf{LV}$ coincides with the logic {\em preserving degrees of truth} of $\LCA$. Let us be more precise.
		
		We call an {\em ordered algebra} a pair $(\alg A, \leq)$, where $\alg A$ is an algebra and $\leq$ is a partial order on its universe. Note that all algebras with a (semi)lattice reduct can be seen as ordered algebras, and thus in particular Boolean algebras and $\vv V$-algebras are ordered algebras.
		\begin{definition}
			Let $\vv K$ be a class of ordered algebras; the logic {\em preserving degrees of truth} of $\vv K$, in symbols $\vdash_{\vv K}^{\leq}$,  is defined as follows: for every set  $\Gamma \cup \{\varphi\}$ of formulas in the language of $\vv K$,
			$$\Gamma \vdash_{\vv K}^{\leq} \varphi$$
			if and only if for all $(\alg A, \leq) \in \vv K$, and assignment $h: Fm(Var) \to \alg A$, $a \in A$,
			$$a \leq h(\gamma) \mbox{ for every } \gamma \in \Gamma \;\;\Rightarrow\;\; a \leq h(\varphi).$$
		\end{definition}
		
		\begin{remark}\label{remark:degrees}
			If $\vv K$ is an elementary class of algebras (in particular, if $\vv K$ is a variety) with a lattice reduct, one can rephrase the above definition and say that $$\Gamma \vdash_{\vv K}^{\leq} \varphi \;\;\mbox{ iff }\;\; \vv K \models \gamma_1 \land \ldots \land \gamma_n \leq \varphi$$ for some $\{\gamma_1, \ldots, \gamma_n\} \sse \Gamma$ for some $n \in \mathbb{N}$ (see \cite[Remark 2.4]{moraschiniforth}).
		\end{remark} 
		\begin{example}
			Intuitionistic logic is the logic preserving degrees of truth of Heyting algebras, and the local consequence of the modal logic $\mathsf{K}$ is the logic preserving degrees of truth of modal algebras. 		
		\end{example}	
		Logics preserving the degrees of truth have been studied in generality in \cite{nowak,font2009}, and in residuated structures in \cite{bouesteva}.
		In order to prove the analogous result for $\Cl$, let us first state a useful lemma.
		
		\begin{lemma}\label{lemma:technical}
			For all $\cc L$-formulas $\varphi, \psi, \gamma$, and $n \in \mathbb{N}$:
			\begin{enumerate}
				\item $\vdash_\mathbf{LV} \Box^{n+1}(\varphi \to \psi) \to \Box^n((\gamma \ba \varphi) \to (\gamma \ba \psi))$;
				\item $\vdash_\mathbf{LV} \Box^{n+1}(\varphi \leftrightarrow \psi) \to \Box^n((\varphi \ba \gamma) \to (\psi \ba \gamma)) $.
			\end{enumerate}
		\end{lemma}
		\begin{proof}
			We will show that the two statements hold in sphere models. We start with (1), and proceed by induction on $n$. Let $n = 0$, we want to prove that: $$\models_l \Box (\varphi \to \psi) \to ((\gamma \ba \varphi) \to (\gamma \ba \psi));$$
			consider a sphere model $\mathcal{M}=\langle W, \mathcal{S}, v \rangle$ and let $w \in W$. Suppose $w \vDash \Box (\varphi \to \psi)$, or equivalently via Lemma \ref{lem:semantic clauses} $\bigcup\mathcal{S}(w)\subseteq v(\varphi \to \psi)$. Now, if $w \vDash \gamma \ba \varphi$, it means that there is $S \in \cc S(w)$ such that $\emptyset \neq S \cap v(\gamma) \sse v(\varphi)$; but since $\bigcup\mathcal{S}(w)\subseteq v(\varphi \to \psi)$, $S \cap v(\varphi) \sse v(\psi)$, i.e. there is $S \in \cc S(w)$ such that $\emptyset \neq S \cap v(\gamma) \sse v(\psi)$, and thus $w \vDash \gamma \ba \psi$. The inductive step follows from the fact that $\Box$ is a modal operator, more precisely that theorems are closed under $\Box$ and $\Box$ distributes over the implication (Proposition \ref{prop:modalbox}).
			
			Let us now show (2), again by induction on $n$; for $n= 0$, we prove:
			$$\models_l \Box(\varphi \leftrightarrow \psi) \to ((\varphi \ba \gamma) \to (\psi \ba \gamma)).$$
			Consider then a sphere model $\mathcal{M}=\langle W, \mathcal{S}, v \rangle$, let $w \in W$, and suppose $w \vDash \Box (\varphi \leftrightarrow \psi)$, or equivalently by Lemma \ref{lem:semantic clauses}, $\bigcup\mathcal{S}(w)\subseteq v(\varphi \leftrightarrow \psi)$. If 
			$w \vDash \varphi \ba \gamma$, it means that there is $S \in \cc S(w)$ such that $\emptyset \neq S \cap v(\varphi) \sse v(\gamma)$; 
			but since $\bigcup\mathcal{S}(w)\subseteq v(\varphi \leftrightarrow \psi)$, we get that $S \cap v(\varphi) = S \cap v(\psi)$, and therefore there is $S \in \cc S(w)$ such that $\emptyset \neq S \cap v(\psi) \sse v(\gamma)$, and thus $w \vDash \psi \ba \gamma$. The inductive step follows again from the fact that $\Box$ is a modal operator.
		\end{proof}

		Once again in close analogy with the modal framework, we have the following result, whose demonstration essentially resembles a standard completeness proof.
		\begin{theorem}\label{th:localcomp}
			$\Cl$ is the logic preserving degrees of truth of $\LCA$.
		\end{theorem}
		
		\begin{proof}
			We need to prove that for all $\Gamma \cup \{\varphi\}\subseteq Fm_{\cc L}$: $$\Gamma \vdash_\mathbf{LV} \varphi \mbox{ if and only if } \Gamma \vdash^\leq_{\LCA}\varphi.$$
			The forward direction is a usual soundness proof; note that the axioms are mapped to $1$ by the algebraizability result (Corollary \ref{cor:alg}), (MP) holds in the form $x \land (x \to y) \leq y$, and (WC) becomes $x \to y = 1$ implies $(z \ba x) \leq (z \ba y)$ which holds in $\LCA$.
			
			For the converse, we reason by contraposition; assume $\Gamma \not\vdash_\mathbf{LV} \varphi$. Then consider the relation $\theta$ defined as follows:
			$$\theta:=\{(\psi, \gamma) \in  Fm_{\cc L}\times Fm_{\cc L} : \Gamma \vdash_\mathbf{LV} \Box^n(\psi \to \gamma)\text{ and }\Gamma \vdash_\mathbf{LV} \Box^n(\gamma \to \psi)\text{ for all $n \in \mathbb{N}$} \}.$$
			We can show that $\theta$ is a congruence relation. In particular, while reflexivity and symmetry are trivial, transitivity follows the fact that $\Box$ distributes over $\to$ (Proposition \ref{prop:modalbox}); the Boolean operations are also easily shown to be respected: for $\land$, as a consequence of the fact that $\Box$ distributes over $\land$ (Proposition \ref{prop:modalbox}), and for $\neg$, by the observation that $\neg \Box^n x = \Box^n (\neg x)$.  Let us then prove that $\theta$ preserves the binary operation $\boxright$. Assume $\psi\theta\gamma$, it suffices to show that:
			$$(\delta \ba \psi,\delta \ba \gamma),(\psi \ba \delta, \gamma \ba \delta) \in \theta;$$
			that is to say, we need prove that:
			\begin{enumerate}
				\item $\Gamma \vdash_\mathbf{LV} \Box^n((\delta \ba \psi) \to (\delta \ba \gamma))$,
				\item $\Gamma \vdash_\mathbf{LV} \Box^n( (\delta \ba \gamma) \to (\delta \ba \psi))$,
				\item $\Gamma \vdash_\mathbf{LV} \Box^n((\psi \ba \delta) \to (\gamma \ba \delta))$,
				\item $\Gamma \vdash_\mathbf{LV} \Box^n((\gamma \ba \delta)\to (\psi \ba \delta))$.
			\end{enumerate}
			Given the assumption that $\Gamma \vdash_\mathbf{LV} \Box^n(\psi \to \gamma)\text{ and }\Gamma \vdash_\mathbf{LV} \Box^n(\gamma \to \psi)$ for all $n \in \mathbb{N}$, and thus also $\Gamma \vdash_\mathbf{LV} \Box^{n}(\varphi \leftrightarrow \psi)$ (since $\Box$ distributes over $\land$ by Proposition \ref{prop:modalbox}), (1)--(4) follow by Lemma \ref{lemma:technical}. Therefore, $\theta$ is a congruence, and we can consider the quotient $Fm_{\cc L}/\theta$; let us verify that the latter is an algebra in $\LCA$. Consider the axiomatization of $\LCA$, we show that for every equation $\varepsilon = \delta$ appearing in it, $(\varepsilon,\delta)\in \theta$. By the algebraizability of $\LCA$, for each such $\varepsilon = \delta$ and the completeness result of Theorem \ref{th: globalompleteness}, $\vdash_\mathbf{GV} \varepsilon \leftrightarrow \delta$; thus also $\vdash_\mathbf{LV}\varepsilon \leftrightarrow \delta$ (since they have the same theorems by Theorem \ref{thm:samethm}), and then $\vdash_\mathbf{LV} \Box^n(\varepsilon \leftrightarrow \delta)$ for all $n \in \mathbb{N}$ by the fact that theorems are closed under $\Box$ (Proposition \ref{prop:modalbox}). Thus, each $(\varepsilon,\delta)\in \theta$, and  $Fm_{\cc L}/\theta \in \LCA$.
			
			We now show that for all (nonempty) finite subsets $\Delta \sse \Gamma$, $\Delta \not\vdash^\leq_{\LCA} \varphi$. By Remark \ref{remark:degrees}, this implies that $\Gamma \not\vdash^\leq_{\LCA}\varphi$ which would complete the proof.	
			Given that $\Delta \vdash^\leq_{\LCA} \varphi$ iff $\LCA \models \bigwedge \Delta \leq \varphi$, it is enough to show that in particular $Fm_{\cc L}/\theta \not\models \bigwedge \Delta \leq \varphi$.

			Consider $\Delta = \{\delta_1, \dots, \delta_n\} \sse \Gamma$, and let $\pi$ be the natural epimorphism $\pi:  Fm_{\cc L} \to Fm_{\cc L}/\theta$; assume by way of contradiction that $\pi(\delta_1)\wedge\dots\wedge h(\delta_n)\leq \pi(\varphi)$; thus $\pi ((\delta_1\wedge\dots\wedge\delta_n) \to \varphi) = 1$. By the definition of $\theta$, this implies that in particular $\Gamma \vdash_\mathbf{LV} (\delta_1\wedge\dots\wedge\delta_n) \to \varphi$. But since $\Delta \sse \Gamma$, it follows that $\Gamma \vdash_\mathbf{LV} \delta_1\wedge\dots\wedge\delta_n$; by modus ponens we would get $\Gamma \vdash_\mathbf{LV} \varphi$, a contradiction. This completes the proof.
		\end{proof}
		Actually, the same proof works for any axiomatic extension of $\Cl$.
		\begin{corollary}
			Let $\Gamma$ be a set of formulas; $\Cl + \Gamma$ is the logic preserving degrees of truth of $\LCA + \tau(\Gamma)$.
		\end{corollary}
		As a corollary of the above theorem, we get that:
		
		\begin{corollary}\label{cor:key}
			Let $\Gamma$ be a set of formulas. For all $\varphi, \psi \in For_{\cc L}$, $\tau(\Gamma) \models_{\LCA}\varphi \app \psi$ if and only if $\varphi \vdash_{\Cl + \Gamma} \psi$ and $\psi \vdash_{\Cl + \Gamma} \varphi$.
		\end{corollary}
		We will now demonstrate that the local consequence relation is not algebraizable. We start by verifying that if it were algebraizable, then it would be algebraizable with respect to a class $\vv K$ of $\vv V$-algebras; the following proof uses standard techniques.
		
		\begin{lemma}\label{lem:basta}
			If $\Cl$ is algebraizable, its equivalent algebraic semantics is a class of algebras $\vv K \sse \LCA$.
		\end{lemma}

		\begin{proof}
			Suppose $\Cl$ is algebraizable with respect to some class of algebras $\vv K$ over the same language. Then there must be a set of equations $\tau(x)$, a set of formulas $\Delta(x, y)$, that witness the	algebraizability of $\Cl$ with respect to $\vv K$ (see the preliminary section). 
			
			We are going to show that $\vv K$ must be a class of $\vv V$-algebras, i.e., every equation that is valid in $\LCA$ must be valid in all the algebras in $\vv K$.
			Let $\varphi \app \psi$ be an equation that holds in $\LCA$; it follows that for all $\delta \in \Delta$, $\LCA \models \delta(\varphi, \varphi) \app \delta(\varphi, \psi).$
			By Corollary \ref{cor:key} the latter implies in particular that $\Delta(\varphi, \varphi)\vdash_\mathbf{LV}\Delta(\varphi, \psi)$; since by algebraizability we also get that $\vdash_\mathbf{LV} \Delta(\varphi, \varphi)$, this yields  $\vdash_\mathbf{LV} \Delta(\varphi, \psi)$. Therefore $\vv K\models \tau(\Delta(\varphi, \psi))$ and thus $\vv K\models \varphi \app \psi$. Hence, every equation that holds in $\LCA$ will also hold in $\vv K$, and then $\vv K \sse \LCA$. 	\end{proof}
		Now, recall that if a logic $\vdash$ is algebraizable with respect to a class of algebras $\vv K$, given every algebra $\mathbf{A}$ in the same signature as the algebras in $\vv K$, the lattice of deductive filter of $\alg A$ with respect to $\vdash$ is in bijective correspondence with the lattice of $\vv K$-relative congruences of $\mathbf{A}$ (i.e. the congruences $\theta$ such that $\alg A/\theta \in \vv K$), see Theorem \ref{thm:isomorphismfilters}. We will use this fact to show that $\Cl$ is not algebraizable, by showing an algebra $\alg A \in \LCA$ where the lattice of deductive filters wrt $\vdash_l$ are not in bijection with the congruences. 
		\begin{example}\label{ex:algebraiz}
			Consider the algebra $\alg A$ in the signature of $\LCA$, with Boolean reduct the 4-element Boolean algebra with domain $A = \{a, \neg a, 0, 1\}$, and the operation $\ba$ is in the following table.
			\begin{figure}[htbp]
				\begin{center}
					\begin{tabular}{c|c|c|c|c}
						$\ba$  & $0$  & $a$ & $\neg a$ & $1$ \\\hline
						$0$  & $1$  & $1$ & $1$ & $1$ \\\hline
						$a$  & $0$  & $1$ & $0$ & $1$ \\\hline
						$\neg a$  & $0$  & $0$ & $1$ & $1$ \\\hline
						$1$  & $0$  & $a$ & $\neg a$ & $1$
					\end{tabular}
					\caption{The operation table of $\ba$.}\label{fig:table}
				\end{center}
			\end{figure}
			\begin{figure}[htbp]
				\begin{center}
					\begin{tikzpicture}
						\draw (0,0) -- (1,1) -- (0,2) -- (-1,1) -- (0,0);
						\node at (0,-0.2) {$0$};
						\node at (0,2.2) {$1$};
						\node at (1.25,1) {$\neg a$};
						\node at (-1.25,1) {$a$};
						\node at (0,-0.2) {$0$};
						\draw [white] (1,1.2)-- (2,2.5);
					\end{tikzpicture}
					\caption{The Hasse diagram of $\alg A.$}\label{fig:algA}
				\end{center}
			\end{figure}
			
			One can see by easy computations that $\alg A \in \LCA$, showing that the four axioms are respected:
			\begin{enumerate}
				\item $x \ba x = 1$ for all $x \in A$, as witnessed by the diagonal of the table in Figure \ref{fig:table}.
				\item $((x \boxright y) \wedge (y \boxright x)) \leq ((x \boxright z)\leftrightarrow(y \boxright z))$ for all $x, y, z \in A$. We first present the calculations for $x = 1$; if $y = 1$ as well, the inequality becomes $1 \leq z \leftrightarrow z = 1$. If $y = a$, we have the following cases:
				\begin{align*}
					z = 1&:\;\; a \leq 1 \leftrightarrow 1 = 1;\\
					z = a&:\;\; a \leq a \leftrightarrow 1 = a;\\
					z = \neg a&:\;\; a \leq \neg a \leftrightarrow 0 = a;\\
					z = 0&:\;\; a \leq 0 \leftrightarrow 0 = 1.
				\end{align*} 
				The case for $y = \neg a$ is symmetric (in fact, the operations on $\alg A$ are symmetric, as attested by Figures \ref{fig:algA} and \ref{fig:table}); lastly, if $y = 0$ then the left-hand side of the inequality is $0$, and thus the inequality holds for all choices of $z$. 
				
				Note that if $y = 1$, the inequality holds since it is symmetric with respect to $x$ and $y$. Consider now the case where $x = a$; we do not need to consider the case where $y = 1$, and in the cases $y \in \{0, \neg a\}$ the left-hand side of the inequality is $0$. The only case left is $y = a$, where the inequality becomes $1 \leq (a \ba z) \leftrightarrow (a \ba z) = 1$ for all $z \in A$. The case $x = \neg a$ is completely analogous. Finally, if $x = 0$, the left-hand side is $0$ for all choices of $y$ except $y = 0$, for which we get $1 \leq 1 \leftrightarrow 1 = 1$. 
				\item $((x \vee y)\boxright x) \vee ((x \vee y)\boxright y) \vee (((x \vee y)\boxright z)\leftrightarrow((x \boxright z)\wedge(y \boxright z))=1$ for all $x,y,z, \in A$. Notice that it suffices to show that at least one of the joinands is $1$ for each case.
				The identity then clearly holds if $x = 1$ or $y = 1$; consider now the case $x = a$. If $y \in \{0,a\}$ then $(x \lor y) \ba x = 1$; if $y = \neg a$ then $$((x \vee y)\boxright x) \vee ((x \vee y)\boxright y) = a \lor \neg a = 1.$$
				The case for $x = \neg a$ is symmetric. Lastly, if $x = 0$, the cases where $y \in \{1, a, \neg a\}$ are as above, and if $y = 0$ the first joinand is $0 \ba 0 = 1$. 
				\item $x \boxright (y \wedge z) = (x \boxright y)\wedge (x \boxright z)$ for all $x,y,z, \in A$. The cases $x \in \{0,1\}$ are straightforward.  Let $x = a$; if $y = 1$ we get $x \ba z = x \ba z$, and the case where $z = 1$ is analogous. Moreover, if $y = z$ the two sides of the identities are the same. If $y = a$ and $z \in \{\neg a, 0\}$, then we get: $$a \ba 0 = 0 = (a \ba a) \land 0,$$ and for the remaining cases where $y = \neg a$ the proof is once again the same. 
				The only cases left are for $x = \neg a$, which is completely analogous to $x = a$.
			\end{enumerate}
			This shows that $\alg A \in \LCA$. Now, if $\vv K$ is a variety of algebras with a lattice reduct, the deductive filters of $\alg A \in \vv K$ with respect to the logic preserving degrees of truths of $\vv K$ are the nonempty lattice filters of $\alg A$ (\cite[Proposition 2.8]{moraschiniforth}); since by Theorem \ref{th:localcomp} $\vdash_\mathbf{LV}$ is the logic preserving degrees of truth of $\LCA$, and $\alg A \in \LCA$, the $\vdash_\mathbf{LV}$-deductive filters are the nonempty lattice filters.

			The congruence filters of $\alg A$ can in turn be identified using the characterization of Proposition \ref{prop:filters}, i.e., congruence filters are nonempty lattice filters closed under $\Box$. Note that:
			\begin{align*}
				\Box 1 &= \neg 1 \ba 1 = 1;\\
				\Box a &= \neg a \ba a = 0;\\
				\Box (\neg a) &= a \ba \neg a = 0;\\
				\Box 0 &= \neg 0 \ba 0 = 0.\\
			\end{align*}
			Therefore the only lattice filters closed under $\Box$ are $\{1\}$ and $A$.
		\end{example}	
		Now, the algebra $\alg A$ in the previous example is in (the signature 	of) $\LCA$, and its deductive filters with respect to $\Cl$ are not in bijection with the $\vv K$-relative congruences for any $\vv K \sse \LCA$; this contradicts Theorem \ref{thm:isomorphismfilters}, so no subclass of $\LCA$ can be an equivalent algebraic semantics for $\Cl$. But no other class can be such by Lemma \ref{lem:basta}, thus we get the following.
		\begin{corollary}
			$\Cl$ is not algebraizable.
		\end{corollary}
		Since it can be directly checked that the algebra $\alg A$ in Example \ref{ex:algebraiz} is actually in $\CA$, by the same reasoning as above it follows that:
		\begin{corollary}
			$\Cul, \Cdl$ are not algebraizable.
		\end{corollary}	
		One can actually show that all the local variably strict conditional logics defined in Lewis's hierarchy are not algebraizable, except for $\mathbf{LVCA}$, which collapses  to classical logic.
		\begin{proposition}
			The strong and weak calculi of $\mathbf{VCA}$ coincide, and they are both algebraizable with respect to the subvariety of $\LCA$ where $x \ba y = x \to y$.
		\end{proposition}
		Now, algebraizability is preserved by axiomatic extensions, and all the logics considered by Lewis (except $\mathbf{LVCA}$) are contained in either $\mathbf{LVCSU}$, $\mathbf{LVWA}$, $\mathbf{LVTSA}$ \cite{Lewis1973}; by proving that the latter are not algebraizable, we will obtain that also all the other logics are non-algebraizable, except for $\mathbf{LVCA}$. 
		\begin{example}\label{ex:algebraiz2}
			Similarly to Example \ref{ex:algebraiz}, we construct two algebras, respectively in $\mathsf{VWA}: = \LCA + \{\tau(\mathbf{W}), \tau(\mathbf{A})\}$ and in $\mathsf{VTSA}:= \LCA + \{\tau(\mathbf{T}), \tau(\mathbf{S}), \tau(\mathbf{A})\}$, where the corresponding deductive filters are not in bijection to the congruences.
			
			In particular, consider the algebras $\alg B$ and $\alg C$, whose lattice reduct is once again given by the 4-element Boolean algebra in Figure \ref{fig:algA}, and whose operations are described respectively by the following two tables.
			
			\begin{figure}[htbp]
				\begin{center}
					\begin{tabular}{c|c|c|c|c}
						$\ba$  & $0$  & $a$ & $\neg a$ & $1$ \\\hline
						$0$  & $1$  & $1$ & $1$ & $1$ \\\hline
						$a$  & $0$  & $1$ & $0$ & $1$ \\\hline
						$\neg a$  & $0$  & $0$ & $1$ & $1$ \\\hline
						$1$  & $0$  & $0$ & $0$ & $1$
					\end{tabular}
					\caption{The $\ba$ table of $\alg B$.}\label{fig:table}
				\end{center}
			\end{figure}
			\begin{figure}[htbp]
				\begin{center}
					\begin{tabular}{c|c|c|c|c}
						$\ba$  & $0$  & $a$ & $\neg a$ & $1$ \\\hline
						$0$  & $1$  & $1$ & $1$ & $1$ \\\hline
						$a$  & $0$  & $1$ & $0$ & $1$ \\\hline
						$\neg a$  & $0$  & $0$ & $1$ & $1$ \\\hline
						$1$  & $0$  & $1$ & $0$ & $1$
					\end{tabular}
					\caption{The $\ba$ table of $\alg C$.}\label{fig:table}
				\end{center}
			\end{figure}
			The reader can check that $\alg B \in \mathsf{VWA}$, and $\alg C \in \mathsf{VTSA}$. As above, we get that the deductive filters are the lattice filters; moreover, once again in both algebras it holds that $$\Box 1 = 1, \Box a = \Box \neg a = \Box 0 = 0.$$
			The latter implies that both $\alg B$ and $\alg C$ have only two congruence filters.
		\end{example}
		\begin{theorem}
			No axiomatic extension of $\Cl$ by the axioms in $\{\alg W, \alg C, \alg N, \alg T, \alg S, \alg U, \alg A\}$ is algebraizable, except $\mathbf{LVCA}$.
		\end{theorem}
		\begin{proof}
			The algebra $\alg A$ in Example \ref{ex:algebraiz} shows that $\mathbf{LVCSU}$ is not algebraizable, since it belongs to $\mathsf{VCSU} :=\CA + \{\tau(\mathbf{U})\}$, as it can easily be checked by direct computation.
			
			Analogously, $\alg B$ and $\alg C$ in Example \ref{ex:algebraiz2} show respectively that $\mathbf{LVWA}$ and $\mathbf{LVTSA}$ are not algebraizable. Since algebraizability is preserved by axiomatic extensions, all other logics in Lewis's hierarchy that are lower in the lattice of axiomatic extensions are also not algebraizable.
		\end{proof}
		\section{Topological dualities}
		Stone representation theorem for Boolean algebras is one of the most important results in the field of algebraic logic, connecting algebraic structures to topological spaces. From the logical point of view, this deep result shows a fundamental connection between algebraic and relational models of a logical system. Stone duality has been extended and enriched in many ways over the years, starting from the classical results about Boolean algebras with operators by Jonnson and Tarski \cite{JonssonTarski51}; in this work we will proceed on these lines to connect algebraic and relational models of Lewis's logics.
		In particular, we first show a Stone-like duality for $\vv V$-algebras, in terms of Stone spaces enriched by a binary function, which are a topological expansion of Lewis's $\alpha$-models \cite{Lewis71}. We will then use the latter  to develop a duality with respect to the topological analogue of sphere models. This will allow us, at the end of the section, to obtain strong completeness of Lewis's logics with respect to all sphere models.
		
	We recall the basics of Stone duality for the sake of the reader. The underlying idea of Stone's result, which is relevant to our investigation as well, is that while finite Boolean algebras are algebras of sets, infinite Boolean algebras are correctly represented by (infinite) sets endowed with a {\em topology}, i.e. a collection of special subsets (called {\em open} sets) that essentially allow to recover the Boolean algebra. Let us be more precise; we will phrase some of the results in the language of category theory, the interested reader can check the textbook references \cite{Davey2002,mac2013} for the unexplained notions. 
	
	To each Boolean algebra $\alg B$, Stone duality associates a topological space $$\STONE(\alg B) =(Ul(\alg B), \tau)$$ where $Ul(\alg B)$ is the set of ultrafilters of $\alg B$, and $\tau$ is the Stone topology generated by the basis of clopen sets (i.e. open sets whose complement is also open) $$Cl(Ul(\alg B), \tau) = \{\stone(a): a \in B\}$$ where $$\stone(a) =\{X \in Ul(\alg B): a \in X\}.$$ 
	In general, the Boolean algebra $\alg B$ can be recovered from $\STONE(\alg B)$, since $\stone$ establishes an isomorphism between $\alg B$ and the Boolean algebra of clopen sets $\{Cl(Ul(\alg B), \tau), \cap, \cup, ^C, \emptyset, Ul(\alg B)\}$. 
	Note that if $\alg B$ is finite, its ultrafilters are in one-one correspondence with its atoms, and so $\STONE(\alg B)$ is the set of atoms of $\alg B$ with the discrete topology (where every subset is a clopen set), and $\alg B$ is isomorphic to the powerset Boolean algebra constructed from its Stone space.
	
	With respect to the maps, let $\alg B, \alg C$ be a pair of Boolean algebras, then given any homomorphism $h: \alg B \to \alg C$, Stone duality associates to $h$ a continuous map (i.e., a map such that the preimage of an open set is open) from $\STONE(\alg C)$ to $\STONE(\alg B)$ defined as $$\STONE(h)(U) = h^{-1}[U].$$
	
	$\STONE$ is a functor from the category of Boolean algebras with homomorphisms to the category of Stone spaces with continuous maps.
	Conversely, the so-called adjoint functor to $\STONE$, which we denote by $\STONE^{-1}$, associates to each a Stone space $(S, \tau)$ the Boolean algebra of its clopen subsets $$(Cl(S, \tau), \cap, \cup, ^C, \emptyset, S);$$
	in order to recover $(S, \tau)$, the following map
	\begin{equation}
		\beta(x) = \{A \in Cl(S): x \in A\}.
	\end{equation}
	yields a homeomorphism (i.e., a bijective continuous map with continous inverse) mapping $S$ to its homeomorphic copy given by the ultrafilters space of the Boolean algebra of its clopen sets.
	Moreover, given a continuous map $\varphi$ from a Stone space $(S, \tau)$ another Stone space $(S', \tau')$, its dual via Stone duality $\STONE^{-1}(\varphi)$, is defined as $$\STONE^{-1}(\varphi)(A) = \varphi^{-1}[A]$$ and it is a homomorphism from $\STONE^{-1}(S', \tau')$ to $\STONE^{-1}(S, \tau)$.
	
	Composing $\STONE$ with $\STONE^{-1}$ result in isomorphic Boolean algebras in one direction, and homeomorphic Stone spaces in the other. The dual equivalence etablished by these two functors testifies a deep connection between the algebraic and topological structures.
	
	
		\begin{notation}
			Henceforth we make a point of using the letters $a, b, c, \ldots$ for the elements of the Boolean algebras, and $X, Y, Z, \ldots$ for their ultrafilters. To allude to the correspondence established by Stone duality, we will use $x, y, z, \ldots$ for the points of a Stone space, and $A, B, C, \ldots$ for its clopen subsets. 
		\end{notation}
		
		The underlying idea of this section is to appropriately extend Stone duality; we observe that $\ba$ as a binary operator does not distribute over joins (nor meets) on both sides, therefore we cannot just apply the well-known Jonsson-Tarski duality for Boolean algebras with operators \cite{JonssonTarski}, which indeed works for operators distributing over joins on both sides and respecting $0$ (or, dually, respecting meets and $1$); however, since $\ba$ on the right does distribute over meets and maps $1$ to $1$, we can use the Jonnson-Tarski approach as an inspiration. Our treatment is also inspired by dualities for bounded distributive lattices with a ternary relation (see e.g. \cite{Celani2004}); we here prefer to use a binary function instead of a ternary relation, in order to highlight the connection with Lewis's $\alpha$-model that indeed use a so-called {\em selection function}.
		\subsection{Topological $\alpha$-models}
		We first show the duality of $\LCA$ with respect to a topological version of Lewis's $\alpha$-models. Given a Stone space $(S, \tau)$, we will make use of a binary function $f : Cl(S) \times S \to \cc P(S)$ to interpret the dualization of $\boxright$.
		We establish the following notation:
		$$f(A, x) = \{y \in S: f(A, x) = y)\}.$$
		
		\begin{definition}\label{def:top-alphamodel}
			Let us call {\em topological $\alpha$-model} a triple $(S, \tau, f)$ where $(S, \tau)$ is a Stone space, and $f : Cl(S) \times S \to \cc P(S)$ is a binary function such that, for all clopen sets $A,B$ of $S$ and $x,y \in S$:
			\begin{itemize}
				\item[($\alpha$1)] $f(A, x) \sse A$;
				\item[($\alpha$2)] $f(A, x) \sse B$ and $f(B, x) \sse A$ implies $f(A, x) = f(B, x)$; 
				\item[($\alpha$3)] Either $f(A \cup B, x) \sse A$, or $f(A \cup B, x) \sse B$, or $f(A \cup B, x) = f(A, x)  \cup f(B, x)$;
				\item[($\alpha$4)] $A \BA B = \{x \in S: f(A,x) \sse B\}$ is clopen;
				\item[($\alpha$5)] $f(A, x)$ is closed.
			\end{itemize}
		\end{definition}
		\begin{remark}
			Topological $\alpha$-models are inspired by Lewis's $\alpha$-models \cite{Lewis71}; here we essentially consider a version based on topological spaces, and then add the last two conditions. ($\alpha$4) is meant to make sure that the counterfactual $A \ba B$ is evaluated in an element of the Boolean algebra associated to the Stone space, while ($\alpha$5) is a more technical condition, necessary for the duality-theoretic machinery.
		\end{remark}
		Starting from a topological $\alpha$-model, we define a $\vv V$-algebra with the same postulates introduced by Lewis \cite{Lewis71}; of course here the underlying Boolean algebra is the Boolean algebra of clopen sets given by Stone duality.
		Given a triple $\mathfrak{S}= (S, \tau, f)$, we write $Cl(\mathfrak{S})$ for the set of clopens.
		Moreover, given any two clopen $A,B$:
		\begin{equation}
			A \BA B = \{x \in S: f(A, x) \sse B\}.
		\end{equation}
		\begin{proposition}\label{prop:COUNTinLCA}
			Let $\mathfrak{S} = (S, \tau, f)$ be a topological $\alpha$-model; $\COUNT (\mathfrak{S}) = (Cl(\mathfrak{S}), \cap, \cup, ^C, \BA, \emptyset, S)$ is in $\LCA$.
		\end{proposition}	
		\begin{proof}
			First we observe that $(Cl (\mathfrak{S}), \cap, \cup, ^C, \emptyset, S)$ is a Boolean algebra via Stone duality. We now show that the defining equations of $\LCA$ involving $\BA$ hold.	
			\begin{enumerate}
				\item Let $A \in C(\mathfrak{S})$, we show that $A \BA A=S$; indeed, $A \BA A = \{x \in S: f(A,x)\sse A\} = S$ by ($\alpha$1).
				\item We now prove that $((A \BA B) \cap (B \BA A)) \sse ((A \BA C)\leftrightarrow(B \BA C))$ for all $A,B,C \in  Cl (\mathfrak{S})$. Since $$((A \BA B) \cap (B \BA A))  = \{x \in S: f(A,x) \sse B\} \cap \{x \in S: f(B,x) \sse A\},$$ if $x \in ((A \BA B) \cap (B \BA A))$ then $f(A,x) \sse B$ and $f(B,x) \sse A$, thus by ($\alpha$2) we get that $f(A,x) = f(B,x)$.
				Now, $((A \BA C)\leftrightarrow(B \BA C)) =$ $$ (\{x \in S: f(A,x) \sse C\} \cap \{x \in S: f(B,x) \sse C\}) \cup (\{x \in S: f(A,x) \sse C\}^C \cap \{x \in S: f(B,x) \sse C\}^C).$$
				Since $f(A,x) = f(B,x)$, the claim follows from the fact that either $x \in  \{x \in S: f(A,x) \sse C\}$ or $x \in  \{x \in S: f(A,x) \sse C\}^C$.
				\item We now prove that for all $A,B,C \in  C (\mathfrak{S})$ $$((A \cup B)\BA A) \cup ((A \cup B)\BA B) \cup (((A \cup B)\BA C)\leftrightarrow((A \BA C)\cap(B \BA C))=S.$$
				It suffices to verify the right-to-left inclusion. Let $x \in S$; by ($\alpha$3) either $f(A \cup B,x) \sse A$, $f(A \cup B,x) \sse B$, or $f(A \cup B,x) = f(A,x)  \cup f(B,x)$. In the first case we get that $x \in (A \cup B)\BA A$, while in the second one it follows that $x \in (A \cup B)\BA B$. Finally, suppose $f(A \cup B,x) = f(A,x)  \cup f(B,x)$. Then 
				\begin{align*}
					&((A \cup B)\BA C)\leftrightarrow((A \BA C)\cap(B \BA C)) \\ = & \{x \in S: f(A \cup B,x) \sse C\} \leftrightarrow (\{x \in S: f(A,x) \sse C\} \cap \{x \in S: f(B,x) \sse C\}) \\ = &
					\{x \in S: f(A,x)  \cup f(B,x) \sse C\} \leftrightarrow \{x \in S: f(A,x)  \cup f(B,x) \sse C\} \\= &  S.
				\end{align*}		 
				\item It is easily seen that $A \BA (B \cap C) = (A \BA B)\cap (A \BA C)$ for all $A,B,C \in  C (\mathfrak{S})$; indeed:
				\begin{align*}
					A \BA (B \cap C) &= \{x \in S: f(A,x) \sse B \cap C\}\\ 
					&= \{x \in S: f(A,x) \sse B\} \cap \{x \in S: f(A,x) \sse C\} \\ &= (A \BA B)\cap (A \BA C).
				\end{align*}
			\end{enumerate}	
			We have shown that $\alg C (\mathfrak{S}) \in \LCA$.
		\end{proof}
		Vice versa, starting from an algebra $\alg C \in \LCA$, we obtain a topological $\alpha$-model on the Stone dual of the Boolean reduct of $\alg C$, $\STONE(\alg C)$, with the function defined by the following stipulation:
		\begin{equation}
			Y \in f_{C}(\stone(a), X) \mbox{ iff } \{c \in C: a \ba c \in X\} \sse Y
		\end{equation}
		where we recall that $\stone(a) = \{X \in Ul(\alg C): a \in X\}$, for all $a \in C, y,z \in Ul(\alg C)$. The reader shall also remember that via Stone duality all clopens of $\STONE(\alg C)$ are of the kind $\stone(a)$ for some $a \in \alg C$.
		
		The following lemma is inspired by \cite[Lemma 2.5]{Celani2004}, and it will be often used in the proofs.
		\begin{lemma}\label{lemma:RC-ba}
			Let $\alg C \in \LCA$, $a, b \in C$, then the following are equivalent:
			\begin{enumerate}
				\item $a \ba b \in X$;
				\item $Y \in f_C(\stone(a), X)$ implies $b \in Y$ ;
				\item $f_C(\stone(a), X) \sse \stone(b)$.
			\end{enumerate}
		\end{lemma}
		\begin{proof}
		(1) implies (2) follows directly from the definition of $f_C$, while the fact that (2) and (3) are equivalent follows from the definition of $\stone(b)$, $\stone(b) = \{X \in Ul(\alg C): b \in X\}$.  We prove (2) implies (1) by contrapositive. Suppose $a \ba b \notin X$, and let $F = \{c \in C: a \ba c \in X\}$; then $F$ is a Boolean filter of $\alg C$. Indeed, it is an upset because of the fact that $\ba$ is order-preserving on the right and $X$ is upwards closed, it is closed under meet since $\ba$ distributes over meets on the right and $X$ is closed under meets, and $1 \in F$ given that $a \ba 1 = 1 \in X$. Moreover, since $b \notin F$ given that $a \ba b \notin X$, by the Boolean Ultrafilter Theorem we obtain that there is a ultrafilter $Y$ of the Boolean reduct of $\alg C$ such that $F = \{c \in C: a \ba c \in X\} \sse Y$, $b \notin Y$; therefore $Y \in f_C(\stone(a), X)$ but $b \notin Y$, which concludes the proof.
		\end{proof}
		\begin{proposition}\label{prop:AC-alphamodel}
			Let $\alg C \in \LCA$; then $\mathscr{A}(\alg C) = (\STONE(\alg C), f_C)$ is a topological $\alpha$-model.
		\end{proposition}
		\begin{proof}
			Since $\STONE(\alg C)$ is a Stone space by Stone duality, we only have to prove that $f_C$ satisfies the properties in Definition \ref{def:top-alphamodel}. Note once again that by Stone duality, every clopen in $\STONE(\alg C)$ is of the kind $\stone(a)$ for some $a \in C$. 
			\begin{enumerate}
				\item[($\alpha$1)] For $a \in C, X \in Ul(\alg C)$, we have to show that $f_C(\stone(a), X) \sse \stone(a)$; suppose $Y \in f_C(\stone(a), X)$, which by definition is equivalent to $a \ba c \in X$ implies $c \in Y$, thus since by the first axiom of $\LCA$: $a \ba a = 1 \in X$, we get that $a \in Y$ and thus $Y \in \stone(a) = \{Z \in Ul(\alg C): a \in Z\}$.
				\item[($\alpha$2)] We prove that for $a, b \in C, X \in Ul(\alg C)$, $f_C(\stone(a), X) \sse \stone(b)$ and $f_C(\stone(b), X) \sse \stone(A)$ implies $f_C(\stone(a), X) = f_C(\stone(b), X)$. Suppose then $f_C(\stone(a), X) \sse \stone(b)$ and $f_C(\stone(b), X) \sse \stone(A)$, which by Lemma \ref{lemma:RC-ba} is equivalent to assuming that $a \ba b \in X$ and $b \ba a \in X$. Then also $(a \ba b) \land (b \ba a) \in X$, and by the property (\ref{C2}) of $\vv V$-algebras, we get that:
				\begin{align*}
					(a \ba b) \land (b \ba a) \leq (a \ba c) \leftrightarrow (b \ba c) &=& ((a \ba c) \to (b \ba c)) \land ((b \ba c) \to (a \ba c))\\
					&\leq& (a \ba c) \to (b \ba c), (b \ba c) \to (a \ba c) \in X.
				\end{align*}
				We want to show that $f_C(\stone(a), X) = f_C(\stone(b), X)$; we prove $f_C(\stone(a), X) \sse f_C(\stone(b), X)$, the reader shall observe that the other inclusion is proved symmetrically. Assume $Y \in f_C(\stone(a), X)$, i.e. $a \ba c \in X$ implies $c \in Y$, we show that then $Y \in f_C(\stone(b), X)$. Indeed, suppose $b \to c \in X$, then since also $(b \ba c) \to (a \ba c) \in X$ and $X$ is a Boolean filter, we get that $a \ba c \in X$; therefore $c \in Y$, which shows  $Y \in f_C(\stone(b), X)$. 
				\item[($\alpha$3)] We now show that, given $a,b \in C, X \in Ul(\alg C)$, either $f_C(\stone(a) \cup \stone(b), X) \sse \stone(a)$, or $f_C(\stone(a) \cup \stone(b), X) \sse \stone(b)$, or $f_C(\stone(a) \cup \stone(b), X) = f_C(\stone(a), X)  \cup f_C(\stone(b), X)$.
				By property (\ref{C3}) of $\vv V$-algebras, for all $a,b,c \in C$: $$((a \vee b)\boxright a) \vee ((a \vee b)\boxright b) \vee (((a \vee b)\boxright c)\leftrightarrow((a \boxright c)\wedge(b \boxright c))=1 \in X.$$
				Since $X$ is an ultrafilter, either $(a \vee b)\boxright a \in X, (a \vee b)\boxright b \in X$, or $((a \vee b)\boxright c)\leftrightarrow((a \boxright c)\wedge(b \boxright c) \in X$. Using again Lemma \ref{lemma:RC-ba}, the first two cases correspond to, respectively, $f_C(\stone(a) \cup \stone(b), X) \sse \stone(a)$ and $f_C(\stone(a) \cup \stone(b), X) \sse \stone(b)$. Consider then the remaining case where $((a \vee b)\boxright c)\leftrightarrow((a \boxright c)\wedge(b \boxright c) \in X$. Then by the definition of $\leftrightarrow$ and the fact that $\to$ distributes over $\land$ we get that also the following  elements are in $X$: $$((a \vee b)\boxright c)\to (a \boxright c),\, ((a \vee b)\boxright c)\to (b \boxright c),\, ((a \boxright c)\wedge(b \boxright c)) \to ((a \vee b)\boxright c).$$
				
				We show that $f_C(\stone(a) \cup \stone(b), X) = f_C(\stone(a), X)  \cup f_C(\stone(b), X)$.
				First, consider $Y \in f_C(\stone(a), X)  \cup f_C(\stone(b), X)$; without loss of generality, we assume $Y \in f_C(\stone(a), X)$, i.e. $a \ba c \in X$ implies $c \in Y$. We prove $Y \in f_C(\stone(a) \cup \stone(b), X)$, i.e. $(a \lor b) \ba c \in X$ implies $c \in Y$. Indeed, if $(a \lor b) \ba c \in X$, given that also $((a \vee b)\ba c)\to (a \ba c) \in X$, we get $a \ba c \in X$ and then $c \in Y$.
				
				For the other inclusion, suppose that $Y \notin f_C(\stone(a), X)  \cup f_C(\stone(b), X)$, we prove that $Y \notin f_C(\stone(a) \cup \stone(b), X)$.
				By assumption, there is $c \in C$ such that $a \ba c \in X, b \ba c \in X$ but $c \notin Y$; in more details, there must be $c', c'' \in C$ such that $a \ba c' \in X, b \ba c'' \in X$ but $c',c'' \notin Y$, then one can take $c = c' \lor c''$. Now, $(a \ba c) \land (b \ba c) \in X$, and also $((a \boxright c)\wedge(b \boxright c)) \to ((a \vee b)\boxright c) \in X$; thus $(a \vee b)\boxright c$ is in $X$ too. But $c \notin Y$, thus $Y \notin f_C(\stone(a) \cup \stone(b), X)$. 
				\item[($\alpha$4)] Let us now prove that for all $a, b \in C$, $\stone(a) \BA \stone(b)$ is clopen.
				Using the definition of $\BA$ and Lemma \ref{lemma:RC-ba} we get that $\stone$ is also a homomorphism with respect to $\ba$:
				\begin{align*}
					\stone(a) \BA \stone(b) &= \{X \in Ul(\alg C): f_C(\stone(a), X) \sse \stone(b)\}\\
					&= \{X \in Ul(\alg C): a \ba b \in X\}\\
					&= \stone(a \ba b),
				\end{align*}
				which implies that $\stone(a) \BA \stone(b)$ is clopen. 
				\item[($\alpha$5)] Lastly, we show that $f_C(\stone(a), X)$ is closed. We start by proving that $$f_C(\stone(a), X) = \bigcap \{\stone(b): f_C(\stone(a), X) \sse \stone(b)\}.$$ While the left-to-right inclusion is clear, we verify the converse. 
				Let $Y \in 	\bigcap \{\stone(b): f_C(\stone(a), X) \sse \stone(b)$, then $Y$ is such that $f_C(\stone(a), X) \sse \stone(b)$ implies $Y \in \stone(b)$ for any $b \in C$; once again using Lemma \ref{lemma:RC-ba}, the latter translates to saying that $a \ba b \in X$ implies $b \in Y$ for any $b \in C$, which is exactly the definition of $Y \in f_C(\stone(a),X)$. Thus, $f_C(\stone(a), X) = \bigcap \{\stone(b): f_C(\stone(a), X) \sse \stone(b)\}$, which means that $f_C(\stone(a), X)$ is an intersection of clopen sets, and therefore it is closed.\qedhere
			\end{enumerate}
		\end{proof}
		Proceeding towards establishing a duality, we show that the two translations we have introduced, $\mathscr{C}$ and $\mathscr{A}$, are inverses of each other.
		\begin{theorem}\label{thm:isoLCA}
			Let $\alg C \in \LCA$, then $\alg C \cong \COUNT(\ALPHA(\alg C))$.
		\end{theorem}
		\begin{proof}
			The isomorphism is given by the Stone map $\stone$ which also preserves the operation $\ba$, indeed as shown in the proof of Proposition \ref{prop:AC-alphamodel} $\stone(a \ba b) = \stone(a) \BA \stone(b)$. 
		\end{proof}
		In order to show the analogous result from the side of the topological spaces, given a Stone space $(S, \tau)$ let us call $\beta$ the homeomorphism that in Stone duality maps $S$ to its homeomorphic copy given by the ultrafilters space of the Boolean algebra of its clopen sets; precisely for $x \in S$:
		\begin{equation}
			\beta(x) = \{A \in Cl(S): x \in A\}.
		\end{equation}
		\begin{theorem}\label{thm:homeo-alpha}
			Let $\mathfrak{S} = (S, \tau, f)$ topological $\alpha$-model; then $\mathfrak{S} \cong \ALPHA(\COUNT(\mathfrak{S}))$.
		\end{theorem}
		\begin{proof}
			Via Stone duality, we only need to show that for all $A \in Cl(S), x,y \in S$: $$y \in f(A,x) \;\mbox{ iff }\; \beta(y) \in f_{\COUNT(\mathfrak{S})}(\beta[A], \beta(x)).$$
			Notice that $$\stone(A)=\{X \in Ul(\COUNT(\mathfrak{S})): A \in X\} = \bigcup_{x \in A} \{C \in Cl(\mathfrak{S}): x \in C\}= \beta[A]\;$$ thus following the definitions, $\beta(y) \in f_{\COUNT(\mathfrak{S})}(\beta[A], \beta(x))$ if and only if $A \BA C \in \beta(x)$ implies $C \in \beta(y)$ for any $C$. Equivalently, if and only if $x \in A \BA C$ implies $y \in C$ for any $C$, if and only if $f(A,x) \sse C$ implies $y \in C$ for any $C$. The latter is equivalent to $y \in f(A, x)$, since $f(A,x) = \bigcap \{C \in Cl(S): f(A,x) \sse C\}$; indeed, by axiom $(\alpha$5), $f(A,x)$ is closed, and in a Stone space clopen sets are a base, so every closed set can be written as the intersection of all the clopen sets in the base that contain it.
		\end{proof}
		We now move to studying the topological counterpart of the algebraic homomorphisms.
		We remind the reader that in Stone duality such counterpart is given by continuous maps between Stone spaces; in our setting, the continuous maps will have to satisfies some properties involving the function $f$.
		\begin{definition}\label{def:alpha-morphism}
			Let $\mathfrak{S} = (S, \tau, f), \mathfrak{S}' = (S', \tau', f')$ be topological $\alpha$-models. An {\em $\alpha$-morphism} is a map $\varphi: S \to S'$ such that:
			\begin{enumerate}
				\item\label{def:alpha-morphism1} $\varphi$ is continuous;
				\item\label{def:alpha-morphism2} $y \in f(\varphi^{-1}[A'], x)$ implies $\varphi(y) \in f'(A', \varphi(x))$;
				\item\label{def:alpha-morphism3} $y' \in f'(A', \varphi(x))$ implies that there exists $y \in S$ such that $\varphi(y) = y'$ and $y \in f(\varphi^{-1}[A'], x)$.
			\end{enumerate}
		\end{definition} 
		Given a homomorphism of $\vv V$-algebras $h: \alg C \to \alg C'$, we will show that its Stone dual $\STONE(h)$, defined as $$\STONE(h)(U) = h^{-1}[U]$$ is an $\alpha$-morphism on the corresponding topological $\alpha$-models.
		
		Moreover, we will show that given an $\alpha$-morphism $\varphi$ from $(S, \tau, f)$ to $\mathfrak{S}' = (S', \tau', f')$, its dual via Stone duality $\STONE^{-1}(\varphi)$, defined as $$\STONE^{-1}(\varphi)(A) = \varphi^{-1}[A]$$ is a homomorphism on the corresponding $\vv V$-algebras.
		\begin{proposition}\label{prop:COUNThisalphamorph}
			Given a homomorphism of $\vv V$-algebras $h: \alg C' \to \alg C$, $\STONE(h)$ is an $\alpha$-morphism from $\ALPHA(\alg C)$ to $\ALPHA(\alg C')$.
		\end{proposition}
		\begin{proof}
			$\STONE(h)$ is a continuous map from $\ALPHA(\alg C)$ to $\ALPHA(\alg C')$ by Stone duality. 
			Let us now show that for all $a' \in C', X,Y \in Ul(\alg C)$,  $$Y \in f_C(\STONE(h)^{-1}[\stone(a')], X) \,\mbox{ implies }\, \STONE(h)(Y) \in f_{C'}(\stone(a'), \STONE(h)(X)).$$
			Suppose then that  $Y \in f_C(\STONE(h)^{-1}[\stone(a')], X)$; we observe that $\STONE(h)^{-1}[\stone(a')] = \stone(h(a))$, indeed:
			\begin{align*}
				\STONE(h)^{-1}[\stone(a')] &= \{ X \in \STONE(C) : \STONE(h)(X) \in \stone(a')\}\\
				&=  \{ X \in \STONE(C) : a' \in \STONE(h)(X)\}\\
				&=  \{ X \in \STONE(C) : a' \in h^{-1}[X]\}\\
				&=  \{ X \in \STONE(C) : h(a') \in X\}\\
				&=  \stone(h(a')).
			\end{align*}
			Thus $Y \in f_C(\STONE(h)^{-1}[\stone(a')], X)$ if and only if $Y \in f_C(\stone(h(a')), X)$ if and only if, by definition of $f_C$, $h(a') \ba c \in X$ implies $c \in Y$ for any $c \in C$. Therefore, consider any $c' \in C'$, $h(c') \in C$, and we get that $h(a') \ba h(c') = h(a' \ba c') \in X$ implies $h(c') \in Y$; equivalently, for any $c' \in C'$, $a' \ba c' \in h^{-1}(X) = \STONE(h)(X)$ implies $c' \in h^{-1}(Y) = \STONE(h)(Y)$. We have showed $\STONE(h)(Y) \in f_{C'}(\stone(a'), \STONE(h)(X))$.
			
			It is left to show that for all $a' \in C', X \in Ul(\alg C), Y' \in Ul(\alg C')$, $Y'\in f_{C'}(\stone(a'), \STONE(h)(X))$ implies that there exists $Y \in Ul(C)$ such that $\STONE(h)(Y) = Y'$ and $Y \in f_C(\STONE(h)^{-1}[\stone(a')], X)$.
			Assume $Y' \in f_{C'}(\stone(a'), \STONE(h)(X))$, i.e. for all $c' \in C'$, $a' \ba c' \in \STONE(h)(X) = h^{-1}[X]$ implies $c' \in Y'$. Equivalently, for all $c' \in C'$, 
			\begin{equation}\label{eq:proof-alphamor}
				h(a' \ba c') \in X \,\mbox{ implies }\, c' \in Y'.
			\end{equation}
			We want to find $Y \in Ul(C)$ such that $\STONE(h)(Y) = Y'$ and $Y \in f_C(\STONE(h)^{-1}[\stone(a')], X)$, equivalently $Y \in f_C(\stone(h(a')), X)$, which means that for all $c \in C$, $h(a') \ba c \in X$ implies $c \in Y$.
			Consider the sets $$I = \downarrow\!\{h(c') : c' \notin Y'\} \quad F = \{c \in C: h(a') \ba c \in X\}.$$
			$F$ is easily seen to be a filter of the Boolean reduct of $\alg C$, $I$ is a Boolean ideal, and we can show that $I \cap F = \emptyset$; indeed, if $d \in I \cap F$, we get $d \leq h(c')$ for some $c' \notin Y'$, and $h(a') \ba d \in X$. But then also $h(a') \ba h(c') = h(a' \ba c') \in X$, which by (\ref{eq:proof-alphamor}) above implies $c' \in Y'$, a contradiction. Thus we can apply the Boolean Ultrafilter Theorem and obtain a ultrafilter $Y \in Ul(C)$ such that $F \sse Y$, and $Y \cap I = \emptyset$. By construction, for all $c \in C$, $h(a') \ba c \in X$ implies $c \in Y$, since $\{c \in C: h(a') \ba c \in X\} = F \sse Y$. We now prove that $\STONE(h)(Y) = Y'$ which will conclude the proof. First, if $c' \in Y'$, $\neg c' \notin Y'$ and then $h(\neg c')\in I$, thus $h(\neg c') = \neg h(c') \notin Y$, which entails that $h(c') \in Y$ given that $Y$ is a ultrafilter. This shows $Y' \sse h^{-1}[Y] = \STONE(h)(Y)$. For the other inclusion, consider $c' \notin Y'$, then $h(c') \in I$, and $h(c') \notin Y$, and thus $c' \notin h^{-1}[Y] = \STONE(h)(Y)$. The proof is complete.
		\end{proof}
		\begin{proposition}\label{prop:ALPHAphiishomo}
			Given an $\alpha$-morphism $\varphi$ from the topological $\alpha$-models $\mathfrak{S}' = (S', \tau', f')$ to $\mathfrak{S} = (S, \tau, f)$, $\STONE^{-1}(\varphi)$ is a homomorphism from $\COUNT(\mathfrak{S})$ to $\COUNT(\mathfrak{S}')$.
		\end{proposition}
		\begin{proof}
			$\STONE^{-1}(\varphi)$ is a homomorphism on the Boolean reducts by Stone duality, thus we only need to show that for all $A,B \in Cl(S)$,
			$\STONE^{-1}(\varphi)(A \BA B) = \STONE^{-1}(\varphi)(A) \BA \STONE^{-1}(\varphi)(B)$.
			
			Now, $$\STONE^{-1}(\varphi)(A \BA B) = \varphi^{-1}[\{x \in S: f(A,x) \sse B\}] = \{ x' \in S': f(A, \varphi(x')) \sse B\}$$
			and
			$$\STONE^{-1}(\varphi)(A) \BA \STONE^{-1}(\varphi)(B) = \varphi^{-1}[A] \BA \varphi^{-1}[B] = \{x' \in S': f'(\varphi^{-1}[A], x') \sse  \varphi^{-1}[B]\}.$$
			Let $x' \in \STONE^{-1}(\varphi)(A \BA B)$, i.e. $f(A, \varphi(x')) \sse B$; we show that $x' \in \STONE^{-1}(\varphi)(A) \BA \STONE^{-1}(\varphi)(B)$, i.e. $f'(\varphi^{-1}[A], x') \sse  \varphi^{-1}[B]$. Indeed, if $y' \in f'(\varphi^{-1}[A], x')$, by Definition \ref{def:alpha-morphism}(\ref{def:alpha-morphism2}) we get $\varphi(y') \in f(A, \varphi(x'))$ and then $\varphi(y') \in B$, that is, $y' \in \varphi^{-1}[B]$.
			
			Finally, let $x' \in \STONE^{-1}(\varphi)(A) \BA \STONE^{-1}(\varphi)(B)$, i.e. $f'(\varphi^{-1}[A], x') \sse  \varphi^{-1}[B]$, we show that $x' \in \STONE^{-1}(\varphi)(A \BA B)$, i.e. $f(A, \varphi(x')) \sse B$. If $y \in f(A, \varphi(x'))$, by Definition \ref{def:alpha-morphism}(\ref{def:alpha-morphism3}) we get that there exists $y' \in S'$ such that $\varphi(y') = y$ and $y' \in f'(\varphi^{-1}[A], x')$; then $y' \in  \varphi^{-1}[B]$, that is, $\varphi(y') = y \in B$. This shows the second inclusion and concludes the proof.
		\end{proof}
		The previous results show that Stone duality extends to a duality between the algebraic category $\LCA$ (with objects the algebras in $\LCA$ and morphism the homomorphisms) and a category whose objectst are topological $\alpha$-models and whose morphisms are $\alpha$-morphisms.
		Let us denote the latter category by $\CATA$.
		
		We define a map $\ALPHA: \LCA \to \CATA$ in the following way, for any $\alg C, \alg C' \in \LCA, h: \alg C \to \alg C'$:
		\begin{align*}
			\ALPHA(\alg C) &= (\STONE({\rm Bool}(\alg C)), f_C);	 \\
			\ALPHA(h) &= \STONE(h);
		\end{align*}
		where ${\rm Bool}(\alg C)$ is the Boolean reduct of $\alg C$.
		Moreover, we define a map $\COUNT: \CATA \to \LCA$ such that, for any $\mathfrak{S} = (S, \tau, f), \mathfrak{S}' = (S', \tau', f') \in \CATA, \varphi: \mathfrak{S} \to  \mathfrak{S}'$:
		\begin{align*}
			\COUNT(\mathfrak{S}) &= (\STONE^{-1}((S, \tau)), \BA);\\	 
			\ALPHA(\varphi) &= \STONE(\varphi).
		\end{align*}
		Both $\ALPHA$ and $\COUNT$ can be easily shown to be well-defined functors since they properly extend the Stone duality functors $\STONE$ and $\STONE^{-1}$ (on objects via Propositions   \ref{prop:COUNTinLCA} and \ref{prop:AC-alphamodel}, on morphisms via Propositions \ref{prop:COUNThisalphamorph} and \ref{prop:ALPHAphiishomo} ). Moreover, the fact that they are adjoint functors and they establish a duality between $\LCA$ and $\CATA$ follows from Stone duality and the isomorphism of Theorem \ref{thm:isoLCA}, and the homeomorphism of Theorem \ref{thm:homeo-alpha}. Therefore:
		\begin{theorem}\label{thm:dual}
			The functors $\ALPHA$ and $\COUNT$ establish a dual equivalence between the algebraic category of $\vv V$-algebras and the category of topological $\alpha$-models $\CATA$.
		\end{theorem}
		With the following lemma we can extend the duality to the most relevant subvarieties of $\LCA$. We follow again Lewis ideas in \cite{Lewis71}.
		\begin{definition}
			Let us call a topological $\alpha$-model $\mathfrak{S} = (S, \tau, f)$ a: 
			\begin{enumerate}
				\item topological $\alpha_1$-model if it is such that for each $A \in Cl(S)$ and $x \in A$, $f(A,x) = \{x\}$;
				\item topological $\alpha_2$-model if it is a topological $\alpha_1$-model such that for each $A \in Cl(S)$ and $x \in S$, $f(A,x)$ contains at most one element.
			\end{enumerate}
			
		\end{definition}
		\begin{lemma}
			Let $\mathfrak{S} = (S, \tau, f)$ be a topological $\alpha$-model. Then:
			\begin{enumerate}
				\item if $\mathfrak{S}$ is a topological $\alpha_1$-model, $\COUNT(\mathfrak{S}) \in \LC$;
				\item if $\mathfrak{S}$ is a topological $\alpha_2$-model, $\COUNT(\mathfrak{S}) \in \CA$.
			\end{enumerate}
		\end{lemma}
		\begin{proof}
			For the first point, let us start by recalling that $\COUNT(\mathfrak{S}) \in \LC$ if and only if
			$$A \cap B \sse A \BA B \sse A^C \cup B.$$
			Moreover, by definition $A \BA B = \{x \in S: f(A,x) \sse B\}$. Let now $x \in A \cap B$;  then in particular $x \in A$ and then if $\mathfrak{S}$ is a topological $\alpha_1$-model, $f(A, x) =\{x\} \sse B$, which shows the first inclusion. For the second one, let $x \in A \BA B$; if $x \notin A^C$, i.e. $x \in A$, once again $f(A, x) =\{x\}$, and then $x \in B$ given that $x \in A \BA B$.
			
			For the second point, we recall that $\COUNT(\mathfrak{S}) \in \CA$ if and only if:
			$$(A \BA B) \cup (A \BA B^C) = S$$
			where $A \BA B = \{x \in S: f(A, x) \sse B\}, A \BA B^C = \{x \in S: f(A, x) \sse B^C\}$. 
			Consider any $x \in S$; if $\mathfrak{S}$ is a topological $\alpha_2$-model, $f(A,x)$ contains at most one element. Thus either $f(A,x) = \emptyset$, and then $x$ is in both $A \BA B$ and $A \BA B^C$, or $f(A,x)= \{z\}$; thus since either $z \in B$ or $z \in B^C$, also $x \in A \BA B$ or $x \in A \BA B^C$.
		\end{proof}
		Conversely:
		\begin{lemma}
			Let $\alg C \in \LCA$. Then:
			\begin{enumerate}
				\item If $\alg C \in \LC$ then $\ALPHA(\alg C)$ is a topological $\alpha_1$-model;
				\item If $\alg C \in \CA$ then $\ALPHA(\alg C)$ is a topological $\alpha_2$-model.
			\end{enumerate}
		\end{lemma}
		\begin{proof}
			Suppose $\alg C \in \LC$, we show that for each $a \in \alg C$ and $X \in Ul(\alg C)$ such that $a \in X$, $f(\stone(a),X) = \{X\}$.
			First, we show $X \in f(\stone(a),X)$; let $a \ba b \in X$, then $a \ba b \leq a \to b \in X$. Now, $a, a \to b \in X$ implies $b \in X$, which proves $X \in f(\stone(a),X)$.
			
			We now show that if $Y \in f(\stone(a),X)$ holds, then $Y=X$. Indeed, let $b \in X$; then since also $a \in X$, $a \land b \in X$, and $a \land b \leq a \ba b \in X$, thus from $Y \in f(\stone(a),X)$ we get $b \in Y$, which shows $X \sse Y$. For the other inclusion, let $b \in Y$. If $b \notin X$, then $\neg b \in X$ and then $a \land \neg b \in X$, $a \land \neg b \leq a \ba \neg b \in X$ and this would give $\neg b \in Y$, a contradiction. Thus $b \in X$ and we have showed $X = Y$.
			
			For the second point, suppose $\alg C \in \CA$. We prove that $f_C(\stone(a), X)$ has at most one element for all $a \in C, X \in Ul(\alg C)$. Suppose that both $Y \in f_C(\stone(a),X)$ and $Z \in f_C(\stone(a),X)$ hold; this means that $a \ba b \in X$ implies $b \in Y \cap Z$. Let $c \in Y$; since $\alg C \in \CA$, $$(a \ba c) \lor (a \ba \neg c) = 1 \in X,$$ and hence $a \ba c \in X$ or $a \ba \neg c \in X$, since $X$ is a Boolean ultrafilter. If $a \ba c \in X$, then $c \in Y \cap Z \sse Z$;  $a \ba \neg c \in X$ instead yields a contradiction, indeed it entails that $\neg c \in Y\cap Z \sse Y$. Thus $c \in Z$ and $Y \sse Z$; the other inclusion can be proved analogously, which implies that $Y=Z$ and therefore $f_C(\stone(a), X)$ has at most one element.
		\end{proof}
		Hence, the functors $\ALPHA$ and $\COUNT$ restrict to the full subcategories of topological $\alpha_1$-models (and $\alpha_2$-model) and $\LC$ ($\CA$), yielding a duality of the subcategories.
		\begin{theorem}
			The algebraic category of $\LC$-algebras (resp., $\CA$-algebras) is dually equivalent to the category of topological $\alpha_1$-models (resp., $\alpha_2$-models) with $\alpha$-morphisms.
		\end{theorem}
		\subsection{Topological spheres}
		We will now show that topological $\alpha$-models are categorically equivalent to a category of {\em topological spheres}. More precisely, we describe the category whose objects are {\em equivalent} in the sense of Lewis \cite{Lewis71} to $\alpha$-models; that is, they are defined on the same set, and they define the same algebra: the Boolean algebras have the same domain, and the operations $\ba$ obtained respectively by the spheres and the function $f$ coincide. 
		
		Lewis himself notices that in order to obtain equivalent models, one needs to assume the {\em limit assumption}: given any point $x$ and its set of spheres, for each formula $\varphi$ that intersects the spheres at $x$, there is a smallest sphere that does so. 
		We will indeed use this assumption to define the category of topological spheres. 
		
		\begin{definition}\label{def:topsphere}
			A \em{topological sphere} is a triple $(S, \tau, \sigma)$, where $(S, \tau)$ is a Stone space and $\sigma: S \to \mathcal{P}(\mathcal{P}(S))$ is such that:
			\begin{enumerate}
				\item[(S1)] for all $x \in S$, $\sigma(x)$ is nested, i.e. for all $U,V \in \sigma(x)$, $U \sse V$ or $V \sse U$;
				\item[(S2)] for all $A,B \in Cl(S)$, $A \BAS B = \{x \in S: A \cap \displaystyle\bigcap_{U \in \sigma(x), U \cap A \neq \emptyset} U \sse B\}$ is clopen;
				\item[(S3)]\label{def:topsphere-limit} for all $x \in S$, $A \in Cl(S)$, if $A \cap \bigcup\sigma(x) \neq \emptyset$, there exists a smallest $U \in \sigma(x), U \cap A \neq \emptyset $;
				\item[(S4)] for all $x \in S$, $A \in Cl(S)$, let $\Sigma(A,x) =  \displaystyle\bigcap_{U \in \sigma(x), U \cap A \neq \emptyset} U$; then $A \cap \Sigma(A,x)$ is closed.
			\end{enumerate}
		\end{definition}
		\begin{remark}
		It is clear that the definition of topological spheres is inspired by Lewis's sphere models. In the next subsection we will see how to obtain sphere models as in Definition \ref{def:lewissphere} from a topological sphere. Moreover, importantly, we will discuss condition (S3), which is Lewis's limit assumption on sphere models.
		\end{remark}
		
		We will now define the maps, inspired by Lewis's work (see both \cite{Lewis1973,Lewis71}), that obtain a topological sphere from an $\alpha$-model and viceversa.
		Consider a topological sphere $(S, \tau, \sigma)$; we associate the triple $\Alpha((S, \tau, \sigma)) = (S, \tau, f_{\sigma})$ where for all $A \in Cl(S)$, $x \in S$:
		\begin{equation}
			f_\sigma(A, x) = A \cap \Sigma(A,x) = A \cap \bigcap\{U \in \sigma(x): U \cap A \neq \emptyset\}.
		\end{equation}
		The following is a straightforward consequence of the definition.
		\begin{lemma}\label{lemma:R_sigma}
			For all $A \in Cl(S)$, $x \in S$, either $\Sigma(A,x) = \emptyset$ or it is the smallest smallest $U \in \sigma(x), U \cap A \neq \emptyset $. Moreover:
			\begin{enumerate}
				\item $f_\sigma(A, x)  = \emptyset$ iff there is no $U \in \sigma(x)$ such that $U \cap A \neq \emptyset$;
				\item $f_\sigma(A, x) \neq \emptyset$ iff $A \cap \bigcup\sigma(x) \neq \emptyset$.
			\end{enumerate}
		\end{lemma}
		\begin{proposition}
			Consider a topological sphere $(S, \tau, \sigma)$; $\Alpha((S, \tau, \sigma)) = (S, \tau, f_\sigma)$ is a topological $\alpha$-model.
		\end{proposition}
		\begin{proof}
			We have to show that $f_\sigma$ satisfies the conditions of Definition \ref{def:top-alphamodel}.
			\begin{enumerate}
				\item[($\alpha$1)] $f_\sigma(A, x) \sse A$ for all $A \in Cl(S)$ by definition.
				\item[($\alpha$2)] We proceed to prove that $f_\sigma(A, x) \sse B$ and $f_\sigma(B, x) \sse A$ implies $f_\sigma(A, x) = f_\sigma(B, x)$; let us then assume that $f_\sigma(A, x) \sse B$ and $f_\sigma(B, x) \sse A$. It follows by the definition and Lemma \ref{lemma:R_sigma} that either both sets are empty (and there is nothing else to prove), or they are both nonempty. Suppose we are in the latter case; thus $f_\sigma(A, x) = A \cap \Sigma(A, x) \sse B$ and $f_\sigma(B, x) = B \cap \Sigma(B, x) \sse A$.
				Since $\emptyset \neq f_\sigma(B, x) = B \cap \Sigma(B, x) \sse A$, we get that $\Sigma(B, x) \cap A \neq \emptyset$; but
				$\Sigma(A, x)$ is the least element in $\sigma(x)$ that intersects with $A$, thus $\Sigma(A, x) \sse \Sigma(B, x)$. Analogously, $\Sigma(B, x) \sse \Sigma(A, x)$ and therefore $\Sigma(A, x) = \Sigma(B, x)$. Thus $f_\sigma(A, x) = A \cap  \Sigma(A, x) = A \cap \Sigma(A, x) \cap B = A \cap \Sigma(B, x) \cap B = \Sigma(B, x) \cap B = f_\sigma(B,x).$
				\item[($\alpha$3)] Given $A, B \in Cl(S), x \in S$ we show that either $f_\sigma(A \cup B, x) \sse A$, or $f_\sigma(A \cup B, x) \sse B$, or $f_\sigma(A \cup B, x) = f_\sigma(A, x)  \cup f_\sigma(B, x)$; in other words:
				$$(A \cup B) \cap \Sigma(A \cup B,x) \sse A \mbox{ or } (A \cup B) \cap \Sigma(A \cup B,x) \sse B \mbox{ or }  (A \cup B) \cap \Sigma(A \cup B,x) =  (A \cap \Sigma(A,x)) \cup  (B \cap \Sigma(B,x)).$$
				Since $\sigma(x)$ is nested, we can assume without loss of generality that $\Sigma(A,x) \sse \Sigma(B,x)$. We consider the following cases:
				\begin{itemize}
					\item $\emptyset \neq \Sigma(A,x) = \Sigma(B,x)$; then necessarily $\Sigma(A \cup B,x) = {\rm min} \{U \in \sigma(x): U \cap (A \cup B) \neq 0\} = \Sigma(A,x) = \Sigma(B,x)$. Thus:
					$$(A \cup B) \cap \Sigma(A \cup B,x) = (A \cap \Sigma(A \cup B,x)) \cup (B \cap \Sigma(A \cup B,x)) = (A \cap \Sigma(A,x)) \cup  (B \cap \Sigma(B,x)).$$
					\item $\emptyset \neq \Sigma(A,x) \subsetneq \Sigma(B,x)$; it follows that $\Sigma(A,x) \cap B = \emptyset$ and then $\emptyset \neq \Sigma(A \cup B,x) = \Sigma(A,x)$. Hence:
					$$(A \cup B) \cap \Sigma(A \cup B,x) = (A \cap \Sigma(A \cup B,x)) \cup (B \cap \Sigma(A \cup B,x)) = A \cap \Sigma(A, x) \sse A.$$
					\item  $\emptyset = \Sigma(A,x) \subsetneq \Sigma(B,x)$; then $\Sigma(A \cup B,x) = \Sigma(B,x)$, and 
					$$(A \cup B) \cap \Sigma(A \cup B,x) = (A \cap \Sigma(A \cup B,x)) \cup (B \cap \Sigma(A \cup B,x)) = B \cap \Sigma(B, x) \sse B.$$
					\item $\emptyset = \Sigma(A,x) = \Sigma(B,x)$; in this case also $\Sigma(A \cup B,x) = \emptyset$, and then trivially the last condition holds, i.e.  $(A \cup B) \cap \Sigma(A \cup B,x) =  (A \cap \Sigma(A,x)) \cup  (B \cap \Sigma(B,x))$.
				\end{itemize}
				\item[($\alpha$4)] The fact that $A \BA B = \{x \in S: f_\sigma(A, x) \sse B\}$ is clopen follows from the fact that $A \BAS B$ is clopen and $\BA$ and $\BAS$ coincide, indeed:
				$$A \BA B = \{x \in S: f_\sigma(A,x) \sse B\} = \{x \in S: A \cap \Sigma(A,x) \sse B\} = A \BAS B.$$
				\item[($\alpha$5)] Lastly, $f_\sigma(A, x) = A \cap \Sigma(A, x)$ is closed  by (S5) in the definition.
			\end{enumerate}
		\end{proof}
		For the converse translation, consider a topological $\alpha$-model $(S, \tau, f)$; we associate the triple $\Sphere((S, \tau, f)) = (S, \tau, \sigma_f)$ such that for all $x \in S$:
		\begin{equation}
			\sigma_f(x) = \left\{\bigcup_{A <_x B} f(A,x) : B \in Cl(S)\right\}
		\end{equation} 
		where 
		\begin{equation}
			A <_x B \;\;\mbox{ iff }\;\; f(B,x) = \emptyset \;\mbox{ or }\;\; \emptyset \neq f(A,x) \sse f(A \cup B,x).
		\end{equation}
		Our definition extends to topological spheres the stipulations in \cite{Lewis71}, where Lewis constructs an $\alpha$-model from a sphere with the limit assumption.
		The following properties are stated without proof in \cite{Lewis71} (and elsewhere in the literature); we include a proof here for the sake of the interested reader.
		\begin{lemma}\label{lemma:preorder}
			Consider a topological $\alpha$-model $(S, \tau, f)$; for all $A, B \in Cl(S)$, $x \in S$:
			\begin{enumerate}
				\item If $A \sse B$ and $A \cap f(B,x) \neq \emptyset$, then $f(A,x) = A \cap f(B,x)$;
				\item if $A \sse B$ and $f(B,x) = \emptyset$, then $f(A,x) =\emptyset$;
			\end{enumerate}
		\end{lemma}
		\begin{proof}
			For the first property, suppose $A \sse B$ and $A \cap f(B,x) \neq \emptyset$. Notice that $$B = B \cap (A \cup A^C) = (B \cap A) \cup (B \cap A^C).$$
			In particular, since $A \sse B$, $B = A \cup (B \cap A^C)$ and then $f(B,x) = f(A \cup (B \cap A^C), x)$. Thus, by ($\alpha3$) we get the following cases:
			\begin{itemize}
				\item $f(B,x) \sse A$, and therefore, since by ($\alpha1$) $f(A,x) \sse A \sse B$, by ($\alpha2$) $f(B,x) = f(A,x)$, which implies that $f(A,x) = A \cap f(B,x)$;
				\item $f(B,x) \sse B \cap A^C$, which yields a contradiction since $A \cap f(B,x) \neq \emptyset$, and thus it is never the case;
				\item $f(B,x) = f(A,x) \cup f(B \cap A^C, x)$, which (together with ($\alpha1$)) directly implies that $f(A,x) \sse A \cap f(B,x)$. For the other inclusion, let $y \in A \cap f(B,x)$; then $y \in f(A,x) \cup f(B \cap A^C, x)$. Since $f(B \cap A^C, x) \sse B \cap A^C \sse A^C$, and $y \in A$, it follows that $y \notin f(B \cap A^C, x)$. Therefore, $y \in f(A,x)$ and we have shown that $f(A,x) = A \cap f(B,x)$.
			\end{itemize}
			Hence, in all the possible cases $f(A,x) = A \cap f(B,x)$.
			
			The second property is easier to show. Indeed, suppose $A \sse B$ and $f(B,x) = \emptyset$. Then $f(B,x) = \emptyset \sse A$ and (using ($\alpha1$)) $f(A,x) \sse A \sse B$; by ($\alpha2$) it follows that $f(A,x) = f(B,x) =\emptyset$.
		\end{proof}
		Using the last lemma, we can show the following facts about the relation $<_x$. The following also appears in \cite{Lewis71} with some hints to a proof, we offer here a complete proof which seems to miss from the literature.
		\begin{lemma}\label{lemma:totalpreorder}
			Consider a topological $\alpha$-model $(S, \tau, f)$; for all $A, B \in Cl(S)$, $x \in S$:
			\begin{enumerate}
				\item $<_x$ is a total preorder on $Cl(S)$;
				\item $A <_x B$ implies $B \cap f(A,x) \sse f(B,x)$;
				\item $A \cap f(B,x) \neq \emptyset$ implies $A <_x B$. 
			\end{enumerate}
		\end{lemma}
		\begin{proof}
			(1) We start by verifying that $<_x$ is a total preorder on $Cl(S)$ for each $x \in S$; reflexivity is straightforward from the definition, so let us prove transitivity. We assume that $A <_x B$ and $B <_x C$, i.e.:
			$$f(B,x) = \emptyset \;\mbox{ or }\;\; \emptyset \neq f(A,x) \sse f(A \cup B,x),$$
			$$f(C,x) = \emptyset \;\mbox{ or }\;\; \emptyset \neq f(B,x) \sse f(B \cup C,x),$$
			and we show that $A <_x C$, that is, $f(C,x) = \emptyset \;\mbox{ or }\;\; \emptyset \neq f(A,x) \sse f(A \cup C, x)$. Suppose then that $f(C,x) \neq \emptyset$; thus since $B <_x C$, $\emptyset \neq f(B,x) \sse f(B \cup C,x)$, and since $A <_x B$ also $\emptyset \neq f(A,x) \sse f(A \cup B,x)$. We need to show that $f(A,x) \sse f(A \cup C, x)$. 
			First, we observe that $(A \cup B) \cap f(A \cup (B \cup C), x) \neq \emptyset$; this follows from $(\alpha 3)$, which yields that either $f(A \cup (B \cup C), x)$ is contained in $A$, and thus intersects with $A \cup B$, or in the two remaining cases it contains $f(B \cup C,x) \supseteq f(B,x)$ (the last inclusion holds by hypothesis), and thus again it intersects with $A \cup B$. Therefore (using Lemma \ref{lemma:preorder} for the last inclusion):
			$$\emptyset \neq f(A, x) \sse f(A \cup B,x) \sse f(A \cup B \cup C, x),$$
			which implies that $A \cap f(A \cup B \cup C, x) \neq \emptyset$. Then also $(A \cup C) \cap f(A \cup B \cup C, x) \neq \emptyset$, and again by Lemma \ref{lemma:preorder} $$f(A \cup C, x) = (A \cup C) \cap f(A \cup B \cup C, x) = (A \cap f(A \cup B \cup C, x)) \cup (C \cap f(A \cup B \cup C, x)).$$
			Thus, since  $A \cap f(A \cup B \cup C, x) \neq \emptyset$, we get that $A \cap f(A \cup C, x)\neq \emptyset$, and by a last application of Lemma \ref{lemma:preorder} we obtain that  $f(A,x) \sse f(A \cup C, x)$ which shows that $A <_x C$; hence $<_x$ is transitive for all $x \in S$. 
			
			For the first point, it is left to show that $<_x$ is total. Consider $f(A \cup B,x)$; if $f(A \cup B,x) = \emptyset$, by Lemma \ref{lemma:preorder} also $f(A,x) = f(B,x) = \emptyset$, and thus $A <_x B$ and $B <_x A$. If $f(A \cup B,x) \neq \emptyset$, then by $(\alpha3)$ either $f(A \cup B,x) \sse A$, which using $(\alpha 2)$ yields $f(A,x) = f(A \cup B,x)$ and then  $A <_x B$, or $f(A \cup B,x) \sse B$, which analogously yields $B <_x A$, or $f(A \cup B,x) = f(A,x) \cup f(B,x)$, and therefore $A <_x B$ and $B <_x A$. In any case, $A$ and $B$ are comparable with respect to $<_x$, which is then a total preorder.
			
			(2) Suppose that $A <_x B$, i.e. $$f(B,x) = \emptyset \;\mbox{ or }\;\; \emptyset \neq f(A,x) \sse f(A \cup B,x).$$
			Assume first that $f(B,x) = \emptyset$, and let us consider two cases: whether $f(A \cup B,x)$ is empty or not.
			If $f(A \cup B,x) = \emptyset$, by Lemma \ref{lemma:preorder} also $f(A,x) = \emptyset$ and then $B \cap f(A,x) =\emptyset = f(B,x)$. If $f(A \cup B,x) \neq \emptyset$, from the fact that $f(B,x) = \emptyset$, using $(\alpha 3)$ and $(\alpha2)$ we get that $f(A \cup B,x) = f(A,x)$ and (necessarily by Lemma \ref{lemma:preorder}(1)) $B \cap f(A \cup B,x) = \emptyset$. Thus $B \cap f(A,x) = \emptyset = f(B,x)$.
			
			Assume now that $\emptyset \neq f(A,x) \sse f(A \cup B,x)$. Then either $B \cap f(A,x) = \emptyset \sse f(B,x)$ or $$\emptyset \neq B \cap f(A,x) \sse B \cap f(A \cup B,x)$$
			and then by Lemma \ref{lemma:preorder}(1) $f(B,x) = B \cap f(A \cup B,x)$ thus 
			$$B \cap f(A,x) \sse B \cap f(A \cup B,x)  = f(B,x).$$
			This completes the proof of (2).
			
			(3) Finally, assume $A \cap f(B,x) \neq \emptyset$, which implies that in particular $f(B,x) \neq \emptyset$. Consider $f(A \cup B,x)$, then by $(\alpha 3)$ we have three cases:
			\begin{itemize}
				\item $f(A \cup B,x) \sse A$; then by $(\alpha 2)$, $f(A \cup B,x) = f(A,x)$ which is necessarily nonempty, since otherwise by Lemma \ref{lemma:preorder}(2) also $f(B,x)$ would be empty, which is a contradiction. Thus we have that $A <_x B$.
				\item $f(A \cup B,x) \sse B$; then again by $(\alpha 2)$, $f(A \cup B,x) = f(B,x)$ and then $A \cap f(A \cup B,x) = A \cap f(B,x) \neq \emptyset$. By Lemma \ref{lemma:preorder}(1) this yields that $\emptyset \neq f(A, x) \sse f(A \cup B, x)$ and then $A <_x B$.
				\item $f(A \cup B,x) = f(A,x) \cup f(B,x)$; thus $f(A, x) \sse f(A \cup B, x)$ and $f(A, x)$ is nonempty (otherwise we would be in the case above), and once again $A <_x B$.
			\end{itemize} 
			We have shown that in every case $A <_x B$, which concludes this proof.
		\end{proof}
		We are ready to show that we can construct a topological sphere from a topological $\alpha$-model.
		\begin{proposition}\label{prop:topsphere}
			Consider a topological $\alpha$-model $(S, \tau, f)$; $\Sphere((S, \tau, f)) = (S, \tau, \sigma_f)$ is a topological sphere.
		\end{proposition}
		\begin{proof}
			We show that $\sigma_f$ satisfies the properties of Definition \ref{def:topsphere}.
			First, notice that $\sigma_f$, defined as $$\sigma_f(x) = \left\{\bigcup_{A <_x B} f(A,x) : B \in Cl(S)\right\}
			$$ is indeed a map from $S$ to $\mathcal{P}(\mathcal{P}(S))$.
			Moreover:
			\begin{enumerate}
				\item for all $x \in S$, $\sigma_f(x)$ is nested, since $<_x$ is total by Lemma \ref{lemma:totalpreorder}. 
				\item We now show (S3), i.e. that for all $x \in S$, $A \in Cl(S)$, if $A \cap \bigcup\sigma_f(x) \neq \emptyset$, there exists a smallest $U \in \sigma_f(x), U \cap A \neq \emptyset$. In particular, we show that such smallest set is exactly $U_A = \displaystyle\bigcup_{C <_x A} f(C,x)$.
				
				First, notice that \begin{equation}\label{eq:inproofalpha-sph}
					f(A,x) = \emptyset \;\; \mbox{ iff }\;\;\not\exists U\in \sigma_f(x): A \cap U \neq \emptyset.
				\end{equation}
				Indeed, the right-to-left direction is obvious since $f(A, x) \in U_A$ and $f(A,x) \sse A$. We show the left-to-right one by contrapositive. Note that if there is $U \in \sigma_f(x)$ such that $A \cap U \neq \emptyset$, then there is a $B$ such that $A \cap f(B,x) \neq \emptyset$ and so $A <_x B$ by Lemma \ref{lemma:totalpreorder}; thus $f(A, x) \neq \emptyset$ by the definition of $<_x$.
				
				Now, assume $A \cap \bigcup\sigma_f(x) \neq \emptyset$, then  $\emptyset \neq f(A, x) \sse A$ (the last inclusion by $(\alpha1)$), and since $A <_x A$ ($<_x$ is reflexive by Lemma \ref{lemma:totalpreorder}),  then $A \cap U_A \neq \emptyset$. 
				
				We now prove that $U_A$ is the smallest $U \in \sigma_f(x), U \cap A \neq \emptyset$; indeed, suppose $A \cap U_B = A \cap \displaystyle\bigcup_{C <_x B} f(C,x) \neq 0$. Then $A \cap f(C,x) \neq \emptyset$ for some $C <_x B$, and then by Lemma \ref{lemma:totalpreorder} $A <_x C$; by the transitivity of $<_x$, it holds that $A <_x B$, and thus $U_A \sse U_B$.
				
				\item We now proceed to prove (S4), that is, for all $x \in S$, $A \in Cl(S)$, $A \cap \Sigma(A,x) =A \cap \displaystyle\bigcap_{U \in \sigma_f(x), U \cap A \neq \emptyset} U$ is closed. By the previous point, either $\Sigma(A,x) = \emptyset$ or $\Sigma(A,x) = U_A = \displaystyle\bigcup_{C <_x A} f(C,x)$. In the first case, the empty set is a closed set in any topological space. In the second case, note that by Lemma \ref{lemma:totalpreorder} for all $C <_x A$, $A \cap f(C,x) \sse f(A,x) = A \cap f(A,x)$. Thus since $A <_x A$:
				$$A \cap \displaystyle\bigcup_{C <_x A} f(C,x) = \displaystyle\bigcup_{C <_x A} A \cap f(C,x) = A \cap f(A,x)$$
				and the latter set is closed, since $f(A,x)$ is closed by definition of a topological $\alpha$-model, $A$ is closed, and closed sets are closed under finite intersection.
				\item Finally, the fact that $A \BAS B$ is clopen for all $A,B \in Cl(S)$  follows from the fact that $A\BAS B = A \BA B$, given that the latter is clopen by definition of a topological $\alpha$-model. Let us verify the equality; it is helpful to first recall the definitions:
				\begin{align*}
					A \BAS B &= \{x \in S: A \cap \displaystyle\bigcap_{U \in \sigma_f(x), U \cap A \neq \emptyset} U \sse B\};\\
					A \BA B &= \{x \in S: f(A,x) \sse B\}.
				\end{align*}
				Notice that by (\ref{eq:inproofalpha-sph}) above, $f(A, x) = \emptyset$ if and only if there is no $U \in \sigma_f(x), U \cap A \neq \emptyset$. If $x$ is such, then clearly $A \BAS B = A \BA B$ for any $B \in Cl(S)$. 
				
				Consider now $x \in S$ such that $f(A, x) \neq \emptyset$, and then $A \BAS B = \{x \in S: A \cap U_A \sse B\}$. Hence it is clear that if $x \in A \BAS B$, also $x \in A \BA B$ since $f(A,x) \sse U_A$; conversely, suppose that $x \in A \BA B$, i.e. $f(A,x) \sse B$. Then $$A \cap U_A = A \cap  \displaystyle\bigcup_{C <_x A} f(C,x) \sse f(A,x) \sse B,$$
				since if $C <_x A$ then $A \cap f(C,x) \sse f(A, x)$ by Lemma \ref{lemma:totalpreorder}.\qedhere
			\end{enumerate}
		\end{proof}
		We now show that if we start from a topological $\alpha$-model $(S, \tau, f)$, and we apply $\Sphere$ and then $\Alpha$, we obtain {\em exactly} the same $\alpha$-model. 
		\begin{proposition}\label{prop:alpha-essentsurj}
			Consider a topological $\alpha$-model $(S, \tau, f)$; then $\Alpha(\Sphere((S, \tau, f))) = (S, \tau, f)$, i.e. $f_{\sigma_{f}} = f$.
		\end{proposition}
		\begin{proof}
			We show that given any $A \in Cl(S), x,y \in S$, $y \in f(A,x)$ if and only if $y \in f_{\sigma_{f}}(A,x)$. Recall that $$f_{\sigma_{f}}(A,x) = A \cap \displaystyle\bigcap_{U \in \sigma_f(X), U \cap A \neq \emptyset} U.$$ 
			First, suppose $y \in f(A,x)$ holds; then since $f(A, x) \neq \emptyset$, as shown in the proof of Proposition \ref{prop:topsphere} (point (2) of the proof), $f_{\sigma_{f}}(A,x) = A \cap \bigcup_{C <_x A}f(C,x)$. We have that $y \in A$, $y \in f(A, x)$, and since $A <_x A$, we get that $y \in f_{\sigma_{f}}(A,x)$.
			
			Vice versa, suppose $y \in f_{\sigma_{f}}(A,x)$; this means that $f_{\sigma_{f}}(A,x)$ is nonempty, and then again $f_{\sigma_{f}}(A,x) = A \cap \bigcup_{C <_x A}f(C,x)$. Moreover, since if $C <_x A$ then $A \cap f(C,x) \sse f(A, x)$ by Lemma \ref{lemma:totalpreorder}, $y \in f_{\sigma_{f}}(A,x)$ implies $y \in f(A, x)$, which concludes the proof.
		\end{proof}
		We now move to describe the morphisms among topological spheres.
		\begin{definition}
			Consider two topological spheres $(S, \tau, \sigma), (S', \tau', \sigma')$; then a map $\varphi: (S, \tau, \sigma) \to (S', \tau', \sigma')$ is a {\em sphere morphism} if:
			\begin{enumerate}
				\item $\varphi$ is continuous;
				\item for all $y \in S, A' \in Cl(S')$ such that $\varphi^{-1}[A'] \cap \bigcup\sigma(x) \neq \emptyset$,  $y \in \varphi^{-1}(A') \cap \Sigma(\varphi^{-1}[A'], x)$ implies $\varphi(y) \in \Sigma(A', \varphi(x))$;
				\item for all $y' \in S', A' \in Cl(S')$ such that $A' \cap \bigcup\sigma'(\varphi(x)) \neq \emptyset$, $y' \in A' \cap \Sigma(A', \varphi(x))$ implies that there is $y \in S$ such that $\varphi(y) = y'$ and $y \in \Sigma(\varphi^{-1}[A'], x)$.
			\end{enumerate}
		\end{definition}
		\begin{proposition}\label{prop:Alpha-func}
			Consider two topological spheres $(S, \tau, \sigma), (S', \tau', \sigma')$, and $\varphi: (S, \tau, \sigma) \to (S', \tau', \sigma')$ a sphere morphism; then $\varphi$ is an $\alpha$-morphism from $\Alpha((S, \tau, \sigma))$ to $\Alpha((S', \tau', \sigma'))$.
		\end{proposition}
		\begin{proof}
			We only need to show the two following properties:
			\begin{enumerate}
				\item $y \in f_\sigma(\varphi^{-1}[A'], x)$ implies $\varphi(y) \in f_{\sigma'}(A', \varphi(x))$;
				\item $y' \in f_{\sigma'}(A', \varphi(x))$ implies that there exists $y \in S$ such that $\varphi(y) = y'$ and $y \in f_\sigma(\varphi^{-1}[A'], x)$.
			\end{enumerate}
			For (1), assume that $y \in f_\sigma(\varphi^{-1}[A'], x)$. 
			Since $$f_\sigma(\varphi^{-1}[A'], x) = \varphi^{-1}[A'] \cap \Sigma(\varphi^{-1}[A'], x)$$
			we get that $\varphi^{-1}[A'] \cap \bigcup\sigma(x) \neq \emptyset$, and $y \in \varphi^{-1}(A') \cap \Sigma(\varphi^{-1}[A'], x)$; by definition of a sphere morphism, this implies $\varphi(y) \in \Sigma(A', \varphi(x))$, and of course $\varphi(y) \in A'$. Thus, since $f_{\sigma'}(A', \varphi(x)) = A' \cap \Sigma(A', \varphi(x))$, we get $\varphi(y) \in f_{\sigma'}(A', \varphi(x))$. 
			
			For (2) we proceed similarly; suppose  $y' \in f_{\sigma'}(A', \varphi(x))$. Thus $A' \cap \bigcup\sigma'(\varphi(x)) \neq \emptyset$, and by definition of $f_{\sigma'}(A', \varphi(x))$, $y' \in A' \cap \Sigma(A', \varphi(x))$. Hence, by definition of a sphere morphism, there is $y \in S$ such that $\varphi(y) = y'$ and $y \in \Sigma(\varphi^{-1}[A'], x)$; this yields exactly $y \in f_\sigma(\varphi^{-1}[A'], x)$.
		\end{proof}
		Conversely:
		\begin{proposition}
			Consider two topological $\alpha$-models $(S, \tau, f), (S', \tau', f')$, and let $\varphi: (S, \tau, f) \to (S', \tau', f')$ be an $\alpha$-morphism; then $\varphi$ is a sphere morphism from $\Sphere((S, \tau, f))$ to $\Sphere((S', \tau', f'))$.
		\end{proposition}
		\begin{proof}
			We need to show the following:
			\begin{enumerate}
				\item for all $y \in S, A' \in Cl(S')$ such that $\varphi^{-1}[A'] \cap \bigcup\sigma_f(x) \neq \emptyset$,  $y \in \varphi^{-1}(A') \cap \Sigma_f(\varphi^{-1}[A'], x)$ implies $\varphi(y) \in \Sigma_{f'}(A', \varphi(x))$;
				\item for all $y' \in S', A' \in Cl(S')$ such that $A' \cap \bigcup\sigma_{f'}(\varphi(x)) \neq \emptyset$, $y' \in A' \cap \Sigma_f(A', \varphi(x))$ implies that there is $y \in S$ such that $\varphi(y) = y'$ and $y \in \Sigma_f(\varphi^{-1}[A'], x)$.
			\end{enumerate}
			Let us show (1). Assume that $\varphi^{-1}[A'] \cap \bigcup\sigma_f(x) \neq \emptyset$,  $y \in \varphi^{-1}(A') \cap \Sigma_f(\varphi^{-1}[A'], x)$. Then we have that:
			$$\Sigma_f(\varphi^{-1}[A'], x) = \bigcup_{C<_x \varphi^{-1}[A']}f(C,x).$$
			Then $y \in \varphi^{-1}(A') \cap f(C,x)$ for some $C<_x \varphi^{-1}[A']$; thus $y \in f(\varphi^{-1}(A'), x)$ by Lemma \ref{lemma:totalpreorder}. Moreover by definition of an $\alpha$-morphism, we get that $\varphi(y) \in f'(A', \varphi(x))$. Thus $$\varphi(y) \in A' \cap \bigcup_{C'<_x A'} f(C', \varphi(x)) = \Sigma_{f'}(A', \varphi(x)).$$
			
			We verify (2) in a similar fashion. Suppose $A' \cap \bigcup\sigma_{f'}(\varphi(x)) \neq \emptyset$, $y' \in A' \cap \Sigma_f(A', \varphi(x))$. Then 
			$$y' \in A' \cap \bigcup_{C'<_{\varphi(x)} A'}f(C', \varphi(x)) \sse f(A', \varphi(x)),$$
			where the last inclusion holds by Lemma \ref{lemma:totalpreorder};
			thus by definition of an $\alpha$-morphism we get that there exists $y$ such that $\varphi(y) = y'$ and $y \in f(\varphi^{-1}[A'],x)$. Thus $$y \in \varphi^{-1}[A'] \cap f(\varphi^{-1}[A'],x)=\Sigma_f(\varphi^{-1}[A'], x)$$
			which concludes the proof.
		\end{proof}
		Let us consider the category of topological spheres with sphere morphisms, which we denote by $\TOPSPH$. 
		Consider the map $\Alpha: \TOPSPH \to \CATA$ defined as: 
		\begin{align*}
			\Alpha((S, \tau, \sigma)) =& (S, \tau, f_\sigma);	 \\
			\Alpha(\varphi) =& \varphi.
		\end{align*}
		$\Alpha$ is easily seen to be a functor between the two categories given Proposition \ref{prop:Alpha-func}.
		\begin{theorem}\label{thm:cateq}
			The functor $\Alpha$ establishes a categorical equivalence between the categories of topological $\alpha$-models with $\alpha$-morphisms and topological spheres with sphere morphisms.
		\end{theorem}
		\begin{proof}
			In order to show that $\Alpha$ establishes a categorical equivalence, it suffices to say that it is full, faithful and essentially surjective (see for instance \cite{mac2013}).
			Fullness and faithfulness (i.e., injectivity and surjectivity on morphisms among corresponding pairs of objects) are obvious since $\Alpha$ is the identity on morphisms. The fact that $\Alpha$ is essentially surjective follows from Proposition \ref{prop:alpha-essentsurj}.
		\end{proof}
		Finally, we show that the categorical equivalence extends to the full subcategories of topological $\alpha_1$-models and  $\alpha_2$-models with respect to the following.
		\begin{definition}
			{\color{white}.}
			\begin{enumerate}
				\item 	Let us call {\em centered} a topological sphere $(S, \tau, \sigma)$ such that $\{x\} \in \sigma(x)$ for all $x \in S$.
				\item Let us call {\em Stalnakerian} a  centered sphere $(S, \tau, \sigma)$ such that if $A \cap \bigcup \sigma(x) \neq \emptyset$, then there exists $U \in \sigma(x)$ and $y \in S$ such that $A \cap U = \{y\}$.
			\end{enumerate}
		\end{definition}
		\begin{lemma}
			Let $(S, \tau, \sigma)$ be a topological sphere. Then 
			\begin{enumerate}
				\item If $(S, \tau, \sigma)$ is a centered topological sphere, $\Alpha((S, \tau, \sigma))$ is a topological $\alpha_1$ model.
				\item If $(S, \tau, \sigma)$ is a Stalnakerian topological sphere, $\Alpha((S, \tau, \sigma))$ is a topological $\alpha_2$ model.
			\end{enumerate}
		\end{lemma}
		\begin{proof}
			{\color{white}.}
			\begin{enumerate}
				\item Suppose $(S, \tau, \sigma)$ is centered, i.e. $\{x\} \in \sigma(x)$ for all $x \in S$.
				We need to show that if $x \in A$, $f_\sigma(A,x) = \{x\}$. Now, $f_\sigma(A,x) = A \cap \Sigma(A, x)$. Recall that $\Sigma(A, x)$ is the intersection of all $U \in \sigma(x)$ that intersects with $A$, and $\Sigma(A, x)$ belongs to $\sigma(x)$ if it is nonempty. Thus if $x \in A$, necessarily $\Sigma(A, x) = \{x\}$, which implies the claim.  
				\item Suppose now $(S, \tau, \sigma)$ is Stalnakerian.
				We want to show that for each $A \in Cl(S)$ and $x \in S$, $f_\sigma(A,x)$ contains at most one element. 
				By hypothesis, if $A \cap \bigcup \sigma(x) \neq \emptyset$, then there exists $U \in \sigma(x)$ and $y \in S$ such that $A \cap U = \{y\}$.
				Since $f_\sigma(A,x) = A \cap \Sigma(A, x)$, either $f_\sigma(A,x) = \emptyset$, or otherwise $A \cap \Sigma(A, x) = \{y\}$.
			\end{enumerate}
		\end{proof}
		Conversely:
		\begin{lemma}
			Let $(S, \tau, f)$ be a topological $\alpha$-model. Then 
			\begin{enumerate}
				\item If $(S, \tau, f)$ is a $\alpha_1$ model, $\Sphere((S, \tau, f))$ is a centered topological sphere.
				\item If $(S, \tau, f)$ is a $\alpha_2$ model, $\Sphere((S, \tau, f))$ is a Stalnakerian topological sphere.
			\end{enumerate}
		\end{lemma}
		\begin{proof}
			{\color{white}.}
			\begin{enumerate}
				\item First, assume that for each $A \in Cl(S)$ and $x \in A$, $f(A,x) = \{x\}$; we show that $\{x\} \in \sigma_f(x)$ for all $x \in S$. Recall that $$\sigma_f(x) = \{\bigcup_{B <_x A} f(B,x) : A \in Cl(S)\}.$$
				Consider a clopen $A$ such that $x \in A$ (since clopen sets are a base of the Stone space, and singletons are closed, there is surely at least one). Then $f(A,x) = \{x\}$; moreover, if $B <_x A$, by definition of the preorder $f(B,x) \sse f(A \cup B,x)$, and since $x \in A\sse A \cup B$, by hypothesis $f(A \cup B,x) = \{x\}$ and thus $f(B,x) \sse \{x\}$. Hence $$\{x\} = \bigcup_{B <_x A} f(B,x) \in  \sigma_f(x).$$
				
				\item Now, assume that for each $A \in Cl(S)$ and $x \in S$, $f(A,x)$ contains at most one element. We prove that if $A \cap \bigcup \sigma_f(x) \neq \emptyset$, then there exists $U \in \sigma_f(x)$ and $y \in S$ such that $A \cap U = \{y\}$.
				Assume $A \cap \bigcup \sigma_f(x) \neq \emptyset$; then $\Sigma_f(A, x)$ is nonempty, and as shown in the proof of Proposition \ref{prop:topsphere}:
				$$\Sigma_f(A, x) = A \cap \bigcup_{B <_x A}f(B,x).$$
				In particular, by Lemma \ref{lemma:totalpreorder}  $f(A, x)$ is nonempty, and thus it contains exactly one element, say $f(A, x) = \{y\}$. Moreover, if $B <_x A$, $A \cap f(B,x) \sse f(A, x) = \{y\}$.
				Therefore, $A \cap \Sigma_f(A, x) = \{y\}$ and $\Sigma_f(A, x)$ is the desired element of $\sigma_f(x)$. \qedhere
			\end{enumerate}
		\end{proof}
	As a consequence of the duality theorem (Theorem \ref{thm:dual}) and the categorical equivalence (Theorem \ref{thm:cateq}) and their restrictions to the most relevant subclasses, we can compose the functors and obtain the following connection between $\vv V$-algebras and topological spheres.
		\begin{theorem}
			The category of topological spheres is dually equivalent to the algebraic category of $\vv V$-algebras, as testified by the composition of the functors $\Alpha$ and $\COUNT$. Moreover, restricting the composed functors to the full subcategories of, respectively, 
			centered topological spheres and Stalnakerian topological spheres, we get a dual equivalence with respect to $\LC$-algebras and $\CA$-algebras.\end{theorem}
		\subsection{Strong completeness and the limit assumption}
		In this last subsection, let us circle back to Lewis's sphere models. 
		We will prove that strong completeness holds with respect to the topological sphere models we have introduced, that correspond to a subclass of Lewis's sphere models. This actually implies the strong completeness with respect to all sphere models, both for the global and the local semantics.
		
		Consider a topological sphere model $\cc M = (S, \tau, \sigma)$, and let $v : Fm_{\cc L} \to \COUNT(\Alpha(\cc M))$ be any assignment of the variables in $Var$ to the Boolean algebra of clopen subsets of $S$.
		Then of course $\cc M_L=(S, \sigma, v)$ is a sphere model satisfying the limit assumption in the sense of Lewis (Definition \ref{def:lewissphere}). 
		
		\begin{definition}
			We define the class $\mathsf{LTS}$ of {\em Lewis topological spheres} to be the class of all the sphere models $\cc M_L$ obtained from a topological sphere model $\cc M$ as above.
		\end{definition}
		Now, $\mathbf{GV}$ is strongly complete with respect to $\LCA$ (Corollary \ref{cor:algebraizability}). As a consequence of the results in this section, every algebra in $\LCA$ arises from a topological sphere model and vice versa (see Theorem \ref{thm:isoLCA}, Proposition \ref{prop:alpha-essentsurj}), which allows us to obtain the following:
		\begin{theorem}[Global strong completeness]\label{thm:gstrong}
			For any $\Gamma \cup \{\varphi\} \sse Fm_{\cc L}$, the following hold:
			$$\Gamma \vdash_\mathbf{GV} \varphi \Leftrightarrow \Gamma \models_{\mathsf{LTS},g} \varphi$$\end{theorem}
		\begin{proof}
			Soundness is easily checked: since $\vdash_\mathbf{GV}$ is finitary, $\Gamma \vdash_\mathbf{GV} \varphi$ if and only if there is a finite subset $\Gamma_0 \sse \Gamma$ such that $\Gamma_0 \vdash_\mathbf{GV} \varphi$; by the finite strong completeness (Theorem \ref{th: globalompleteness}), this is equivalent to $\Gamma_0 \models_g \varphi$; in particular the consequence holds in the subclass of sphere models $\mathsf{LTS}$, and we can conclude that $\Gamma \models_{\mathsf{LTS},g} \varphi$.
			
			We prove completeness by contrapositive; assume that $\Gamma \not\vdash_\mathbf{GV} \varphi$. Then by the algebraizability result (Theorem \ref{thm:Cgalg}), we get that $\tau(\Gamma) \not\models_{\LCA} \tau(\varphi)$, i.e. there is $\alg A \in \LCA$, and an assignment from $Fm_{\cc L}$ to $\alg A$ such that $\alg A, h \models \tau(\Gamma)$, but $\alg A, h \not\models \tau(\varphi)$. Precisely, $h(\gamma) = 1$ for all $\gamma \in \Gamma$, but $h(\varphi) \neq 1$. Consider now $\cc M = \Sphere(\ALPHA(\alg A))$, and take $\cc M_L$ where the evaluation $v$ is induced by $h$, i.e. $v(\varphi) = \stone(h(\varphi))$. Then we get that $\cc M_L \Vdash \Gamma$, since $\varphi$ is mapped to the whole universe, but $\varphi$ is not so $\cc M_L \not\Vdash \varphi$; the latter yields $\Gamma \not \models_{\mathsf{LTS},g} \varphi$ and concludes the proof.
		\end{proof}
		Let us now show the similar result for the weaker calculus.
		\begin{theorem}[Local strong completeness]\label{thm:lstrong}
			For any $\Gamma \cup \{\varphi\} \sse Fm_{\cc L}$, the following hold:
			$$\Gamma \vdash_\mathbf{LV} \varphi \Leftrightarrow \Gamma \models_{\mathsf{LTS},l} \varphi$$\end{theorem}
		\begin{proof}
			Soundness can be proved analogously to the case of the global semantics. 
			Suppose now that $\Gamma \not\vdash_\mathbf{LV} \varphi$, we prove that $\Gamma \not\models_{\mathsf{LTS},l} \varphi$; since $\vdash_\mathbf{LV}$ is the logic preserving the degrees of truth of $\LCA$, it follows that there exists $\alg A \in \LCA$, an assignment $h$ of the variables $Var$ to $\alg A$ and an element $a \in A$ such that \begin{equation}\label{eq:proofcomp}
				a \leq h(\gamma) \mbox{ for all } \gamma \in \Gamma,\; a \not\leq h(\varphi).
			\end{equation}
			Without loss of generality, we can consider $\alg A$ to be the algebra obtained from a topological sphere model: precisely, every algebra in $\LCA$ arises as an algebra defined on a topological $\alpha$-model by Theorem \ref{thm:isoLCA}, and every topological $\alpha$-model arises by a topological sphere model on the same Stone space and defining the same algebra by Proposition \ref{prop:alpha-essentsurj}. One can then see $\alg A$ as the algebra defined on a  Lewis topological sphere $\cc M = (S, \sigma, v)$, where all elements are (clopen) subsets of $S$. Thus (\ref{eq:proofcomp}) above can be rephrased as
			$$
			w \in a \mbox{ implies } w \in h(\gamma) \mbox{ for all } \gamma \in \Gamma.
			$$
			and 
			$$
			\mbox{ there exists } w' \in S \mbox{ such that } w' \in a, w' \notin \varphi.
			$$
			Thus, there exists $w$ such that $w\Vdash\Gamma$, but $w\not\Vdash\varphi$; hence, $\Gamma \not\models_{\mathsf{LTS},l} \varphi$, which completes the proof.
		\end{proof}
		
		We conclude with a final remark about the models we are considering. Note that all Lewis topological spheres satisfy the limit assumption (we indeed included it in the axiomatization of topological sphere models). While Lewis argues that one ``has no right to assume that there always are smallest antecedent-permitting sphere and, within it, a set of closest-antecedent worlds'' \cite[p.20]{Lewis1973}, his logics actually are perfectly captured by the logical behavior of the sphere models that do satisfy the limit assumption, as witnessed by the above completeness theorem. 
		
		Nonetheless, of course the above results imply that Lewis's logics are strongly sound and complete with respect to the larger class of all sphere models, including those that do not satisfy the limit assumption:
		
		\begin{theorem}
			For any $\Gamma \cup \{\varphi\} \subseteq Fm_\mathcal{L}$, the following hold:
			\begin{equation}
				\Gamma \vdash_\mathbf{GV} \varphi \Leftrightarrow \Gamma \vDash_g \varphi 
			\end{equation}
			\begin{equation}
				\Gamma \vdash_\mathbf{LV} \varphi \Leftrightarrow \Gamma \vDash_l \varphi
			\end{equation}
			
			\noindent
			
		\end{theorem}
		
		\begin{proof}
			Since both logics are finitary, soundness can be proved as in Theorem \ref{th: globalompleteness} and Theorem \ref{thm:lewiscomplete}. 
		Completeness follows directly from Theorem \ref{thm:gstrong} and Theorem \ref{thm:lstrong}, since if a deduction fails in $\vdash_\mathbf{GV}$ (or $\vdash_\mathbf{LV}$), its failure is already testified in sphere models in $\mathsf{LTS}$, and therefore also in the larger class of all sphere models. 
		\end{proof}
		And once again, the above results can be extended to every logic in the family of Lewis variably strict conditional logics and the corresponding class of sphere models.
			\begin{theorem}
		Let $\Sigma$ be a subset of the axioms $\{\alg W, \alg C, \alg N, \alg T, \alg S, \alg U, \alg A\}$.
			For all subsets $\Gamma \cup \{\varphi\}\subseteq For_{\cc L}$,
			\begin{equation}
			 	\Gamma \vdash_{\mathbf{GV\Sigma}} \varphi \Leftrightarrow \Gamma \models_{\vv K_{\Sigma}, g} \varphi. 
			\end{equation}
			\begin{equation}
				\Gamma \vdash_{\mathbf{LV\Sigma}} \varphi \Leftrightarrow \Gamma \models_{\vv K_{\Sigma}, l} \varphi.
			\end{equation}
			
		\end{theorem}
		
		The strong completeness results with respect to topological spheres suggests that Lewis sphere models possess an excess of expressiveness relatively to Lewis's logic. This richness, as evidenced by the limit assumption, arises from the models' ability to encode information that exceeds the expressive boundaries of the considered propositional logics. In particular, in order to isolate the subclass of sphere models that satisfy the limit assumption, one would likely require a more expressive language, such as an infinitary one. Let us elaborate more on this.
		
		As Fine \cite{Fine2012}, Starr \cite{sep-counterfactuals}, and Lewis himself \cite{Lewis71} have observed, the absence of the limit assumption implies the failure of a very intuitive infinitary rule. More precisely, let us consider now an infinitary language over a denumerable set of propositional variables $\{p_0, p_1, p_2, \dots\}$ which admits infinitary conjunction and disjunction; for instance  $p_0 \wedge p_1 \wedge \cdots$ (where the $p_i$'s are propositional variables) would be a formula in this language. Furthermore, let us assume that this language can be interpreted within sphere models in the expected way; namely, given a sphere model $(W, \mathcal{S}, v)$:
	$$v(\bigwedge\limits_{i \in I} \varphi_i)=\bigcap\limits_{i \in I} v(\varphi_i),\;\; v(\bigvee\limits_{i \in I} \varphi_i)=\bigcup\limits_{i \in I} v(\varphi_i)$$
	and the remaining connectives are interpreted as usual. Then, with $\Gamma \cup \{\psi\}$ being a set of formulas in this infinitary language, the following meta-rule turns out to be not valid in the class of sphere models without the limit assumption (while it does hold in those that do satisfy it): 
	\begin{center}if $\Gamma \models \psi$, then $\{\varphi \boxright \gamma : \gamma \in \Gamma\} \models \varphi \boxright \psi$.\end{center} 
	The failure is testified by Lewis's counterexample in \cite{Lewis71}, which we replicate here for the sake of the reader. 
	\begin{example}
		Consider the sphere model $\cc M^* = (W, \mathcal{S}, v)$ where
		\begin{itemize}
			\item $W$ is the set of real numbers;
			
			\item for all $x \in W$, $\mathcal{S}(x)=\{[x, i] : x \leq i \}$ where $\leq$ is the natural order over the real numbers and $[x, i]$ is a usual closed interval in the Euclidean topology;
			
			\item $v$ maps the propositional variables to the following open intervals: $v(p_0)=(0, \infty)$, and for each $p_i$ with $1 \leq i$, $v(p_i)=(-\infty, 2^{-i})$, so for instance $v(p_1)=(-\infty, \frac{1}{2})$.
		\end{itemize}
			Now, it is easy to realize that given any set of formulas $\Gamma$, $\Gamma \models \bigwedge \Gamma$ is a valid \emph{local} logical consequence with respect to all the sphere models, indeed, for every sphere model $\mathcal{M}$, for every world $w$ in $\mathcal{M}$, if $w \Vdash \gamma$ for all $\gamma \in \Gamma$, then clearly $w \Vdash \bigwedge \Gamma$. However, consider the sphere model $\cc M^*$. It is clearly the case that $\{p_1, p_2, \dots\} \models \bigwedge\limits_{i=1}^\infty p_i$. Moreover, the counterfactual formulas $p_0 \boxright p_1, p_0\boxright p_2, p_0\boxright p_3, \dots$ can be checked to hold at $w = 0$. However, $p_0 \boxright \bigwedge\limits_{i=1}^\infty p_i$ is not true at $0$ since there are certainly points in the sphere associated to 0 where $p_0$ is true, but there is no point in the model where both $p_0$ and $\bigwedge\limits_{i=1}^\infty p_i$ are true. The failure of the meta-rule depends on the fact that there is no minimal sphere in $\mathcal{S}(0)$ where $p_0$ holds, and this entails the failure of the limit assumption. 
	\end{example}
		Let us conclude the section by stressing once again that if one considers the finitary propositional language,  the logic corresponding to all sphere models satisfying the limit assumption is exactly the same as the logic defined over the broader class of all sphere models. In this sense, the models without limit assumption are redundant, as testified also by the duality that does not ``see'' them.
		
		In fact, the sphere models satisfying the limit assumption occupy a privileged position. They provide a complete and faithful representation of variably strict conditional logic, as their topological counterparts are dually equivalent to the equivalent algebraic semantics of the logic. As far as our study of the logics and their properties is concerned, there exists no independent justification for considering models without the limit assumption. On the other hand, there exists a compelling mathematical rationale for investigating variably strict conditional logics through the lens of models satisfying the limit assumption: these models offer a comprehensive and faithful representation of the logic, as exemplified by the duality results.
		
		\section{Conclusions}
		The primary objective of this paper was to provide a comprehensive logico-algebraic analysis of Lewis variably strict conditional logics, a subject that has been notably lacking in the literature. To achieve this goal, we embarked on a multifaceted investigation that encompassed both syntactic and semantic aspects of these logics.
		
	First, we sought to refine the traditional logical framework employed to represent Lewis's logic by introducing new, simpler axiomatizations of his logical systems, based on the counterfactual arrow as their main primitive connective.
	
	 Next, we delved into an investigation of the \emph{logics} associated to those systems by distinguishing between a weaker and stronger consequence relations definable from those systems. We then explored the corresponding semantical consequences associated to those logics, by analyzing their properties and connections in a way that parallels the analysis of the global and local consequence of modal logic. Throughout this investigation, we introduced model-theoretic tools that, to the best of our knowledge, were not employed before to analyze these logics (such as the generated sub-model for spheres). These powerful tools proved instrumental in establishing several intriguing logical results, including two versions of the deduction theorem and a representation of one logical consequence in terms of the other.
		
		Subsequently, we shifted our focus to the algebraic investigation of the logics introduced at the beginning; we defined a new class of Boolean algebras equipped with a counterfactual operator, demonstrating their structural properties, such as the characterization of their deductive filters. The algebraic investigation culminated in the algebraizability results, highlighting the connection between the algebras and the corresponding logics. In more detail, the algebras we introduced served as equivalent algebraic semantics for the global version of the logic. The local version, while not strictly algebraizable, was characterized as the logic that preserves degrees of truth over these algebraic structures.
		
		Our exploration culminated in a topological study of the logics, revisiting Lewis's standard sphere semantics through a novel perspective. In particular, we consider a class of Stone spaces equipped with a new operation, and establish their dual equivalence with the algebraic  category of the algebras introduced earlier. These duality results unveiled a deeper mathematical understanding of Lewis variably strict conditional logics through the lens of topological spaces. Moreover, they shed new light on the status of the limit assumption and its intricate relationship with the logic.
		
		In conclusion, we hope that this work has provided a extensive logico-algebraic treatment of Lewis variably strict conditional logics. Our efforts have clarified several ambiguities in the literature surrounding these logics, explicitly defined and refined their properties and theorems, and introduced a novel general algebraic framework for their technical analysis. By doing so, this research aims to prove once more the fruitful synergy between a classical theme in formal philosophy and the advancements in abstract algebraic logic.

		\section*{Acknowledgements}
		
		Giuliano Rosella acknowledges that this work has been funded by a grant from the Programme Johannes Amos Comenius under the Ministry of Education, Youth and Sports of the Czech Republic, CZ.02.01.01/00/23\_025/0008711.

		
		\nocite{*}
		
		\printbibliography

	\end{document}